\documentclass[11pt]{article}

\usepackage{epsfig}
\usepackage{graphicx}
\usepackage{color}
\usepackage{psfrag}

\usepackage{mathptmx,amsmath,amssymb}

\newcommand{\DS}{\LARGE \renewcommand{\baselinestretch}{1.67} \normalsize}

\newtheorem{theorem}{Theorem}[section]

\title{Algorithms for Deforming and Contracting Simply Connected Discrete Closed Manifolds (III)}

\author{
Li Chen \\
Department of Computer Science and Information Technology\\
University of the District of Columbia\\
Email: lchen@udc.edu\\
}

\begin{document}
\DS

\maketitle

\abstract
{In a recent paper, {\it Algorithms for Deforming and Contracting Simply Connected Discrete Closed Manifolds (II)}, we discussed two algorithms for deforming and contracting a simply connected discrete closed manifold into a discrete sphere. The first algorithm was a continuation of work that began in {\it Algorithms for Deforming and Contracting Simply Connected Discrete Closed Manifolds (I)}, the second algorithm contained a more direct treatment of contraction for discrete manifolds. In this paper, we clarify that we can use this same method on standard piecewise linear (PL) complexes on the triangulation of general smooth manifolds. Our discussion is based on triangulation techniques invented by Cairns, Whitehead, and Whitney more than half of a century ago.
In this paper, we use PL or simplicial complexes to replace certain concepts of discrete manifolds in previous papers. Note that some details in the original papers related to discrete manifolds may also need to be slightly modified for this purpose.  In this paper, we use the algorithmic procedure (a contraction process) to prove the following theorem: For a finite triangulation of a simply-connected closed and orientable 3-manifold $M$ in Euclidean space, if a (simply-connected closed) 2-cycle which was made by 2-cells of this triangulation separates $M$ into two connected components, then each of the components will also be simply-connected. In addition, we can algorithmically make $M$ to be homeomorphic to a 3-sphere. The relationship of the theorem to a generalized Jordan separation problem, the general Jordan-Schoenflies theorem, and other important problems are also discussed.
This revision adds more figures and explanations as well as two more special cases. We will post the detailed algorithm/procedure for practical triangulations in the following papers.
}



\section{Introduction}

According to the Whitney embedding theorem and the Nash embedding theorem,
any smooth real $m$-dimensional manifold can be smoothly embedded in $R^{2m}$ and any $m$-manifold with a Riemannian metric (Riemannian manifold) can be isometrically embedded in
an $n$-Euclidean space $E^n$ where $n\le c\cdot m^3$.
Therefore, we only need to consider smooth manifold $\cal M$ in $E^n$.

In 1940, based on Cairns's work, Whitehead first proved
the following theorem~\cite {Whit40}:  Every smooth manifold admits an (essentially unique) compatible piecewise linear (PL) structure.

A PL structure is a triple $(K,{\cal M},f)$ where $K$ is a PL-complex and $f: K\to \cal M$ is piecewise
differentiable.  In addition, Whitehead proved (roughly speaking)
that there is a sequence of homeomorphic mappings $f_i$ on the finite complex $K_i$ that are arbitrarily good approximations of $f$. ($K_i$ is a collection of finite cells, each of these cells is a subdivision of some cells in $K_{i-1}$.)

In 1957, Whitney designed a simpler procedure to get a triangulation~\cite{Whitney57}. In 1960, Cairns found an even simpler method for compact smooth manifolds\cite{Cairns35,Cairns}: A compact smooth manifold admits a finite triangulation, i.e. there is a finite complex that is diffeomorphic to the original manifold. A finite complex means that a complex has a finite number of simplices
or (PL) cells.
Therefore, for simplicity, a compact smooth manifold can be triangulated by a finite complex~\cite{Cairns35} (while preserving the diffeomorphism between them).

The more general research is related to {\it combinatorial triangulation conjecture} that is
stated as follows:

\hskip 1 in {\it Every compact topological manifold can be triangulated by a PL manifold.}

Rado in 1925 proved that this conjecture holds in two dimensional cases. Moise followed by proving the conjecture in three dimensions.
However, the conjecture is false in dimensions $\ge 4$ ~\cite{LeeBook}.

In this paper, we usually always assume that a manifold {\cal M} is embedded in a Euclidean space. And, it is smooth, closed, compact, orientable, and simply connected. A discrete manifold is
slightly different than a PL manifold. A discrete manifold is a computer representation of a PL manifold. The cells in discrete manifold is in the abstract form. For example, a 2-cell in
PL manifold has a centroid that can be calculated. A 2-cell in a discrete manifold might not have such a centroid, in order to add a centroid, we need to embed a 2-cell into a Euclidean or other space to
get one.  Therefore, the discrete manifold is more general than the PL manifold.  Every PL manifold has its discrete manifold representation. A PL manifold can be naturally equivalent to its discrete manifold
representation when embedding to Euclidean space. We will not distinguished them in this paper. In fact, we only need to use $k$-simplexes instead of $k$-cells in most part of this paper.

In \cite{Che17I},  we discussed two algorithms for
 deforming and contracting a simply connected discrete
closed manifold into a discrete sphere.  The first algorithm was a continuation of work initiated  in
 \cite{Che15I}.  The second algorithm was a more direct treatment of contractions for discrete manifolds.

This paper focuses on the second algorithm and its relationship to smooth manifolds. We clarify that we can use the same method for standard PL complexes from the triangulation of general smooth manifolds.
Our discussion is based on
triangulation techniques invented by Cairns,  Whitehead, and Whitney more than half of a century ago.
This paper focuses on the 3-manifold, the most important case.  However, we also discuss the more general $k$-manifold.

In this paper{\footnote {This is the third revision of the paper. In this revision, we mainly add some figures to show the proof in
Section 4 of this paper. We also made a separated short paper to explain a theorem used as Part I of the proof of Theorem 4.3. \cite{Che18SummerTop}.
In this paper, we use PL or simplicial complexes to replace certain concepts of discrete manifolds in previous papers.
Note that some details in the original papers related to discrete manifolds may also need to be slightly modified for this purpose.
}
we use the algorithmic procedure to prove the following theorem: {\it For a finite triangulation of a simply-connected closed 3-manifold $M$ in Euclidean space, if a
(simply-connected closed) 2-cycle which was made by 2-cells of this triangulation separates $M$ into two connected components, then each of the components will also be simply-connected.}
The relationship of this theorem to a generalized Jordan separation problem, the general Jordan-Schoenflies theorem, and other important problems are also discussed.

In this paper, we assume that $M$ is orientable. In addition, a 2-cycle means that it is homeomorphic to a 2-sphere in this paper, and this homeomorphism can be determined in a reasonable time such as in the polynomial time.

We use $\cal M$ to denote a smooth or general topological manifold in Euclidean space and $K$ to denote a complex that is a
collection of cells based on a triangulation or decomposition. We also use $M$ to represent a restored manifold of $K$ in Euclidean space.

$M$ is essentially $K$, but $M$ is a discrete-like manifold (with $K$ as the decomposition) that will serve better for our discussion. $M$ is now embedded in a Euclidean space too, and each cell of $M$ with regarding to $K$ has the geometric meaning. For instance, every $i$-cell in $M$ with regarding to $K$, $i> 0$, has the centroid point that is different from all 0-cells contained in this $i$-cell in $K$. When we need, we can use this type of points
to refine the decomposition of $M$ ($K$) and the new one will be refinement of $K$. This refinement process will be performed in finite times only.

The Jordan Separation Property in higher dimensions has great deal of relationship to geometric topology. For 3D, if we admit the Poincare conjecture and  the general Jordan-Schoenflies theorem, then we will have the
Jordan Separation Property for 3D manifold. This is because, we can first make a compact, closed, simply connected 3D manifold to be 3D sphere. Then use the general Jordan-Schoenflies theorem.

In this paper, we admits 3D Jordan Separation Property for 3D manifolds, we will give an algorithmic proof for each compact, closed, simply connected 3D manifold is homeomorphic to a 3D sphere.

In fact the 3D Jordan Separation Property was proved by Chen and Krantz in 2015~\cite{Chen-Krantz}.

This paper is written based on the previous version.
{\footnote{
 In a recent paper, {\it Algorithms for Deforming and Contracting Simply Connected Discrete
Closed Manifolds (II)},  we discussed two algorithms for
deforming and contracting a simply connected discrete
closed manifold into a discrete sphere.  The first algorithm was a continuation of work that began in
{\it Algorithms for Deforming and Contracting Simply Connected Discrete
Closed Manifolds (I)}, the second algorithm contained a more direct treatment of contraction for discrete manifolds.
In this paper, we clarify that we can use this same method on standard
 piecewise linear (PL) complexes on the triangulation of general smooth manifolds. Our discussion is based on
triangulation techniques invented by Cairns,  Whitehead, and Whitney more than half of a century ago.}  }

\section{Some Basic Concepts of Manifolds}

A topological manifold is defined as follows \cite{LeeBook}:  A topological space $\cal M$ is a manifold of dimension $m$ if:
\newline
\hskip 1in (1) $\cal M$ is Hausdorff, and
\newline
\hskip 1in (2) $\cal M$ is second countable, and
\newline
\hskip 1in (3) $\cal M$ is locally Euclidean of dimension $m$.

We know that every manifold is locally compact, which means that every point has a neighborhood (containing this point) that is compact \cite{LeeBook}.

Let $x$ be a point in $\cal M$ and an open set containing $x$ is a neighborhood $\theta(x)$. Mapping $\theta(x)$ locally homeomorphic to an $E^m$ is called a $chart$. A collection of charts whose domain ($\theta(.)$ ) covers $\cal M$ is called an $atlas$ of $\cal M$.

Roughly, a manifold is smooth if for any point $x\in \cal M$ and each nearby point $y$ of $x$, the chart transition (called transition-map) from $\theta(x)$ to $\theta(y)$ is smooth. This means that $chart_y^{-1}(chart_x (\theta(x)\cap \theta(y)))$ is smooth. We omit the formal definition here for simplicity. In this paper, we only want to include local compactness as a property for conditions.

On the other hand, in \cite{Che15I}, we observed that the Chen-Krantz paper (L. Chen and S. Krantz,  A Discrete Proof of the General Jordan-Schoenflies Theorem,
{\tt http://arxiv.org/abs/1504.05263}) actually proved the following result essentially: {\it For a simply connected (orientable) closed manifold $M$ in discrete space $U$, if $M$ is a supper submanifold, (the dimension of $M$ is smaller than the dimension of $U$ by one), $M$ separate $U$ into two components.} Here, $U$ can be a decomposition or triangulation of Euclidean Space, a sphere, or a simply connected orientable manifold.   Further details can be found directly from the Chen-Krantz paper. If $U$ is a 4-dimensional Enclidean space, that is partitioned into grids\/cubic. $M$, a PL 3-manifold, separates $U$ into two components. The Chen-Krantz paper also provided a method that contracts
$M$ to be the boundary of a 4-cell. This means that $M$ is homeomorphic to the 3-Sphere in 4-dimensional Euclidean space.

However, when the dimension of $U$ is much bigger than
the dimension of $M$, then the contraction discussed in Chen-Krantz will not be applicable. We may need a new method. For instance, we could fill an $(m+1)$-manifold bounded by $M$, where $m$ is the dimension of $M$. Some algorithms have been discussed in \cite{Che15I}, but these algorithms might not work to certain cases.
In \cite{Che17I}, we continued the discussions by
\newline
\hskip 1in (1) Refining the algorithm by introducing a tree structure for filling, and
\newline
\hskip 1in (2) Designing another direct method of contraction.

This paper continues for the second part in \cite{Che17I}. In this paper, we attempt to use concepts such as the simplicial complex in combinatorial topology instead of the discrete manifolds used in \cite{Che15I,Che17I}, which makes this paper more readable for researchers and professionals in mathematics. We also attempt to draw a stronger connection between our results and the Poincare conjecture.

According to the official description of the Poincare conjecture made by John Milnor~\cite{Milnor}:

{\bf Question.} If a compact three-dimensional manifold $M^3$ has the property that every simple closed curve within the manifold can be deformed continuously to a point, does it follow that $M^3$ is homeomorphic to the sphere $S^3$.

According to Moise, any 3D topological manifold is homeomorphic equivalent to a compact smooth manifold ${\cal M}^3$, and there is only one diff-structure for any ${\cal M}^3$. So this question is the same to determine for a compact smooth manifold.

Morgan wrote a much clear explanation~\cite{Morgan06}: ``Now to the notion of equivalence. For topological manifolds it is homeomorphism.
Namely, two topological manifolds are equivalent (and hence considered the
same object for the purposes of classification) if there is a homeomorphism, i.e., a
continuous bijection with a continuous inverse, between them. Two smooth manifolds
are equivalent if there is a diffeomorphism, i.e., a continuous bijection with a
continuous inverse with the property that both the map and its inverse are smooth
maps, between the manifolds. Milnor’s examples were of smooth 7-dimensional
manifolds that were topologically equivalent but not smoothly equivalent to the 7-sphere. That is to say, the manifolds were homeomorphic but not diffeomorphic to the 7-sphere. Fortunately, these delicate issues need not concern us here, since in
dimensions two and three every topological manifold comes from a smooth manifold
and if two smooth manifolds are homeomorphic then they are diffeomorphic. Thus,
in studying 3-manifolds, we can pass easily between two notions. Since we shall be
doing analysis, we work exclusively with smooth manifolds.''

On the other hand, Lott explained the Poincare conjecture as ~\cite{Lott06}: ``The version stated by Poincare is equivalent to the following.

{\bf Poincare conjecture :}  A simply-connected closed (= compact boundaryless) smooth 3-dimensional manifold is diffeomorphic to the 3-sphere.''

Therefore, both Milnor and Lott used the compact manifold. It is certain that there is an equivalence between the compact manifold and general manifold for a simply-connected smooth 3D manifold.

In order to cover the topic more generally, we do not restrict our discussion to 3D compact manifolds.  Instead, we directly use  the
local compactness property of a topological manifold \cite{LeeBook}. We would like to obtain a sequence of PL complexes that approaches
the original manifold with any degree of approximation. Using this method, we can extend
the results to any type of high dimensions of smooth manifolds.

\section{Cell Complexes, Discrete Space, and Their Relationship}

In this paper, we want to make some modifications to the previous paper~\cite{Che17I} that was mainly based on discrete manifolds. Now, we would like to use cell complexes, specifically PL complexes in order to convey the same results. In this section, we review some related results.

Let $K$ be a cell complex that is a decomposition (triangulation) of a compact smooth $m$-manifold $\cal {M}$ in $n$-dimensional Euclidean space $E^n$. $K$ only contains a finite number of cells~\cite{Cairns35}. $M$ is the union of all cells in $K$, and $M$ can be viewed as a polyhedron in Euclidean space.
$M$ can also be viewed as a discrete manifold in ~\cite{Che15I,Che17I}. In other words, $M$ is a discrete manifold where each cell in $M$ can be found in $K$. Since $M$ is a polyhedron, $M$ is the restoration of $K$ in Euclidean space $E^n$.

Therefore, $K$ is a set of 0-cells, 1-cells, $\dots$, $m$-cells, where the set of $i$-cells is called the $i$-skeleton of $M$.
We must assume that there is no $(m+1)$-cell in $K$. We also assume that each $i$-cell in $M$ is an $i$-simplex, $i$-cube, or $i$-convex polyhedron. Without losing generality, we can assume that the intersection of any two cells in $K$ or $M$ is empty or a single $j$-cell where $j\ge 0$. To summarize, $K$ is the set-theoretic representation (just like the computer graphics input of a 3D object) of $M$, $M$ is the geometric representation of $K$.

We can also easily see that: An $i$-cell is an $i$-complex that contains a simply connected closed $(i-1)$-manifold plus the inner part of this manifold.  The simplest example is the $i$-simplex or $i$-cube.
In other words, an $i$-cell is bounded by a set of $(i-1)$-cells, and it contains the inner part of the cell in the $i$-dimension.
The inner part of an $i-$cell usually refers to an open $i$-cell in $E^n$.

The construction of $K$ is interesting, practically. We can see that
$K$ contains all inner parts of cells contained in $K$. However, to represent $K$ (in computer or any other form), we only need to attach the 0-cells with the inner parts of cells.

In fact, $K$ only contains the inner part of $i$-cells (open $i$-cells) when $i\ge 1$ and all 0-cells (points). In order to make our discussion easier, let $K$ contain all open and closed cells from the triangulation or decomposition of $\cal M$.

Computationally, let $a$, $b$, $c$ be 0-cells.
To represent a complex in computers or any other explicit form, we can use $[a,b,c,...]$ to express a cell that also contains the boundary of
the cell.  $(a,b,c,...)$ is an open set that contains the inner part of $[a,b,c,...]$. Let $<a,b,c,...>$ represent the boundary of $[a,b,c,...]$. Then, we have

  $\partial ([a,b,c,...]) = <a,b,c,...> = [a,b,c,...]\setminus (a,b,c,...)$,

\noindent and

  $[a,b,c,...] =  <a,b,c,...> \cup (a,b,c,...)$.

Because each element of $K$ in $E^n$, such as $[a,b,c,...]$,  $(a,b,c,...)$, and $<a,b,c,...>$, has real meaning, we can
assume that $K$ contains not only open cells but all forms including “$[.]$” “$(.)$” and”$<.>$” of cells for convenience.
The advantage of these types of representations is easy for coding in computers. For instance, to determine the centroid of a triangle, we can map the corner points to 2D plane and then calculate that. Each point has
a real Euclidean location.

A curve in $M$ is a path of 1-cells in $K$, so we can define the simple connectivity of $M$. $K$ contains a finite number of cells. Usually no cell can be split into two or more if it is not specifically mentioned, which means when $K$ is completely defined, we will usually not make any further changes. (In fact, to make finite numbers of refinements or merges on cells will not affect much of the proofs in this paper. )

Since every pair of two vertices (0-cells) in finite $K$  have a certain distance in $E^n$, local flatness will not be a problem for a polyhedron since it is always locally flat in $E^n$. To recall the definition of local flatness: $S$ is called locally flat in $E^n$ (or other space) if for any point $x$ in $S$, there is a neighborhood
$U$ of $x$ in $E^n$ such that $U\cap S$  is a neighborhood of $x$ in $S$. It roughly means that $U\cap S$  is homeomorphic to a (deformed) disk of $S$ containing $x$.  Therefore

{\bf Proposition 3.1 }: $M$ as a polyhedron is locally flat in $E^n$.

{\bf Proposition 3.2}: Every $i$-submanifold $M(i)$ determined by $i$-skeleton of $K$ is locally flat in $E^n$.

Local flatness here is not the same as {\it discrete local flatness} in ~\cite{Che17I}. When we require discrete local flatness in this paper, we can always to use the barycentric (or other type of) subdivision of triangles to create one. However, since we use embedded manifolds in $E^n$, there
is no need for such a process.  We know that the inner part of a $j$-cell ($j>0$) can be determined in $E^n$ as the centroid of the cell (or
the counterclockwise direction of $(j-1)$-cells in the $j$-cell; just think about algebraic representation of triangulations in algebraic topology ~\cite{Hat}).  Unlike in pure discrete space, there is no centroid point in discrete space. Discrete space only contains the rotational directions (one can indicate inner and another can indicate outer).

Therefore, the case of $K$ with its embedding counterpart $M$ in $E^n$ is much simpler than the discrete manifolds for pure discrete cases.

One concept that is called {\it cell-distance } was developed in the Chen-Krantz paper (L. Chen and S. Krantz,  A Discrete Proof of The General Jordan-Schoenflies Theorem, {\tt http://arxiv.org/abs/1504.05263}).  the cell-distance is for two points in a manifold: Two points $p$ and $q$, $q\ne q$ in $K$ or $M$ is said to have a $j$-cell-distance $k$ means that the length of the shortest $j$-cell-path that contains both $p$ and $q$ is $k$.
We will repeat this definition when we actually use it in this paper.

Again, we assume that all manifolds discussed in this paper are orientable.

Some of the concepts in this paper are also included in \cite{Che14}. In addition, we include the concept of the cell distance from the following, see \cite{Chen-Krantz}:

In a graph, we refer the concept of “distance” as the length of the shortest path between two vertices.
The concept of {\it graph-distance} in this paper refers to edge distance, meaning how many edges are required to get from one vertex to another. We usually use the length of the
shortest path between two vertices to represent distance in graphs. In order to distinguish between this definition of distance and distance in Euclidean space,
we use graph-distance to represent length in graphs in this paper.

Therefore, graph-distance is the same as edge-distance or 1-cell-distance (how many 1-cells are needed to travel from $x$ to $y$).  We can generalize
this concept to define 2-cell-distance by counting how many 2-cells are needed to go from point $x$ to point $y$. In other words,
2-cell-distance is the length of the shortest path of 2-cells that contains $x$ and $y$. In this path, each adjacent pair of
2-cells shares a 1-cell and this path is 1-connected.

We can define $d^{(k)}(x,y)$, the $k$-cell-distance from $x$ to $y$, as the length of the shortest path (or the minimum number of $k$-cells in such a sequence)
where each adjacent pair of two $k$-cells shares a $(k-1)$-cell. This path is $(k-1)$-connected.

We can see that the edge-distance or graph-distance, called $d(x,y)$, is $d^{(1)}(x,y)$.

We can also define $d^{(k)}_i(x,y)$) as a $k$-cell path that is $i$-connected. However, we do not need to use such a concept here.

\section{Contraction of the Finite PL Complex $K$}

Let $M$ be a smooth compact simply connected closed manifold in $E^n$. The complex $K$ is a triangulation or PL decomposition of $M$. $K$ only contains finite number of cells.  In other words, assume that all partitioned PL cells of $M$ are included in $K$ where $K$ is a finite simplicial complex or PL complex of $M$. We will only deal with $K$ and its elements, not every point in $M$ that is a subset of $E^n$.

Based on a theorem, Theorem 5.1, proved in Chen-Krantz paper, we concluded that a simply connected closed $(m-1)$-manifold $S$ split a simply connected closed $m$-manifold $M$ into two components in discrete space. We required both $S$ and $M$ be locally flat in discrete space. For continuous space, it is a true property for both $M$ and $S$ in $E^n$ as we discussed before. In other words, $M$ in $E^n$ is locally flat. Any $i$-submanifold $S$ of $M$ is locally flat if each cell of $S$ is an element in $K$ ($K$ was defined previously). As we discussed in last section, $M$ is locally flat in $E^n$.

On the other hand, one of the key facts in this paper is that we allow simplexes or cells with open and closed-cells in $E^n$, we do not need to consider local flatness for discrete cases. This is because we can easily partition a cell (few times) into more cells to meet the requirement of discrete local flatness.

We restate this theorem, Theorem 5.1, in \cite{Chen-Krantz}, as follows:
\begin{theorem}
(The Jordan theorem for the closed surface on a 3D discrete manifold) Let $M$
be a simply connected 3D manifold (discrete or $PL$); a closed discrete surface $S$ (with
local flatness) will separate $M$ into two components. Here, $M$ can be closed.
\end{theorem}

\subsection{The 3D Jordan Theorem in Euclidean Space}

Note that $M$ in this paper is slightly different from $M$ in Chen-Krantz paper. $M$ in Theorem 4.1 is a discrete manifold; $M$ in the current paper
is the combined form of $M$ and $K$. Local flatness is required in Theorem 4.1 (Theorem 5.1 in Chen-Krantz paper),
However, in this paper, $(M,K)$ is naturally locally flat. This is because a (finite) $K$ can be subdivided a finite number of times: Every cell (simplex), except 0-cells, must have a centroid point
that is distinguished from the boundary of the cell. When we consider
$(M,K)$ in this paper, any simple curve (meaning that no 0-cell appears twice in the curve) or simple (sub-)manifold is locally flat.
Therefore, based on Theorem 4.1, we can easily get:

\begin{theorem}
Let $M$ be a simply connected closed 3D manifold associated with the finite complex $K$ where $K$ is a triangulation or a PL decomposition of $M$. A simple closed 2D-surface $S$ in $M$ (where each 2-cell of $S$ is a 2-cell in $K$) will separate $M$ into two components.
\end{theorem}

\noindent{\bf \it  The Brief Proof:}

We just prove this theorem here briefly.
There is a little difference from Theorem 4.1 and Theorem 4.2. This is because Theorem 4.1 deals with discrete manifold, and Theorem 4.2 is for an embedded manifold in Euclidean space. From the following explanation, we can also see that the embedded 3-manifolds are simpler than the discrete manifolds.

Since $S$ is finite in Theorem 4.2 and $S$ is naturally locally flat in Euclidean space, so $S$ separates $M$ into two parts as we have proved in Chen-Krantz paper with at most two- or three-times refinements by (barycentric or other type) subdivision  if it is needed when original $S$ is not discrete locally flat in $M$. After the refinement, we can directly use the Theorem 4.1.

                                                                                                                                                                $\hfill$ $\square$ \\

The separation by using $S$ to $M$ here is more significate, the connectedness about the two parts after the separation could be somewhat easier in our intuition.
We also could have the following observation for connectedness without thinking to perform refinements using subdivision: Since $S$ is a simply closed surface by 2-cells in $K$, then

(1) If $S$ is not the boundary of just one 3-cell in $M$,  any 3-cell in $K$ in Theorem 4.2 must contain a face (2-cell) that is not totally on $S$.  This 3-cell can connect to other 3-cells without using any 2-cells in $S$.
Since each (open or closed) 2-cell of a 3-cell will be included in two 3-cells in $M$ ($M$ is closed), let 2-cells (not in $S$) to glue 3-cells. We will get one connected part. So, each of the two components (parts) in  Theorem 4.2 is connected.

(2)If $S$ is the boundary of a 3-cell in $M$. Then $S$ still separates $M$ into two components since one is the inner part of this 3-cell and the another is the rest of 3-cells of $M$. \\

\noindent A complete proof can be found in Chen-Krantz paper. Here just provide some intuition of the proof.  Thus, we can see that the embedded manifold in Euclidean space is simpler than the discrete cases.



Now we will give an algorithmic proof for the following statement: If $M$ is closed in Theorem 4.2, then $M$ is homeomorphic to a 3-sphere.

\subsection{The Main Theorem}

\begin{theorem}
If $M$ is closed in Theorem 4.2, then we can algorithmically make $M$ to be homeomorphic to a 3-sphere.
\end{theorem}

\noindent{\bf \it The Algorithmic Proof:}

According to Theorem 4.2, a simply connected closed $2$-manifold (orientable) $B$ will split the $3$-manifold $M$ into two components $D$ and $D'$.
$B$ is the common boundary of $D$ and $D'$. ($B$ was $S$ in Theorem 4.2.) $B$ only contains cells in $K$ which is the simplicial complex or  PL-complex corresponding to $M$ in Euclidean space as defined previously.
(Note that, elements of $B$ will only select from the simplicial complex $K$.)

In this proof, there two major tasks: In the first part, we prove that $D$ and $D'$ are simply connected. In the second part, we prove that $D$ (and $D'$) is homeomorphic to a 3D-cell (or an $m$-cell) based on the first part of the proof.

\noindent {\bf Part 1:} We first show that each of the components ($D$ and $D'$) will be simply connected.

According to the definition of simple connectedness, if $M$ is simply connected then a simple closed curve $C$ on $M$ is contractible to a point $p\in C$ on $M$. Here, $C$ only contains edges (1-cells) in $K$.

Let $\Omega$ be the contraction sequence in the discrete case, which we call gradual variation in \cite{Che14}. We might as well to assume that $p$ is not on $B$.

Again, this contraction sequence only takes edges (1-cells) in $K$. Each curve here is
a closed path of edges in $M$. In addition, $B$ is a subset of $K$, so $B$ only contains a finite number of cells including a finite number of edges as well as 0-cells and 2-cells.

\begin{figure}[h]
	\begin{center}
    \psfrag{Omega}{{\bf $\Omega$}}
   \epsfxsize=4.2in
   \epsfbox{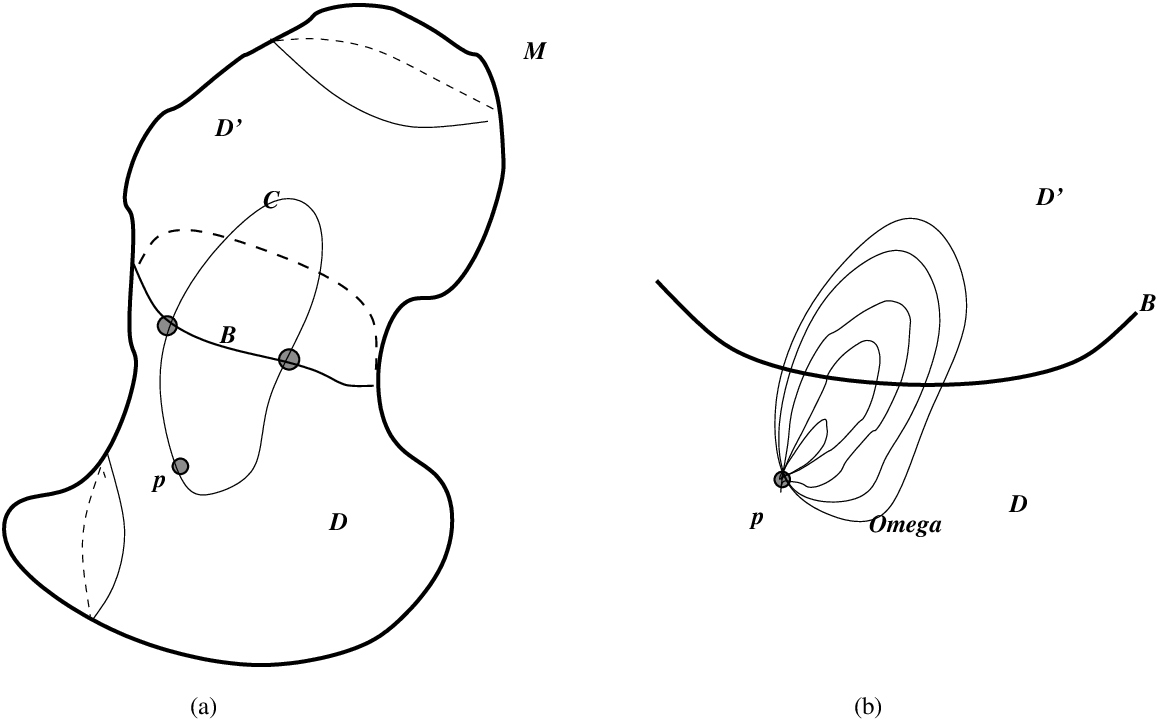}   
\caption{An example of $M$ that is separated by $B$.}
\end{center}
\end{figure}

The contraction sequence $\Omega$ may contain some point on $B$ and pass $B$ to another component (and come back), see Fig. 1(a)(b). In such a case, we want to modify this contraction sequence that will replace
$\Omega$. It means to get a new contraction sequence in $D\cup B$ (or $D'\cup B$) to show that $D$ (or $D'$) is simply connected.  We start the detailed proof in the following.

First, a curve $C$ on $\Omega$ may intersect with $B$ but not crossover $B$ (this means that $C \subset D\cup B$ or $C \subset D'\cup B$). We will not make changes for this case here, see Fig.2 (a).
(That is to say if the curve only stays on $B$ and returns to the original component $D$, then this will not count as ``a pass.'')

Note that we usually think that $D$ contains $B$ and $B=\partial D$. We use $D\cup B$ is to ensure that no matter what we think $D$ is open or closed in Euclidean space, our proof is still clear.

Second, we assume that $p\in D$, some curves in $\Omega$ pass $B$ to $D’$ and then back to $D$.

Let $p\in D$. We want to modify the case that a closed curve $C$ containing $p$ in $\Omega$ starts at $p$ and passes $B$ to reach a point in $D'$ then come back to $p$ in $D$, see Fig. 2(b).

We define {\it a pass} regarding to $B$ to be entering $B$ from one component and exiting $B$ to another component. For instance, a sequence of points $d, b_1,\cdots,b_k,d'$ in $C$ is {\it a pass} if $d\in D$, $d'\in D'$ (or {\it vise versa}), and $b_i\in B$ where $i=0,\cdots,k$ and $k\ne 0$.  (We also think about that the curve $C$ travels counterclockwise.) See Fig. 2 (b). Note that $C$ is always 1D. $B$ can be 1D, 2D or high-dimensional.

\begin{figure}[h]
	\begin{center}
  \psfrag{Omega}{{\bf $\Omega$}}
   \epsfxsize=5in
   \epsfbox{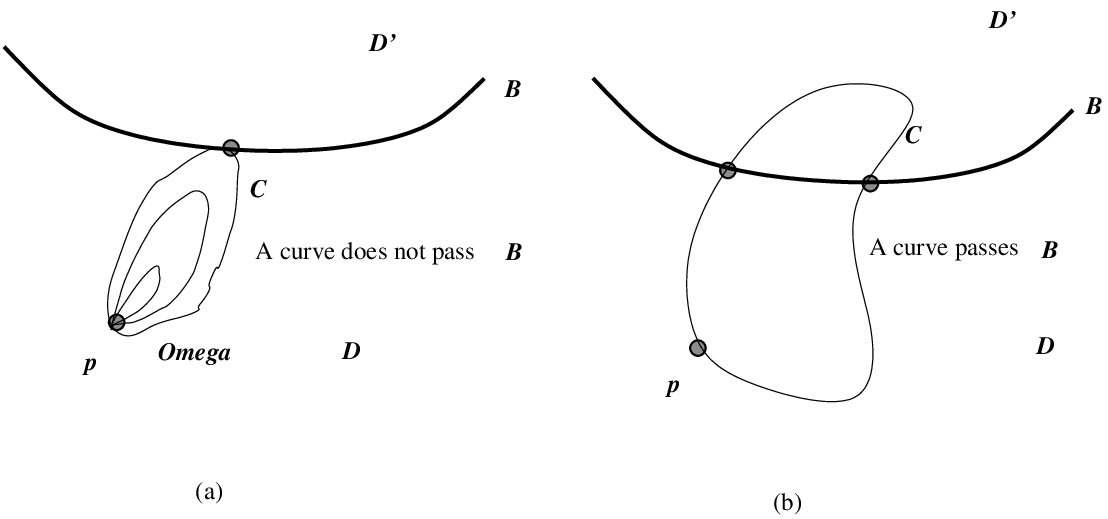}   
\caption{The curve $C$ and the boundary $B$: (a) A curve does not pass $B$. (b) A curve passes $B$. }
\end{center}
\end{figure}

\noindent Note that $C$ contains only a finite of edges (1-cells) in $K$, thus $C$ can only pass $B$ finite times.
In fact, since $C$ starting at $p$ and ending at $p$, $C$ pass $B$ an even number of times including 0 pass (no pass).

Let $C'\in \Omega$ containing $p$ always start from $D$ (meaning that
$p\in D-B$).  When $C'$ passes $B$, it will (later) enter $B$ from $D'$ in order to reach $p\in D-B$ again in the cycle $C'$ (according to Theorem 4.2, $B$ separates $D$ and $D'$). Therefore, the last-out and first-in in $C'$ (regarding to $B$) are always a pair. See Fig. 3.

\begin{figure}[h]
	\begin{center}
   \epsfxsize=6in
   \epsfbox{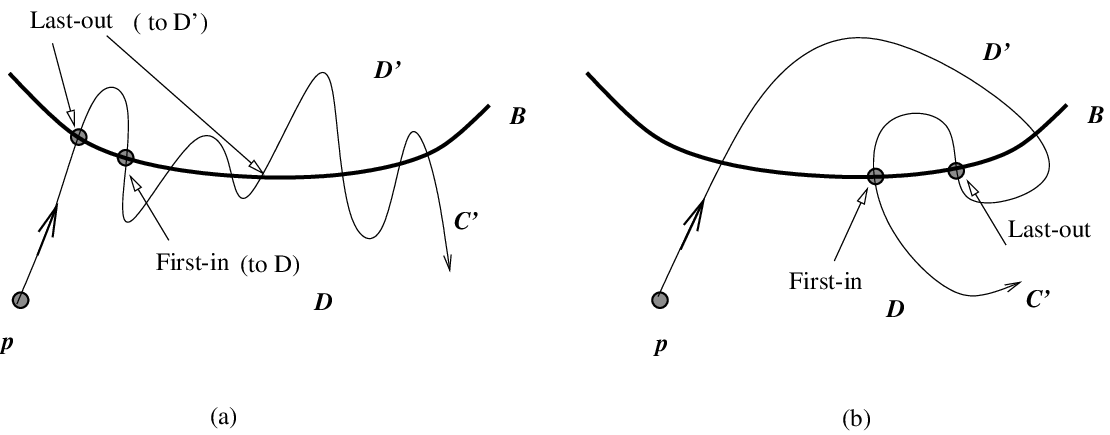}   
\caption{ Two cases of the passing pairs: (a) A simple case with multiple passes. (b) A more complex case. }
\end{center}
\end{figure}

We now present some observations and the idea for proofs:

{\it We can find a curve or an arc in Fig. 3(a), $C'_B(0)$, on $B$ to link the following two points (since $B$ is also simply connected): the last point leaving $B$ and the fist point entering $B$ from $D'$ (as a pair of points ).   See Fig.3 and Fig. 4(a). Using $C'_B(0)$ to replace the corresponding arc in $D'$ (they form a cycle in $D'\cup B$), if there are multiple passes, then we can use $C'_B(1),...,C'_B(i)$ to replace the whole $C’$ (sequence). Therefore, we obtain a new $C'_B$ that only contains points in $D\cup B$. See Fig. 4 (b).

Note that, an arc is an unclosed simple curve in this paper.

Following the same process for all curves in $\Omega$, we will get an $\Omega_D$ that is a contraction (sequence) to $p$. Therefore, $D$ is simply connected. For the same reason, $D'$ is also simply connected.}

\begin{figure}[h]
	\begin{center}
    \psfrag{alpha}{{\bf $\alpha$}}
    \psfrag{beta}{{\bf $\beta$}}
   \epsfxsize=6in
   \epsfbox{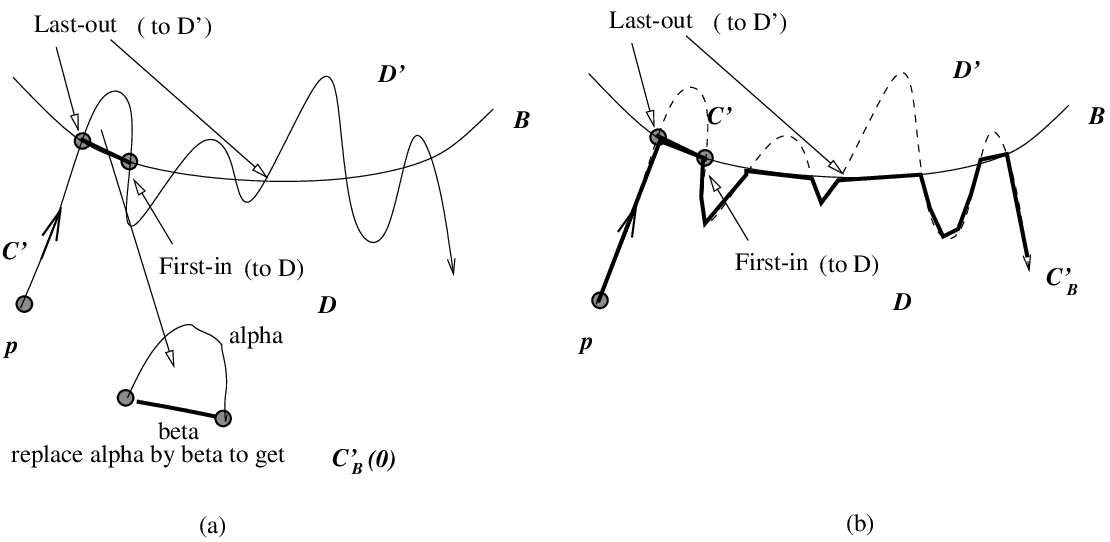}   
\caption{A modified curve only in $D\cup B$: (a) The original curve $C'$, and (b) The modified curve $C'_B$   }
\end{center}
\end{figure}

Now, we start the detail construction procedure for above discussion: Since ``the last-out and first-in in $C'$ (regarding to $B$) are always a pair,'' we can assume that they are: $d^{(1)} b^{(1)}_1\cdots b^{(1)}_{k_1} d'^{(1)}$  and $d'^{(2)} b^{(2)}_1\cdots b^{(2)}_{k_2} d^{(2)}$.  Here, $d^{(1)},  d^{(2)}$ are in $D-B$
and  $d'^{(1)}, d'^{(2)}$ are in $D'-B$;   $b^{(1)}_i$ and $b^{(2)}_j$ are in $B$. See Fig. 5.

\begin{figure}[h]
	\begin{center}

    \psfrag{alpha}{{\bf $\alpha$}}
    \psfrag{beta}{{\bf $\beta$}}
    \psfrag{d1}{$d^{(1)}$}
    \psfrag{d2}{$d^{(2)}$}
    \psfrag{d'1}{$d'^{(1)}$}
    \psfrag{d'2}{$d'^{(2)}$}
   \includegraphics [width=0.8\columnwidth] {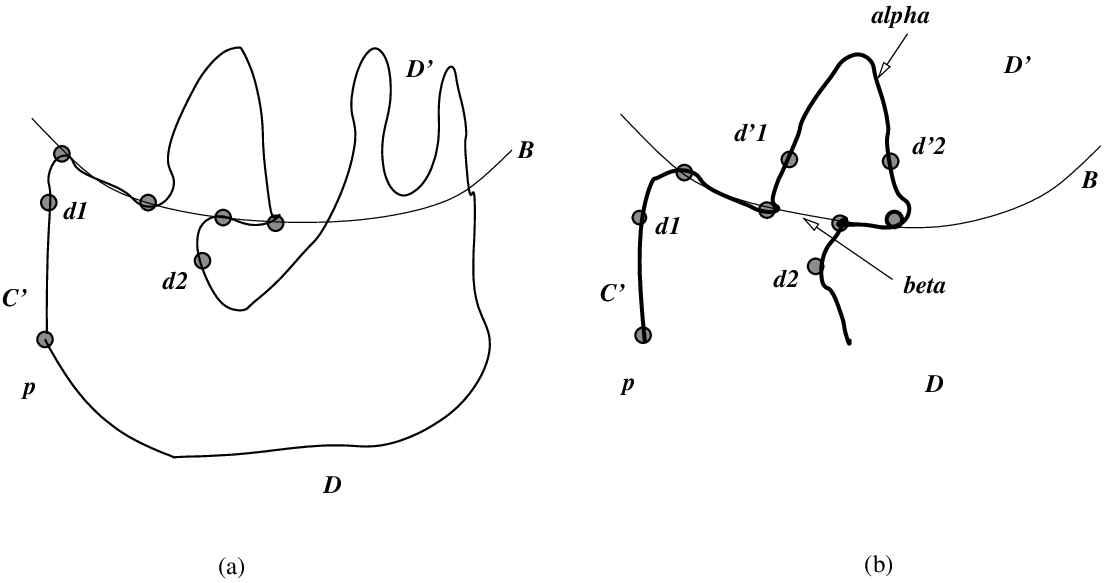}   
  \caption{More details in the process of replacing every arc in $D'$ with an arc in $B$: (a) A curve, and (b)Identifying arcs in $D'$ and $B$. }
\end{center}
\end{figure}


We can see that $b^{(1)}_{k_1}$ is the leaving point in $B$ and $b^{(2)}_1$ is the entering point in $B$ in Fig.5. Since $B$ is also simply connected, there is a path of edges (1-cells) between $b^{(1)}_{k_1}$ and $b^{(2)}_1$ inside of $B$. We denote this path as $\beta$ (See Fig. 4 (a) and Fig. 5(b)). Let $C' =\cdots b^{(1)}_{k_1} d'^{(1)} \cdots \cdots \cdots d'^{(2)} b^{(2)}_1 \cdots$  be denoted as $\cdots b^{(1)}_{k_1} d'^{(1)} \alpha d'^{(2)} b^{(2)}_1 \cdots $. Here, $\alpha$ is a path (an arc) in $D'$.  $\alpha$ could be empty and $d'^{(1)}$ could be  $d'^{(2)}$.

Note also, $\alpha$ and $\beta$ form a cycle (a simple closed curve made by 1-cells in $K$).  We use $\beta$ to replace $d'^{(1)}\alpha d'^{(2)}$. (Regardless if there is a gradually varied sequence to make $d'^{(1)}\alpha d'^{(2)}$ to be $\beta$ without using points in $D$, we do not need that. Just do the simple replacement. It is important in this paper. We also refer that $\beta$ is a projection
of $\alpha$ to $B$. )

Repeat the above process to find the new last point leaving $B$ and the new fist point entering $B$ (from $D'$) until all points in $C'\cap D'$ are replaced by points in $B$. We obtain a new $C'_B$. Use $C'_B$ to replace $C'$ in $\Omega$, we will get a $\Omega_D$.  The main task is now to prove that $\Omega_D$ is a gradually varied sequence (discretely ``continuous'' sequence) or to make
that $\Omega_D$ is a gradually varied sequence.

In the above construction (or algorithm), we can see that (1) we did not change $C'$ in $D-B$ and all entering points in $B$ (from $D$), and (2) Since $\Omega$ is a gradually varied sequence in $M$ (each curve must contain $p$ in the sequence based on the definition of contraction in discrete case, see Chen-Krantz paper), then $\Omega_D$ is a gradually varied sequence in $D-B$.
However, we still have a problem: {\it $C'_B$ might not be a simple curve inside of $B$,} see Fig. 6.

%
\begin{figure}[h]
	\begin{center}
    \psfrag{d1}{$d^{(1)}$}
    \psfrag{d2}{$d^{(2)}$}
   \epsfxsize=3in
   \epsfbox{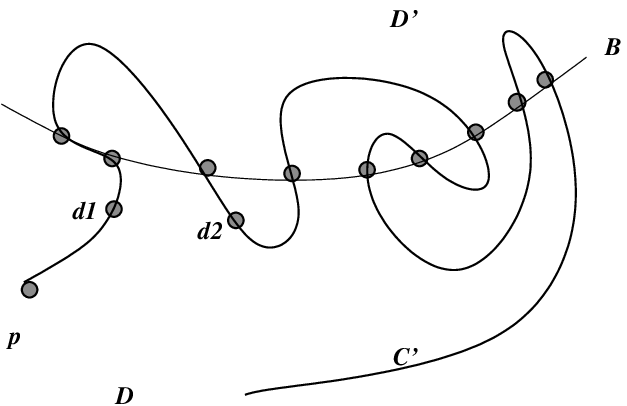}   
\caption{ The new modified curve $C'_B$ might not be a simple curve especially when $B$ is a 1-cycle. }
\end{center}
\end{figure}

Now we want to solve this problem as follows:

What we want is that the replacing arcs or curves in $B$ will form the simple arcs or curves in $B$ and they should not intersect each other.  We have two cases here:

(1) $B$ is 1-dimensional, and

(2) $B$ is 2-dimensional or higher-dimensional.

When $B$ is 1-dimensional then $M$ is 2-dimensional, we have already used a contraction method in \cite{Che14} to prove that $D$ and $D'$ are both simply connected (or we just use 2D solution for 2D Poincare case and the Jordan-Schoenflies Theorem to resolve it.). By the way, in this case, we use direct contraction based on the 2-cell distance. This technique also used in Chen-Krantz paper \cite{Chen-Krantz}. As long as a manifold is a supper submanifold of an ambient space, we can use it.  {\footnote {We are working on the research for algorithmic contraction of $k$-cycle on $(k+1)$-manifolds with a referencing $k$-D $B$ as the separation cycle.}}

Now, we just discuss the case that {\it $B$ is 2-dimensional or higher dimensional}.  Let's just assume that $B$ is simply connected closed 2-manifold (similar to higher dimensional cases). Since it is also a part of $M$. $B$ is embedded in Euclidean space.

We know that $C'$ may travel to $B$ multiple times. However it only take different point in $B$ since $C'$ is a simple curve. Every time we use an arc $\beta$ in $B$ to replace an arc in $D'-B$, (note that in above discussion we have used ``$b^{(1)}_{k_1}$ is the leaving point in $B$ and $b^{(2)}_1$ is the entering point in $B$'' as a pair.) we want to keep  $\beta$ do not intersect any of the existing part of $C'$ (and its modification) in $B$ at the time of we travel with $C'$ (starting from $p\in D-B$ along with counterclockwise order).
In other words, the projection from $\alpha$'s to $\beta$'s may cause intersections among $\beta$'s. Some new $\beta$'s may cross each other.   We want to design a procedure to guarantee them do not cross each other in $B$.

The previous arc of $C'$ (and its modification) we have modified (so far) is assumed that there is no intersection. Assume that $C'=C'_B(0)$, $\cdots$, $C'=C'_B(i)$ were processed.

Then we make the current $\beta$ do not contain any $C'_B(i)$ for all $i$ previously processed. We know two ending points of  $\beta$ are at the two distinguished points in $B$. We also know some of the points in $B$ are taken by previous process.

Let $\rho_0$ and $\rho_1$ be two ending points in  (current) $\beta$. We want to prove that there is a path in $B$ that can link $\rho_0$ and $\rho_1$ without intersect with any points of $C'$ and its modification in $B$. See Fig. 7 (a) and (b). That is our objective of the rest of the proof for Part I.

\begin{figure}[h]
	\begin{center}
    \psfrag{rho1}{$\rho_0$}
    \psfrag{rho2}{$\rho_1$}
   \epsfxsize=6in
   \epsfbox{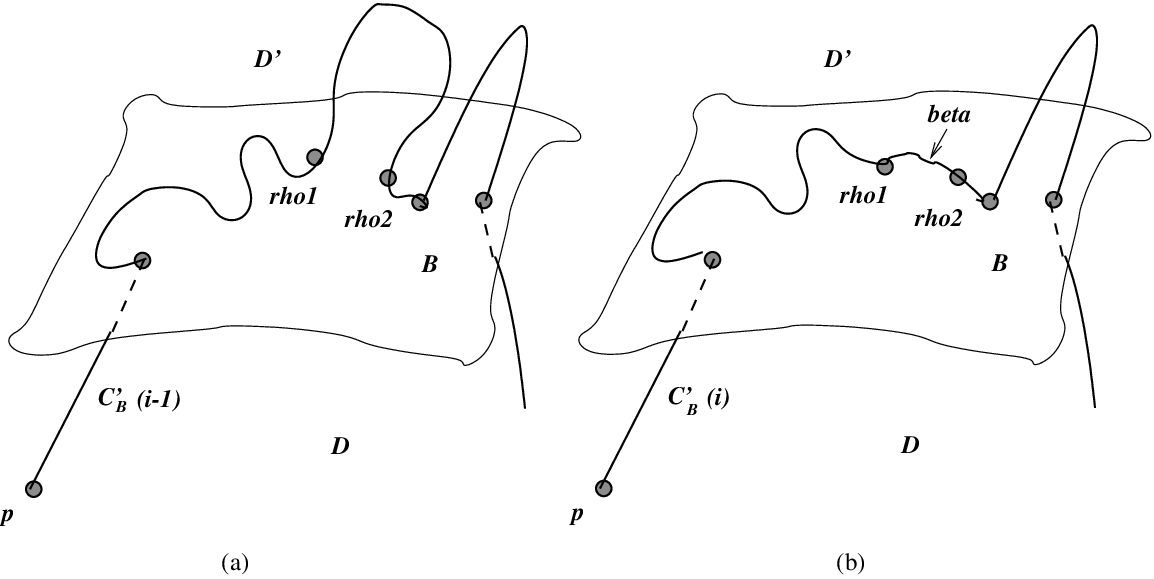}   
\caption{The new curve to make: (a) the old curve, and (b) the new curve.}
\end{center}
\end{figure}

Note that the above two paragraphs are important for the idea of the our proof. We might as well to assume now $B$ is a 2-cycle in the following proof while $M$ is a 3D manifold. Actually we mentioned it
a few paragraphs earlier. A 2-cycle is defined as a 2-manifold that is homeomorphic to 2-sphere. We know that the 2D Poincare conjecture is true. So $B$ is a 2-cycle if $M$ is 3-dimensional.


\noindent {\bf Subpart 1.1: The detailed construction for the projection of a single curve to $D\cup B$ } 


According to the Jordan curve theorem, a simple closed curve will separate $B$ into two components if $B$ is 2D. We will prove that between  $\rho_0$ and $\rho_1$, there is no such a closed curve (a path of 1-cells) in $C'$ (not $C'$ itself) and its modification in $B$. (Note that this is based on our special construction of the new curve $\beta$ which we will show the details below.) See Fig. 7. (a).


On the contrary, we assume there is a such closed curve (or subcurve).

{\bf (1.1.a):} Assume that this closed curve was only made by a pair of a leaving and an entering point respecting to $D'$. (Meaning that there is only one $\alpha$ in the closed curve.)  In order to make a closed curve (a path of 1-cells ) in $C'$ in $B$, there must be a point in $B$ that was traveled twice by {\it $C'$ and its modification}. It is impossible since $C'$ is a simple closed curve that contains $p\in D-B$. $C'$ cannot contain another closed curve in $B$. There are must be a (open) 1-cell in $B$ that is not in $C'$  and its modification.

On the other hand (as another explanation), if there is an arc in $C'$ from the leaving point $\rho_0$ to the entering point $\rho_1$ totally in $B$, (another one inside of $B$, but we have one in $D'$ already.) there will be a closed path (cycle) in $C'$ . Note that this closed path is also totally in
$B\cup D'$. It is impossible since $C'$ contains a point $p$ in $D-B$, an proper subset of $C'$ can not contain any closed path since $C'$ itself is a simple path.

The above paragraph means that $C'$ may cover many 1-cells in $B$; however, $C'$ could not cover all 1-cells of $B$ since
$C'\cap B$ can not contain any closed curve. We can see that it is possible to find an arc  $\beta$ in $B$ to replace $\alpha$ in $D'$ (even though we may need to refine some cells, we will discuss it later), see Fig. 7.

{\bf (1.1.b):} Assume that this closed curve was  made by two or more pairs of leaving and entering points respecting to $D'$. (Meaning that there are multiple $\alpha$'s in the closed curve.) Since no point in $B$ was used twice. So there must be a 1-cell (an edge) was not taken by {\it $C'$ and its modification}.  It is possible for $C'$ to take all 0-cells in this path. However, there are must be a (open) 1-cell that is not in $C'$ and its modification. In other words,  $C'$ could cover all original vertices or points in $B$ but the travel path of $C'$ will not be able to take all 1-cells in $B$.

Let us continue based on (1.1.a) and (1.1.b).
Since we embedded $M$ in Euclid space, the middle point of such a 1-cell mentioned above is not on {\it $C'$ and its modification}. (We want to make that) $\rho_0$ and $\rho_1$ are connected in $B$ without intersect with  {\it $C'$ and its modification} even though we may
need to use some the middle points of 1-cells in $B$.

Since $B$ is also a closed discrete 2-manifolds, any two 2-cells in $B$ are 1-connected. See~\cite{Che14}. The concept of 1-connectedness of 2-cells means that there is a sequence of 2-cells, the two adjacent two 2-cells in the sequence share a 1-cell.

Assume there is a path of 2-cells in $B$ that link $\rho_0$ and $\rho_1$, each adjacent pair in the path share one 1-cell that is not in $C'$ and its modification.
We only need to do subdivision for these 2-cells in $B$ (in the case of most of edges are taken by $C'$. Otherwise we do not need to do many subdivision.). We can make a simplest subdivision by cutting a 2-cell into two 2-cells. Or at most, we use the center point of the 2-cell to cut it into finite pieces
(e.g. barycentric subdivision, it can also be done for any simplex in any dimension. See Fig. 8).   After the subdivision, we need to update $K$ first. Then, we link $\rho_0$ and $\rho_1$ without intersect with  $C'$ and its previous modification. The following algorithm will construct such a path:

Algorithmically, to find this special path from $\rho_0$ and $\rho_1$, we can start with {\it the star of $\rho_0$} in $B$, there must be an edge
in {\it the link of $\rho_0$} in $B$ that is not in {\it $C'$ and its modification}. Use this edge, $aEdge$, to find the two 2-cells containing this edge. Actually we only need to find another
2-cell that contains this edge $aEdge$. Repeat this process at an new edge not in {\it $C'$ and its modification} and this process will mark all cells in $B$. We could stop if a 2-cell containing $\rho_1$ is reached or marked. In computation, using breadth-first-search will find the almost shortest path that do not intersect with $C'$ and its modification. Or using depth-first-search, we can also make it.  So we find the path of 2-cells with 1-connectedness. Except two ending 2-cells, we will partition the middle edge and its associated 2-cells. For the ending two 2-cells, we will use centroid of
the cell to do a subdivision. So we can link $\rho_0$ and $\rho_1$ without intersect {\it $C'$ and its modification} in $B$. We will give an example in a figure later.

Therefore, with the subdivision of some 2-cells (finite times, at most the number of pair of passes of $C'$ regarding to $B$), we can find $C'_{new}$ in $D\cup B$ that is a simple curve. The part of $C'\cap (D-B)$ is unchanged.

When doing subdivision of $B$, we also need to make associated 3D (or corresponding $m$-cells) in $D$ to be divided accordingly. It is easy to see that if $C'$ and $C''$ only differs from a 2-cell $\Delta \in K$ (or any $k$-cell in $K$). It means that $(C'\cap \Delta)\cup (C''\cap \Delta)$ is a 1-cycle (e.g. the boundary of the 2-cell, or 1-cycle on the $k$-ell.), and $(C'\cap \Delta)\cap (C''\cap \Delta)$ is two 0-cells. Then for any finite subdivision of $\Delta$, there will be a sequence of gradually varied path in between $C'$ and $C''$ for the updated $K$.

To summarize, this part is to find a new closed curve in $D\cup B$. This curve is modified based on $C'$. Plus, $C'$ in $D$ (i.e. $C' - {B\cup D'}$) is unchanged. In the next, we will find a sequence of such new curves that are contracting to the point $p$ without using any points in $D'$. \\

\begin{figure}[h]
	\begin{center}
   \epsfxsize=1.2in
  \epsfbox{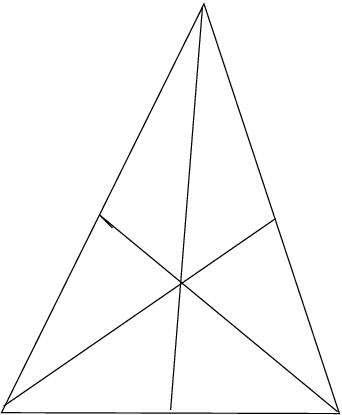}   
\caption{ Barycentric subdivision. (a) 2D, and (b) 3D.}
\end{center}
\end{figure}

\noindent {\bf Subpart 1.2: The detailed construction for a sequence of simple curves}

Now, we need to check for the contraction sequence $\Omega$. We know that if $C'$, $C''$ are two adjacent simple curves in the sequence,
we can also assume that there is only one 2-cell that are changed in the contraction, i.e. $Modulo2(C', C'')$ is the boundary of only one 2-cell. Since $B$ is a PL 2-manifold, there are only three cases: (a) this 2-cell is not intersect with $B$, (b) it is entirely in $B$,  or (c) this 2-cell has a part (an arc) in $B$ and the other part is in $D$ or $D'$ not both. (See Chen-Krantz paper , or just think it is a triangle.)

For simplicity of the proof, let us just assume that a cell is a simplex in the rest of the proof. In other words, each 2-cell is a triangle.  For example, $b_0 b_1 b_2$.
We can do analysis for the three cases:

{\bf (1.2.a):} If the change is not related to $B$, we just use $C'\cap B$ to be $C''\cap B$. See Fig. 9.

\begin{figure}[h]
	\begin{center}
   \epsfxsize=4.5in
   \epsfbox{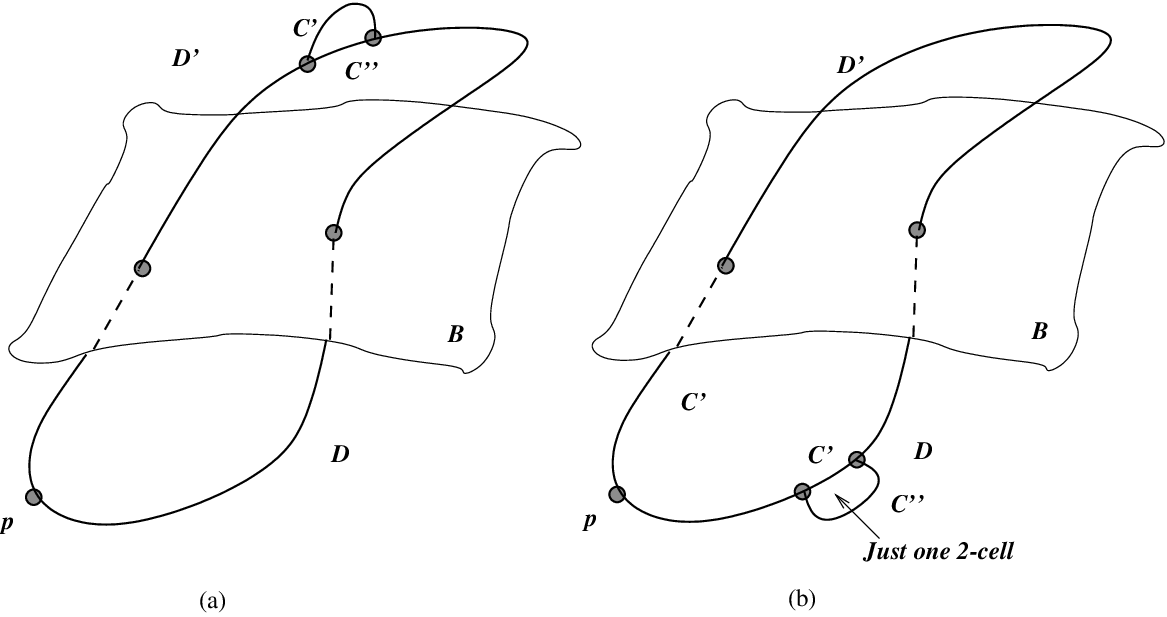}   
\caption{Two adjacent curves in $\Omega$: (a) Changes in $D'$ , (b) Changes in $D$ }
\end{center}
\end{figure}

{\bf (1.2.b):} If the change is totaly on $B$,  only an arc (on the boundary of a 2-cell) on $B$ was changed at most (See Fig.10). We now just use a triangle as the example, $b_0 b_1 b_2$ where $b_0$ and $b_2$ in both $C'$ and $C''$. See Fig.11 (a).  If $b_1$ is in $C'$ we just use edge $(b_0,b_2)\in C''$ to replace $(b_0,b_1)$,$(b_1,b_2)$, see Fig. 11(a).

The problem we must consider is that if $(b_0,b_2)$ was subdivision-ed when constructing $C'_{new}$. A part of (sub-)edges of (the open set of) $(b_0,b_2)$  might be taken by $C'_{new}$. However, this case could not be happened since three points of triangle $b_0$, $b_1$, $b_2$ are all in $C'$ , no path need the middle point of edge $(b_0,b_2)$ to get to any other places in $B$. (Recall the construction process of $C'_{new}$.) Therefore, this case is resolved. \\

If $b_1$ is in $C''$, i.e. $(b_0,b_2)$ is in $C'$, we just use  $(b_0,b_1)$,$(b_1,b_2)$
to replace $(b_0,b_2)$ in $C'_{new}$ to get  $C''_{new}$ (See Fig. 11 (b)). Again, the only problem is that (if or when) a part of edge in $(b_0,b_1]$ and $[b_1,b_2)$  is in
$C'_{new}$.
\newline
We know that $C'$ does not contain $b_1$. So to use $b_1$ (or related subdivision) is to make that $b_1$ will be included in $C'_{new}$ at
some (time) point. As we discussed before, when $b_1$ is available (for $C'_{new}$) , there is no need to subdivision the middle part of edge $(b_0,b_1)$ or $(b_1,b_2)$ in constructing$\/$finding $C'_{new}$. In fact, (if we know $b_1$ is in $C''$) we do not have to use $b_1$ to make $C'_{new}$---there must be another way to do it.

This is because
that we can independently make a $CC''_{new}$ for $C''$. From  $CC''_{new}$, we can make $CC'_{new}$ just like discussed in the above paragraph (Subpart 1.1).
Assume that the arc $\alpha \in C'$ in $D'$ caused to use $b_1$, but this $\alpha$ is exactly the same as it in $C''$ in $D'$. It does not need to use  $b_1$ regarding to $C''$. So we choose a path for $C'_{new}$ that do not use $b_1$. So it is valid.
(To explain more, $b_1$ is in $C''$ and  $\alpha$ is in both $C'$ and $C''$. To replace $\alpha$ in $C''$ must not use $b_1$ if we
make an independent $CC''_{new}$. So we use the one do not use $b_1$ . Based on the only difference of a triangle in $B$ between $C'$ and $C''$. We can choose the path does not need $b_1$ for $C'$ . See, Fig. 11 (b))

Note if two adjacent points are in $C'\cap B$, these two points are also in $C'_{new}\cap B$. We will not change that. There is no
reason to modify them. In other words, an edge in $C'\cap B$ will be kept in $C'_{new}\cap B$. The original arcs of $C'$ in $B$ will be kept in $B$.

The strategy of making $C''_{new}$ (the next new curve in ) is the following:
We do not want to make $C''_{new}$ independently, then do the modification based on $C'_{new}$. We just do it from $C'_{new}$ to $C''_{new}$.
(This sentence means that we do not want to make $C''_{new}$ first, we want to make $C''_{new}$ based on $C'_{new}$. It is step by step, as similar to an iterated process. One by one move and look forward for the part in $D'$. )

We will give more examples of this case later.

\begin{figure}[h]
	\begin{center}
   \epsfxsize=5in
   \epsfbox{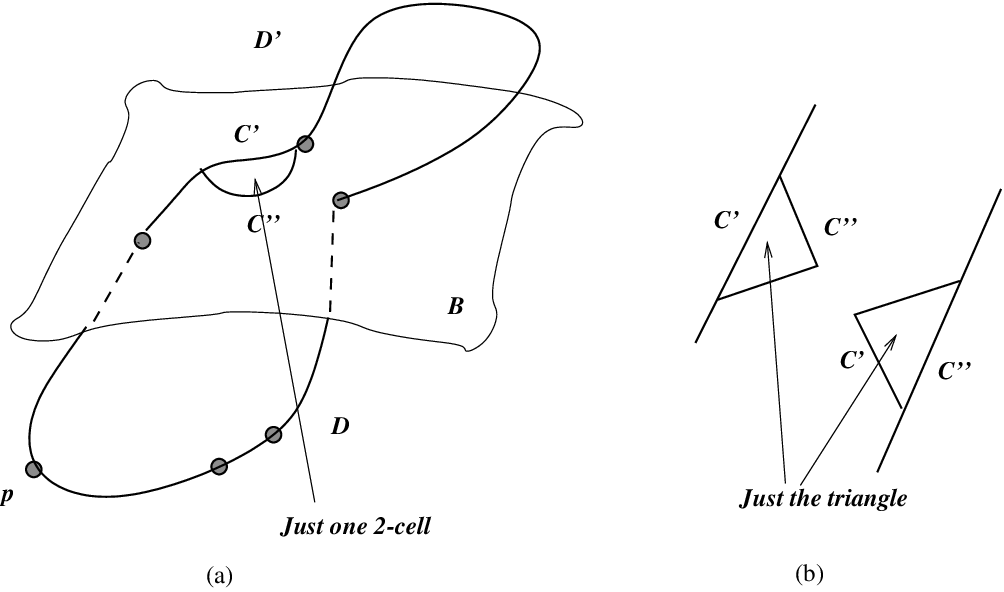}   
\caption{Two adjacent curves in $\Omega$, changes on $B$.  }
\end{center}
\end{figure}

\begin{figure}[h]
	\begin{center}
   \epsfxsize=5in
   \epsfbox{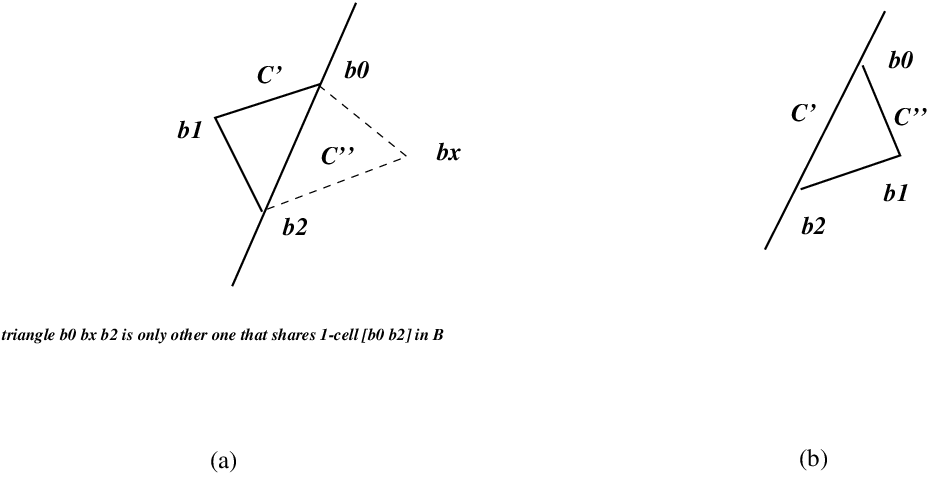}   
\caption{The detailed construction of two adjacent curves in $\Omega$, changes on $B$.  }
\end{center}
\end{figure}

{\bf (1.2.c):} If the triangle $b_0$, $b_1$, $b_2$ has a part in $D$ or $D'$, cannot be in both. This is also based on the general Jordan Theorem again (Theorem 4.2). Or just think about that a triangle 2-cell has three edges and three vertices. See Fig. 12.

Then, (we can assume that) $B$ contains an 1-cell $[b_1, b_2]$, two points in the triangle. Otherwise, if $B$ only contains one 0-cell of the triangle, then no change in $B$ is needed. See Fig.12(b).
So assume  that $b_1$ is in $C'$ and $b_2$ is in $C''$.  If $b_2$ is also in $C'$, (see Fig. 13 (a)) then  $(b_1, b_2)$ is an edge in $C'$. It is also in $C'_{new}$ (as we talk before that a 1-cell in $C'$ will not be changed in $C'_{new}$ based on the construction of $C'_{new}$). Replace arc $b_0 b_1 b_2$ by $b_0 b_2$ in $C'_{new}$ to get $C''_{new}$. See Fig. 13 (b).

If $b_1$ is in $C''$, (see Fig. 14 (a)), $b_2$ could be taken by the process of making $C'_{new}$ as similarly we discussed before in (1.2.b). For the same reason,  $b_0$, $b_1$, $b_2$ are in $C''$. A same arc $\alpha$ (it is in $C'$ and it takes $b_2$ in making $C'_{new}$)  are in both $C'$ and $C''$, there must be a path do not need to pass $b_2$ regarding to $C''$ (in making an independent $C''_{new}$). Use it for $C''$ to get the new $C'_{new}$. The new $C'_{new}$ does not contain $b_2$. (Note that the only difference between $C'$ and $C''$ is the boundary of the triangle $b_0$, $b_1$, and $b_2$. In real algorithm design and implementation, we actually did not go back to change it---the corresponding ``$\beta$'' to $\alpha$ for $C''_{new}$, we look forward to find this triangle in $C'$ and $C''$ before we make $C'_new$. We do this to prevent this case, i.e. $b_2$ could be taken, to be happening.)  Therefore, we can use $(b_1, b_2)$ and $(b_2, b_0)$ to replace $(b_1, b_0)$ in $C'_{new}$ to obtain $C''_{new}$.

To summarize, we use the independent process of making $C''_{new}$ first in Subpart (1.2.a), to avoid to use $b_2$ for $C'_{new}$, then use $C'_{new}$ to get $C''_{new}$. The independent process of making $C''_{new}$ are the same as using $C'$ for $C'_{new}$ just before reaching $\alpha$. We get $\beta$ and use it in $C'_{new}$. Then continue. (there might be other $\alpha$'s to be processed before we reach to the triangle $b_0$, $b_1$, $b_2$ that need to be processed.)  See Fig. 14. We have discussed a similar case in (1.2.b).
This discussion shows us there must be one that does not use $b_2$ to get current $C'_{new}$ while gets the projection of $\alpha$ in $B$. We use this one.

(It is important to mention that we only care about current $C'$ and its previous one and its next one in $\Omega$.)

\begin{figure}[h]
\begin{center}
   \epsfxsize=5in
   \epsfbox{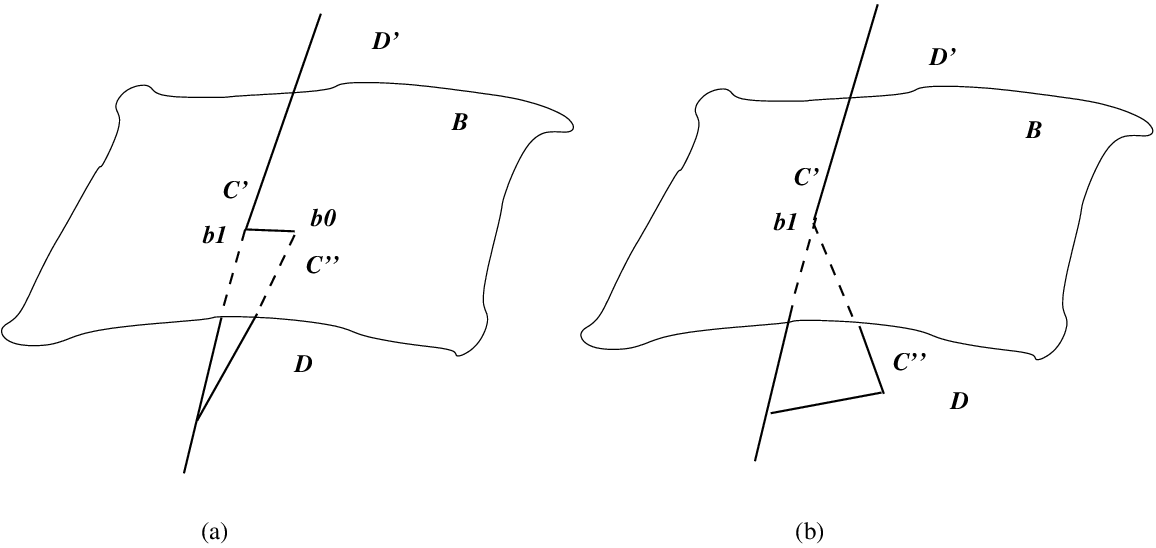}   
\caption{Two adjacent curves in $\Omega$, changes partly in $B$.  }
\end{center}
\end{figure}

\begin{figure}[h]
\begin{center}
   \epsfxsize=5in
   \epsfbox{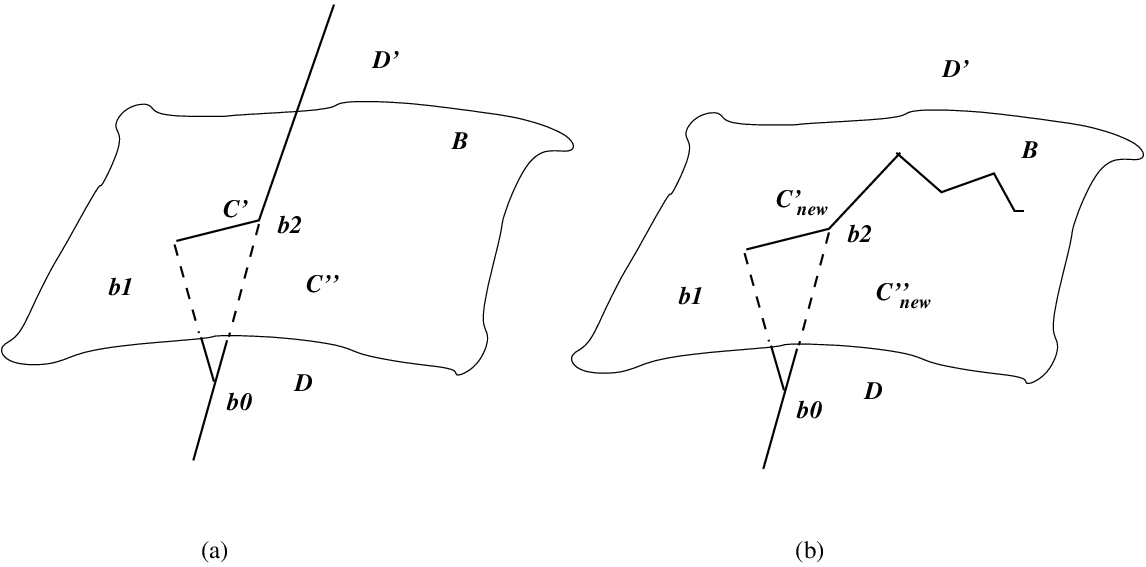}   
\caption{Changes partly in $B$ from $C'_{new}$ to $C''_{new}$ in $D\cup B$.  }
\end{center}
\end{figure}

\begin{figure}[h]
\begin{center}
   \epsfxsize=3in
  \epsfbox{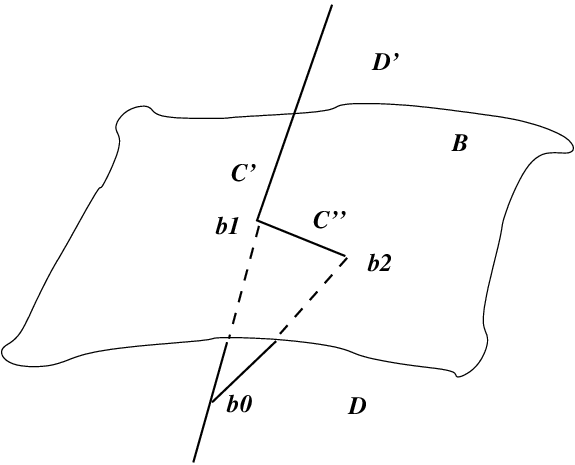}   
\caption{Look forward to get $C''_{new}$ in $D\cup B$ for $C'_{new}$ in order to avoiding to use $b_2$ in $C'_{new}$.  }
\end{center}
\end{figure}

Next, we face the case that the triangle $b_0$, $b_1$, $b_2$ has a part in $D'$ and another part in $B$. This one should be simpler since
we only care the part in $B$ and we will cut out all curves in $D'-B$.  If the triangle $b_0$, $b_1$, $b_2$ has only one point in $B$, $C''_{new}$ will be made to be equal to $C'_{new}$, no need to change it.

Let the triangle $b_0$, $b_1$, $b_2$ have a part $[b_1, b_2]$ in $B$ (See Fig. 15).

{\it Subcase 1:}  We first assume that $b_1$, $b_2$ are in $C'$ (Fig. 15 (a)). Then $[b_1,b_2]$ is in $C'$ and also in $C'_{new}$. $[b_0,b_2]$ is in $C''$. We need to determine if $b_1$ is a leaving point or entering point (from $p$ on counterclockwise) for making $C'_{new}$.

\begin{figure}[h]
\begin{center}
   \epsfxsize=5in
  \epsfbox{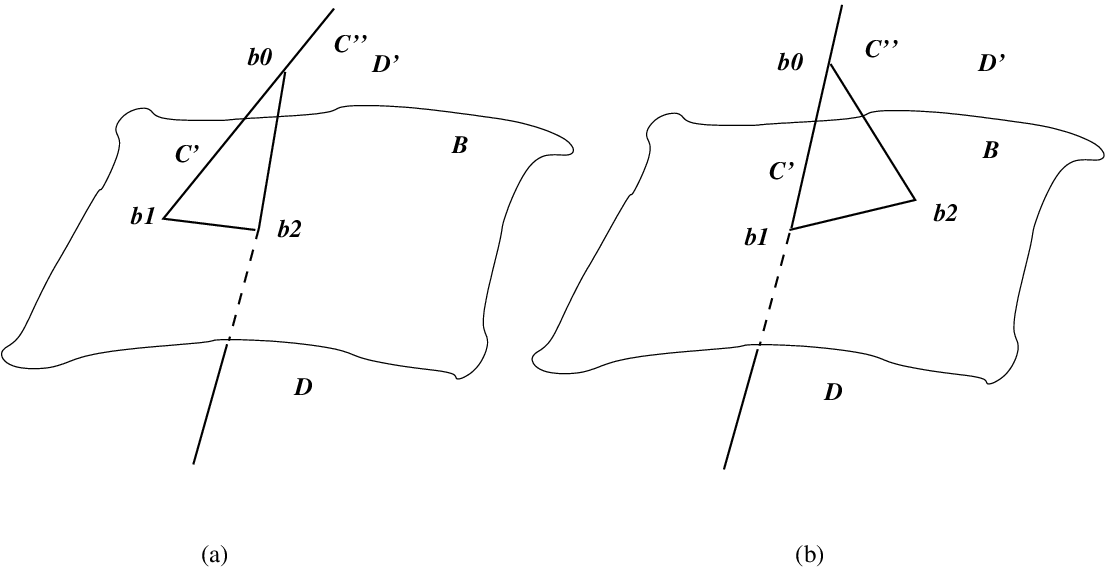}   
\caption{Changes partly in $B$ from $C'$ to $C''$ in $D'\cup B$. (a) Subcase 1, (b) Subcase 2.}
\end{center}
\end{figure}



If $b_1$ is a leaving point to $D'$ and $b_2$ is  in $C''$ (Fig. 16 (a)).  Because both arcs containing $b_0$ in $D'$ will enter the same point $x$ in $B$. $C'_{new}$ contains an arc in $B$ from $b_1$ to $x$, we just use this arc for $C''_{new}$ with the edge $[b_1,b_2]$ that is also in $C'_{new}$, i.e. just let $C''_{new}= C'_{new}$ (at this time). Note that $C''_{new}$ might not have to pass $b_1$, however, pass $b_1$ is fine as long as $C''_{new}$ contains $b_2$.   See Fig. 16 (a).    If $b_1$ is an entering point from $b_0 \in D'$.  So $b_2$ is an entering  point of $C''$ from $b_0$ in $D'$. We can still use the same $C'_{new}$ to be $C''_{new}$ since $b_2$ is also in $C'$ as well as in $C'_{new}$. We let $C''_{new}$ pass $b_1$ is fine in $B$. See Fig. 16 (b).
(It dose not matter to make $b_1$ in the middle of a path.)

\begin{figure}[h]
\begin{center}
   \epsfxsize=6in
  \epsfbox{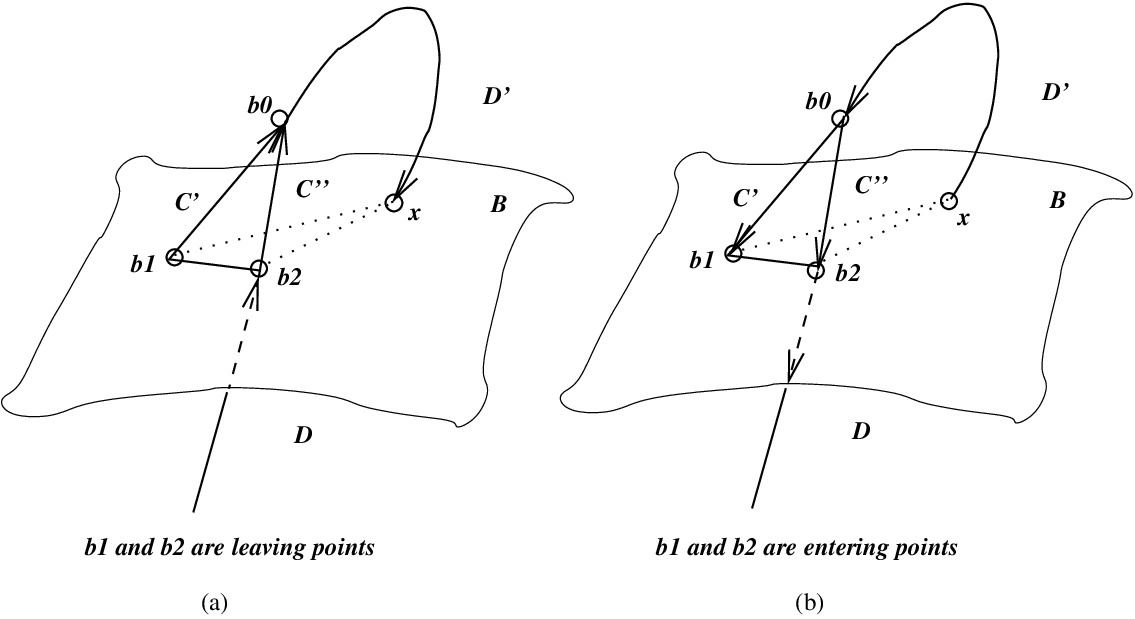}   
\caption{Subcase 1: changing partly in $B$ from $C'$ to $C''$ in $D'\cup B$ while making $C'_{new}$ and $C''_{new}$ . }
\end{center}
\end{figure}

{\it Subcase 2:} Assume that $b_1$, $b_2$ (in $B$) are in $C''$ and only $b_1$ in $C'$ (Fig. 15 (b)). If $b_1$ is a leaving point regarding to $B$ (Fig. 17(a)), then $b_2$ is the leaving point of $C''$.
There is a path in $B$ can pass $b_2$ (and $b_1$)in making $C''_{new}$ because of $C''$ having the edge $[b_1,b_2]$ . Use that one for $C'_{new}$ as well (this path contains the arc from $b_1$, $b_2$, then to $x$ in $B$). This is to say  $C''_{new} = C'_{new}$ here.

If $b_1$ is the entering point of $C'$ and $b_2$ is the entering point of $C''$ (Fig. 17(b)),
$C'_{new}$ can have $b_2$. Plus $[b_2 b_1]$ must be in the $C''_{new}$ as we wanted before (do not change any edge of $C'$ in $B$). Again, if there is an $\alpha$ of $C'$ that could make $C'_{new}$
to take $b_2$. We can look ahead to make a $\beta$ for $C''_{new}$ first and let $C'_{new}$ the same as $C''_{new}$), see Fig. 17 (b).
Thus, we can also make an arc in $B$ that pass $b_2$ to $b_1$ in $C'_{new}$ as well.

(These arrangement are just for the local with regarding to this triangle.
When $C'$, $C''$, and $C'''$ are three consecutive curves having the similar settings, we will use the same strategy to construct $C'''_{new}$ first, and so on so forth. We will give a general consideration in the following and some more examples and figures later.)

\begin{figure}[h]
\begin{center}
   \epsfxsize=6in
  \epsfbox{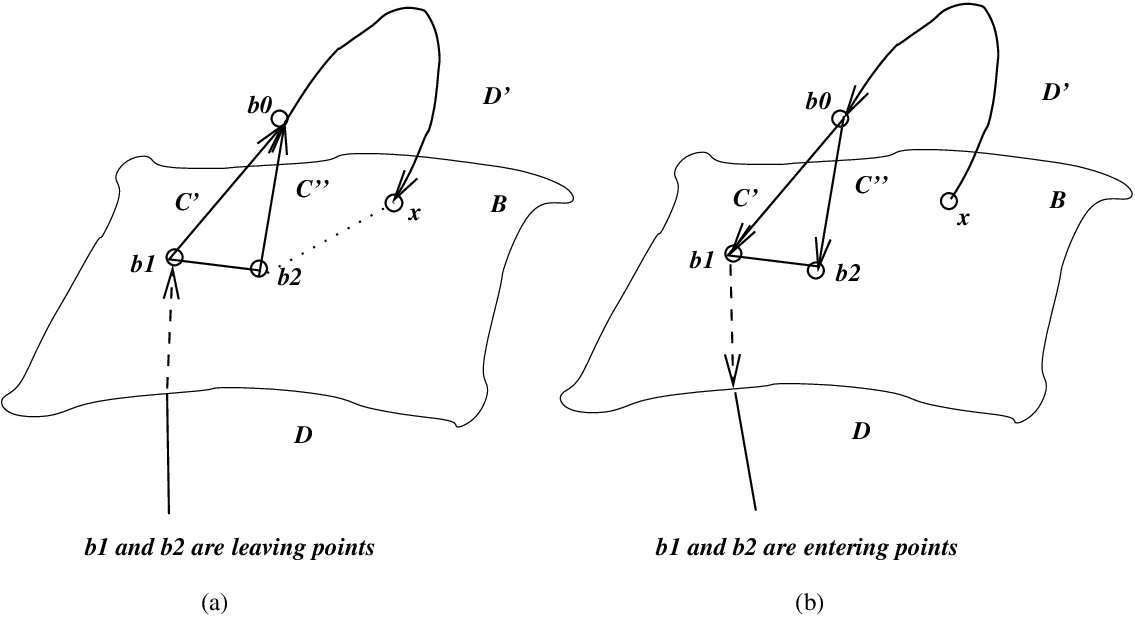}   
\caption{Subcase 2: of changing partly in $B$ from $C'$ to $C''$ in $D'\cup B$ while making $C'_{new}$ from $C''_{new}$: (a) Use the path $b_1\rightarrow b_2 \rightarrow x$ for $C'_{new}$, (b) Let $b_2$ in $C'_{new}$. }
\end{center}
\end{figure}

\noindent Now, we have completed the construction for (1.2.c). \\

Now we recall the cases that need to look forward consecutively. In other words, we here give a general consideration of the consecutive cases below:

In both (1.2.b) and (1.2.c), we meet the cases that $C'_{new}$ may take the point on $C''$ in $B$. And it is possible that this case will consecutively occur. We want to provide a unified solution for this situation. Let us assume that there are $t$ curves $C'$,$C''$, $\cdots$, $C^{(t)}$ in $\Omega$ where $Modulo2 (C^{(i-1)},C^{(i)})$ is the boundary of a triangle $\Delta_i$. The characteristics are:
\newline
1) $C^{(i-1)}$ shares one edge with $\Delta_i$ and $C^{(i)}$ shares two edges with $\Delta_i$.
\newline
2) $C^{(i)}$ has one more 0-cell than  $C^{(i-1)}$; $C^{(i)}$ contains all 0-cells of $C^{(i-1)}$.
\newline
3) $C^{(i)}$ can be viewed as always extending a $\Delta_i$ by adding a new point from $C^{(i-1)}$.

\noindent $t$ is finite.   $C^{(t)}$ will pass (include) all points of  $\Delta_i$, $i=2,\cdots, t$. Also, $C^{(t)} - C'$ contains the union of all 0-cells of  $\Delta_i$, $i=2,\cdots, t.$  If any of points in $C^{(t)}$ in $B$ was used in $C'_{new}$ , it must be in $\Delta_i$, $i=2,\cdots, t.$  So if we use $C^{(t)}$ to determine $\beta$ for that particular $\alpha$ mentioned in (1.2.b) and (1.2.c), none of points $\Delta_i$, $i=2,\cdots, t.$,  will be included in  $C'_{new}$. So this $\beta$ in $B$ can (will) be used for $C'$ for getting $C'_{new}$. Also $\beta$ in $B$ will be used all $C^{(i)}$, $i=2,\cdots,t$. So this issue is totally resolved (proved). See Fig. 20. \\

After all three cases considered, we can construct a new contraction sequence without using any points in $D'$.
It is a kind of projecting method to project the part (arcs) in $D'$ to $B$ for each curve ($C'$, made by 1-cells) and keeping the new curves
to be simple curves.

We repeat the above process for all simple closed curves (1-cycles) in the contraction sequence $\Omega$, we will get a new gradually varied sequence $\Omega_{new}$ in $D\cup B$. This means that the changes were made along with the boundary of only one 2-cell at a time for two adjacent new curves in $D\cup B$. Therefore, $D\cup B$ is simply connected.   In the same way, we can prove that $D'\cup B$ is simply connected too.

Note that the curves in the contraction sequence $\Omega$ in $D-B$ are not changed if the curves are completely in $D-B$. However, since we may do subdivision of $B$, so we have to subdivide some
3-cells (or $m$-cell) intersecting $B$. As we discussed before, we can find some gradually varied curves inserting to some curves in original $\Omega$. The inserted subsequence only have changes inside of a 3-cell (or $m$-cell).

In addition, for any curve in $\Omega$, the arcs inside of $D-B$ are not changed even though we may change its part in $D'\cup B$. In above proof, we projected the arcs of $\Omega$ in $D'\cup B$ to $B$ without moving points intersecting with $B$.  We also make sure that they are gradually varied in $B$. $\Omega$ as a sequence could travel to $D'$ multiple times from $D$, we proved for the single travel. For the case with multiple travels, it is the same. So we complete the proof for Part 1.\\

\noindent We can summarize the entire proof of {\bf Part 1}  as follows: Our goal is to prove that if $M$ is simply connected so are $D$ and $D'$.

(1) First, assume a point $p$ in a simple 1-cycle $C$, made by 1-cells in $K$. We might as well to assume that $C$ is not touch $B=\partial D$. We want to prove that if there is a contracting sequence $\Omega$ for C in $M$ , then there will be one in $D$ (or $D\cup B$). A step of contraction of a curve in finite sense can be viewed as to move a curve from one side of a cell in $K$ to another side of this cell. Note that 1-cycle $C$, defined that $C$ is homeomorphic to a 1-sphere, is equivalent to a simply connected simple path of 1-cells.  (This statement is true for 2-cycles that is equivalent to
simply connected 2-manifold, we are especially prove this statement for 3-manifolds in this paper.)

We only need to consider the case that some cycles of the contraction sequence $\Omega$ will go to $D'$ then go back to $D$ then shrink to one point $p\in C$. Let $\Omega = \{ C_1, C_2, \cdots, C_N \}$ be such a sequence. We know that every $C_i$, $i=1,\cdots,N$, contains $p$. ($C_N$ can be the boundary of the last 2-cell containing $p$.)

(2) We might as well to assume here that this sequence starting in $D$ only goes $D'$ once than goes back to $D$. (The sequence that  travels multiple times from $D$ to $D'$ will be similar. We just need to repeat the following procedure.) There must be an $k$ so that $C_k\in \Omega$ is the last 1-cycle in the sequence that has intersection with $B$. $C_{k+1}, \cdots, C_N$ are totally in $D-B$.

(3) Assume also $C_{t+1}, \cdots, C_k$ have 0-cells in $B$, and $C_{t}$ is the last one that does not have any intersection with $B$.
Note that $C_{t+1}$ is gradually varied to $C_{t}$. $C_{k}$ is gradually varied to $C_{k+1}$. They all contain $p$.

(4) So there are three parts in $\Omega$: $\{ C_1, C_2, \cdots, C_t \}\in (D-B)$ , $\{ C_{t+1}, \cdots, C_k \}\in M)$ , and $\{C_{k+1}, \cdots, C_N\}\in (D-B)$.  It is not hard to know that $C_{t+1}$ and $C_k$ are in $D \cup B$. These two are not be changed. The modifications
were made only for $C_{t+2}, \cdots, C_{k-1}$ in the above proof.

(5) No entering or leaving points between $D-B$ and $B$ are changed in $\Omega$ . The process (of projection)  will guarantee the arcs in $B$ are gradually varied in the new sequence. The new contraction sequence will be totally in $D\cup B$. So we proved that $D$ is simply connected. Again, the concept of gradually varied sequence of curves means that each adjacent pair of curves only differs a triangle (or a 2-cell) in $K$, i.e. $Modulo(2)$ of two such curves will be the boundary of a 2-cell (triangle). A step of contraction is a move of gradual variation between two adjacent curves in the sequence. \\

\noindent {\bf Part 2:} We design an algorithm to show that $D$ is homeomorphic to an $m$-disk as well as $M$ is homeomorphic to $m$-sphere if $m=3$ or greater.

First, according to the definition of simple connectedness to $M$, a closed 1-cycle will be contractible on $M$. Since we assume that
there is no boundary in $M$ ($M$ is closed and orientable), then there are no holes in $M$. To understand this, we can use examples of 2D manifolds.  If $M$ is a 2D manifold and there is a hole in $M$,
then the boundary of a hole will be a 1-cycle and it is no longer contractible. For another instance, if $M$ is a 2D torus, then some cycles are not contractible as well.

Second, any closed and simply connected manifold (orientable) is a 2-cycle meaning that this manifold is homeomorphic to a 2-sphere. In this part of the proof, we will keep the boundary of $D$ to be a 2-cycle.

Because we only have finite numbers of cells, we can calculate the cell-distance for each pair of cells. We will select two special points (0-cells) they have the far-most distance among each pairs in the manifolds.
one of the point will be in a 3-cell $e$ that will be removed first. Another one will be the origin $O$ that will remain in the contracting process.


Let $m=3$ for now, i.e., $M$ is a 3-manifold. We will see that $m$ can be any number.

\noindent {\bf Step (1):}

Remove an $m$-cell $e$ from $M_m=M$. This will result in an $(m-1)$-cycle on the boundary of $e$ and $M-e$ called $B$. This cycle $B$ is always simply connected since its boundary is also the boundary of $e$, a single
cell. Let us call this cell is $e_{0}$. See Fig. 18. (a). We will use this property in addition to the theorems and properties we presented above including the Part 1 to prove  this second part.

Note: As for using the mathematical induction, we can assume that any closed simply connected (and orientable) $(m-1)$-manifold is an $(m-1)$-cycle. It means  that it is homeomorphic to the $(m-1)$-sphere when proving for $m$.
This is the condition of the mathematical induction (the induction hypotheses):  Assume that when $k=m-1$, the statement is true, then prove it for $k=m$.

Let $D=M_m-e$. Note that the subtraction is to remove the inner part of the $m$-cell $e$, and not the $(m-1)$-faces of $e$. $B=\partial e$ at this time (moment).

\begin{figure}[h]
\begin{center}
   \epsfxsize=6in
  \epsfbox{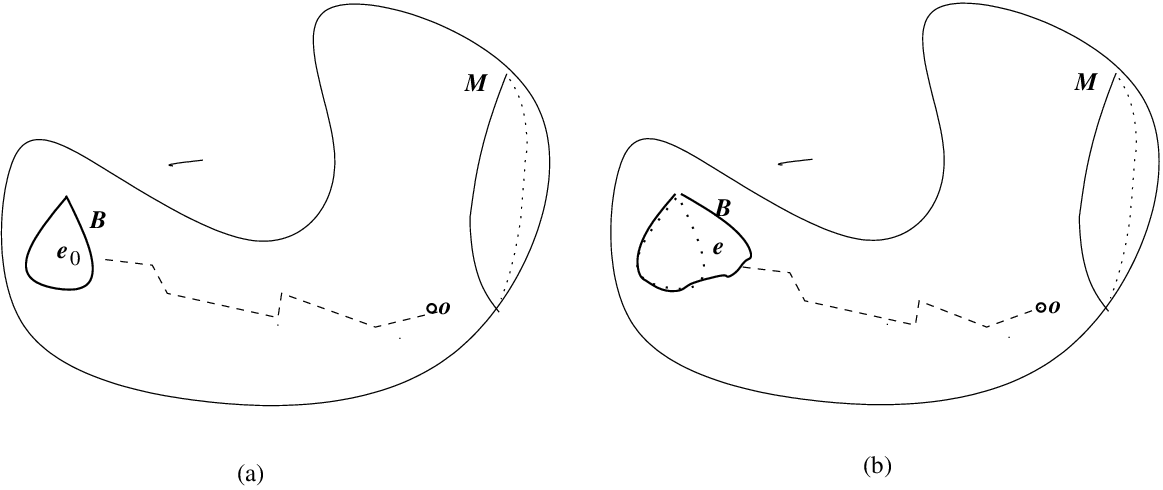}   
\caption{Remove a 3-cell $e$ from $M$ to get $D$ and then contract $D$ to point $o$. $M$ only contains a finite number of 3-cells. In this figure, (a) remove $e$ by leaving all boundary cells in $D=M-\{e\}$ ,
(b) contract a new $e$ by using part of boundary of $e$ to replace the boundary on old $B$ }
\end{center}
\end{figure}

\noindent {\bf Step (2):}

Next, we remove more $m$-cells of $M_m$ ($D$) if they have an edge (face in 3D) on $B$. This type of removing is contraction called the contracting removal. The operation is to use the part of boundary of a new $m$-cell $e$ in $D$ to replace its boundary on old $B$ to form a new $B$. See Fig. 18. (b). We will prove that $D$ is contractible to a $m$-cell in Part 2. Then we can see that attaching the boundaries of 2 $m$-cells will be an $m$-sphere.

A contracting removal of an $m$-cell $e$ is to use a part of the boundary of $e$ to replace another part of $e$. The union of these two parts are
the total boundary of $e$; the intersection of these two parts is a 1-cycle when $m=3$. When $e$ is a 3-simplex, any partition of
the boundary of $e$ is 1-connected (sharing a 1-cell, not only a point or 0-cell). The boundary of the partition is 1-cycle. See Fig. 19.

In other words, we can first assume $e\cap B$ is a 2-cell, or two 2-cells, or three 2-cells.  Contracting removal of $e$ is to use
$\partial e - e\cap B$ to replace $e\cap B$ in $B$ to get a new $B$. The boundary of $e\cap B$ is the same as the boundary of
$\partial e - e\cap B$ that is a 1-cycle. (This cycle separates $\partial e$ into two components as the Jordarn curve theorem says.
If we view each component is an open set, the two components are not connected when e is embedded in Euclidean space.)

\begin{figure}[h]
\begin{center}
   \epsfxsize=3.5in
  \epsfbox{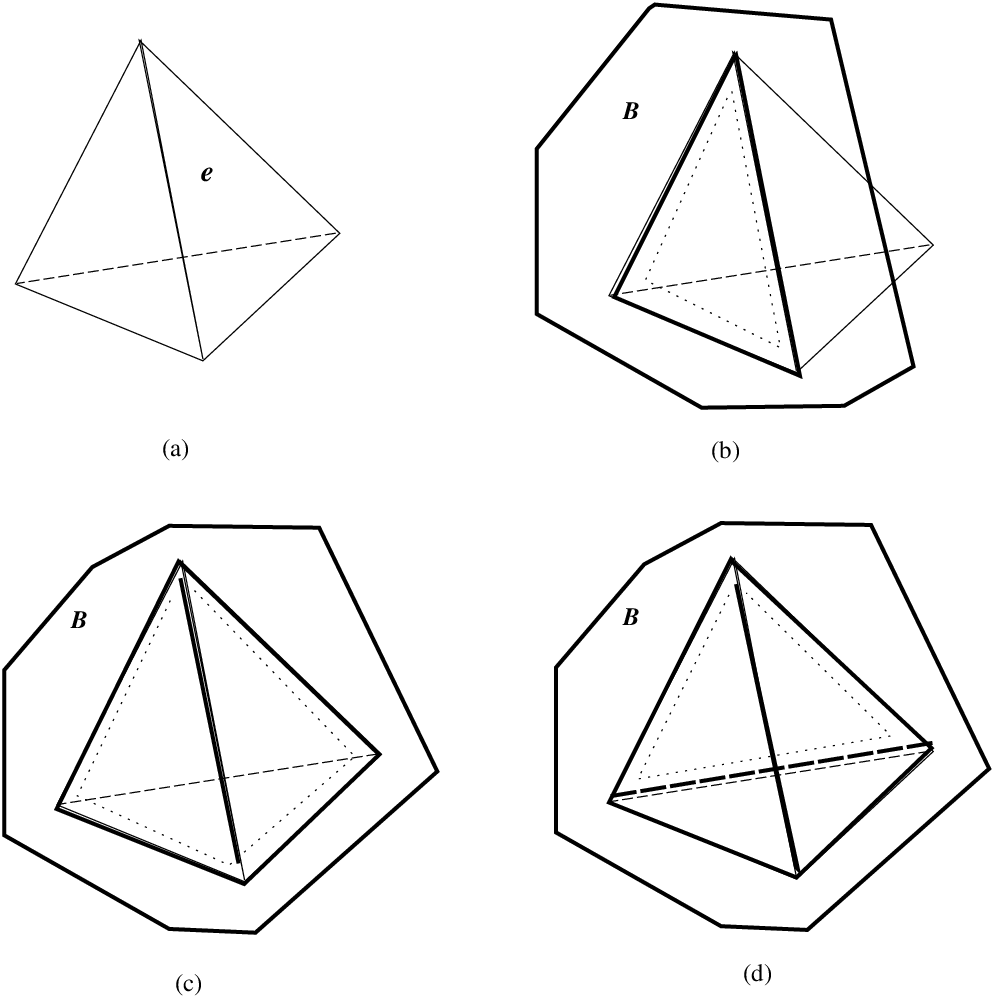}   
\caption{Different kinds of partitions of the boundary of the 3-simplex $e$ on $B$, i.e. $e\cap B$: (a) A 3-simplex (3-cell) $e$, (b) a 2-face of $e$ on $B$, (c) two 2-faces on $B$, and (d) three 2-faces on $B$.
The dotted line-segments represent the boundary of $e\cap B$.}
\end{center}
\end{figure}

Algorithmically, since $B$ is changing as time goes on, we remove another $e\in D$ when $e$ has an $(m-1)$-face (cell) on $B$ (at least one on $B$).
We always denote the new edge (boundary) of $D$ as $B$. For detailed algorithm design, we need to first design a data structure to save $K$ as a linked list. We denote $K_D$ representing the
complex of $D$.  So remove an $e$ is to delete the simplex $e$ from $K_D$. So it is actually programable in computers. We will present the algorithm steps at the end of the proof of Part 2.

Therefore, $B$ will be a new boundary face set of $D=M-\{\mbox{removed $m$-cells}\}$.
$B$ is still the $(m-1)$-cycle (a pseudo-manifold in \cite{Che14}) since we use the boundary of the new $e$ in $D-B$ to replace $e\cap B$ .

Let us look at the Fig. 19 once again:
Basically, we can assume $e$ is always a 3-simplex (or $m$-simplex). Removing $e$ from $D$ is to use a part of boundary of $e$ to replace $e\cap B$. This part is $\partial{e} - e\cap B$ that contains all (boundary) faces in $e$ not in $B$. There are three cases:

(a) $e\cap B$ is just one 2-cell. See Fig.19 (b). Its boundary (of this 2-cell) is a simple closed curve which is also the boundary of three other 2-cells of $\partial {e}$. Use the three 2-cells to replace  $e\cap B$ in $B$ will result also a 2-cycle. Assume the original one (the old $B$) is a 2-cycle,
a closed simple curve (the boundary of $e\cap B$)  will separate original $B$ into two parts. According to the Jordan curve theorem and Jordan-Schoenflies theorem, each part is homeomorphic to a 2-disk. After replacing, the part which left in $B$ stays the same. The three 2-cells is also homeomorphic to a 2-disk. So the connected-sum of two parts (two 2-disks) is homeomorphic to 2-sphere. So this new $B$ is a 2-cycle. Note that, for higher dimensions, we can use the $(m-1)$-dimensional Jordan theorem and Jordan-Schoenflies theorem as the condition of the mathematical induction--the induction hypothesis.

(b) $e\cap B$ contains exact two 2-cells. They must share a 1-cell. So the boundary of these two cells is still a simple closed curve. Removing $e$ means to use other two 2-cells to replace the previous two ($e\cap B$) in $B$. $\partial{e} - e\cap B$ is other two 2-cells where these two 2-cells still shares a 1-cell. The boundary of these two 2-cells is a 1-cycle.   For the same reason given in (a), the new $B$ is a 2-cycle as well.

(c) $e\cap B$ contains exact three 2-cells. These three are three 2-faces in $e$, every pair of two 2-cells shares a 1-cell. The boundary of the three 2-cells is the boundary of one 2-cell. Use this 2-cell to replace the previous three 2-cells of $B$ will result a new $B$. This case is just like (a) reversed in proving the new $B$ is a 2-cycle. We will always preserve that the new $B$ is 2-cycle.

For the case of an (m-1)-cycle, $m>3$, as long as we use simplex, the prove will be the same. For general PL $m$-cells (or the cell we restricted in \cite{Chen-Krantz}), the method is also similar.
The key is to just keep removing the $m$-cell $e'$ so that the boundary $B \cap e'$ is homeomorphic to $(m-1)$-disk. But for simplicity, we just deal $m$-simplex now.

Based on Part 1 of this proof, the new $D$ by removing $e$ from old $D$ is simply connected. The new $B$ is always the boundary of the new $D$ .


For the actual design of the algorithm for contracting removal(for Step 2), we first calculate the cell-distances to all points in $M_m$ from a fixed point $o$,
$o \notin e$ which is determined beforehand. We can also assume $e$ is the one who has the longest distance to $o$. Or we first find the pair of two points that has the longest distance in $M$. Then we determine that one is $o$ and the another point for find an $e$ containing that point. As we did in Step 1. We let $e_0=e$, and remove inner part of $e$ from $M$ and we got $D=M-e_0$.

We always select a new 3-cell (or $m$-cell) $e$ that is adjacent to $B$ (meaning that it has a 2-face on $B$) with the greatest cell-distance to $O$. This means that one point (0-cell) in the new $e$ has the greatest $m$-cell-distance to $o$. When more than one such cells are found, we will select
the (new) cell $e$ ($e\cap B$ is an $(m-1)$-cell) that has the most of the greatest distance points (0-cells) to $o$. This is a strategy for balancing the selection of the new $e$. We will always attempt to make cycle $B$ the equal distance, in terms of each point on $B$, with respect to $o$. In other words, each cell on $B$ should have almost equal cell-distance to $O$ when it is possible. It is easy to know that $D$ and $D-e$ are homeomorphic.

Note that, for simplexes, it is easy to know that the cell-distance is the same as the edge-distance (or just the distance in graph theory).

In some case, it is possible that the cell-distance of some cells on the new boundary $B$  could be increase when remove a 3-cell ($m$-cell). It does not matter, we still keep to remove the furthermost $m$-cell on the new boundary based on this distance rule.
There are only finite many of cells. So the process will be contracting to a cell containing $O$ finally. As we know, for this contraction, there is no need to require a fixed point $O$. As long as the contraction will be ended at one 3-cell (or $m$-cell), this process is correct. Shrink a 3-cell ($m$-cell) to a point (0-cell) is simple. In the actual algorithm, we will be ending at $Star(o)$ in $M$ since we remove all $m$-cells that has a distance larger than 2 from point $O$.  $Star(O)$ in $M$ is homeomorphic to an $m$-disk,  and $Star(O)$ is very easy to be contract to an $m$-cell.

In fact, there is no need to keep a $Star(O)$ in $D$ for our theorem proof, just keep one 3-cell is also fine in the contraction process. But in actual procedure or computer programming, we could keep $Star(O)$.

There are some extreme cases
that need to consider the order of contractions for 3-cells ($m$-cells) in $D$. Specifically, we could not contract a 3-cell $e$ when

{\it Case 1}: A 2-face ($f$) of $e$ is on $B$ and any of the three other 2-faces is not on $B$, but a 0-cell $P$ not in $f$ is on $B$.

{\it Case 2}: A 2-face ($f$) of $e$ is on $B$ and any of the three other 2-faces is not on $B$, but a 0-cell $P$ and a 1-cell (containing $P$) not in $f$ is on $B$.

\noindent It is easy to see that contracting this $e$ could make $B$ that is not a 2-manifold. In other words,  we just could not contract this 3-cell $e$ at this moment. We need to remove its neighboring cells before contracting this one ($e$).  This is because in such two cases, all 0-cells of $e$ are in $B$, $e$ only contain four 0-cells if $m=3$, the boundary of $e\cap B$ is not 1-cycle. So directly remove $e$, will make new boundary not a 2-cycle.

Note that these two cases usually not happen, but if it happens we will do the following special process to find $e$'s neighbor or its neighbor's neighbor to remove it beforehand.


In any of the two cases, we can use the following process: We can first find a 3-cell $e'$ that has two or more 2-faces on $B$; denote it as $e'$ (i.e. we contract $e'$ first).  We will give a direct and complete proof in the {\bf Appendix A} for the following theorem: Case 1 or Case 2 cannot be true for all 2-cells on $B$. Some discussions will be given in the next subsection where we tried to have a simple proof this theorem.

Note that it is still possible to have few inner points in $D$ even each forth 0-cell of $e$ for all $e$ adjacent to $B$ ($B\cap e$ is a 2-cell) is on $B$.

In summary, remove $e$ that has two or more 2-faces on $B$.

\noindent {\bf Step (3):}

Since $B$ is an $(m-1)$-cycle, $D$ is always simply connected based on Theorem 4.2. $M_m$ only contains a finite number of $m$-cells, which means that this process will end at $Star(o)$ in $M_m$. Therefore, we algorithmically show that $D=M_m-e$, $D\leftarrow D-e_{new}$, $\dots$ is continuously shrinking (homeomorphic) to $Star(o)\cap M_m$, since $Star(o)\cap M_m$ is an $m$-disk.  Therefore, reversing the steps determines a homeomorphic mapping from an $m$-disk to $M_m-e_0$ where $e_0$ is the first cell removed from $M_m$.

Thus, both $e_0$ and $M_m-e_0$ are homeomorphic to an $m$-cell. Then $M_m$ is homeomorphic to an $m$-sphere ($m=3$ or $m>3$) since $M_m$ is homeomorphic to the manifold that is formed by attaching the boundaries of two $m$-cells. Therefore, we complete the proof of Part 2.

                                           $\hfill$ $\square$ \\

\subsection{Is Possible to Have a Simple Proof of Theorem A in Appendix A?}

Assume that the contractional-removal have removed many 3-cells ($m$-cells) until the following two conditions hold for all 2-cells on the 2-cycle $B$ that is the boundary of $D$ in the proof of Part 2 of Theorem 4.3:

For each 3-cell $e$ having a 2-face $f$ on $B$:

{\it Case 1}:  The 0-cell $P$ not in $f$ is on $B$. Or

{\it Case 2}:  The 0-cell $P$ and a 1-cell (containing $P$) not in $f$ is on $B$.

We want to make sure the following statement is true: There must be a 3-cell has two 2-cells on $B$.

We now give some validations here. The following discussions are helpful to understand above question.

Since Case 2 is a special case of Case 1, so we could only discuss Case 1. We know that the Euler Characteristic of a 2-cycle $B$ is 2. Euler characteristic of $S^m$ is $1+(-1)^ {m}$. We could not use the Euler Characteristic for $D$ since we want to prove $D$ is contractible.

For an $m$-dimensional simplicial complex, the Euler characteristic is the alternating sum

$ \chi =k_{0}-k_{1}+k_{2}-k_{3}+\cdots = \Sigma_{i=0}^{m} (-1)^i k_{i} , $

\noindent where $k_i$ denotes the number of $i$-simplexes in the $m$-complex.

This simple validation  uses the mathematical induction (This simple discussion is not a complete proof it is a validation. This is because that we assumed that inserting a 3-cell or removing a 3-cell from middle of $D$ will only affect cells locally, one could argue that we assume something to prove something. But here is a good observation at least.):

(i) First we check $D$ just have a 3-cell. This cell has two 2-faces on $B$. And $B$ is a 2-cycle. Let us check the  $\chi$. Since $k_0=4$, $k_1=6$,$k_2=4$,   $\chi =k_{0}-k_{1}+k_{2}=2$.

(ii) Assume that if (any) $D$ contains only $t-1$ (or less) 3-cells, there must be a 3-cell having two 2-faces on the boundary $B$ of $D$. Note that $B$ is 2-cycle as assumed too.

(iii) We now prove that if $D$ contains $t$ 3-cells, under the condition that $\partial D = B$ is a 2-cycle, $D$ must also contain a 3-cell having  two 2-faces on the boundary $B$ of $D$.

Note that, each $m$ simplicial complex with $t$ $m$-cells can be made by adding an $m$-cell and some $i$-cells, $ 0\le i<m $, to some $m$ simplicial complex with $(t-1)$ $m$-cells. Therefore,
 the new $D$ that has $t$ cells can be made by adding a 3-cell $\Delta_3$ to some (previous) $D$ with $(t-1)$ 3-cells to make the current $D$ that has  $t$ 3-cells.
{\footnote {We can even prove the statement directly using mathematical induction: (1) if ($D$ has $t$ 3-cells) the case 1 is true, and no case 2 appearance, then we can find a 1-cell simple closed-path (a 1-cycle) on $B$ that that splits $B$ into two parts, then determine a surface no 2-cells in $B$ that is bounded by the 1-cycle. Split $D$ into two parts, each split-part is a 3-manifold with fewer 3-cells. Each part of $D$ with less than $t$ 3-cells is contractable.  (2) If there is a Case 2, we just remove $Star_M(Q)\cap D$ where $Q$ is a point on the 2-cell and the 1-cell that is not in the 2-cell on $B$.  Using the mathematical induction hypothesis, each such 3-manifold $D$ with less than $t$ 3-cells is contractable. So, we did our proof.}}

We know that to make new $D$ to be 3-manifold, we must attach at least one 2-face of $\Delta_3$ to $D$. So we have four cases: (a) $D\cap \Delta_3$ is one 2-face.
(b) $D\cap \Delta_3$ is two 2-faces. (c) $D\cap \Delta_3$ is three 2-faces. and (d) $D\cap \Delta_3$ is four 2-faces. For (a) and (b), $\Delta_3$ will have two 2-faces on $B$ after attaching, so we've done these two cases.
For (d), if we insert  $\Delta_3$ inside of old $D$, if $\Delta$ was the 3-cell in old $D$ that has two 2-faces on $B$. It will not change the fact. If $\Delta_3$ is to fill the blank space of $D$, than $B$ of old $D$ will have two disconnected 2-cycles or new $D$ does not have boundary. So they are impossible. So the only case left here is (c), we still have two possibilities: one is to subdivide a 3-cell and another one is to fill.
This action will not affect to $\Delta$ that already has two 2-faces on $B$ if the subdivision do not generate an internal point in $D-B$. It will not change the fact for  $\Delta$. We know it is usually not easy to make the a simple subdivision for just adding one 3-cell generally. It is not our job here, we just assume here if there is a possible subdivision. Then it does not change the fact of having two 2-faces on $B$. We do not have to allow subdivision. It is just a nice action in practice.

Now we check the case for filling, we use $k_i^{t-1}$  for old $D$ and  $k_i^{t}$ for new $D$ for the number of $i$-simplexes in $D$.
$k_1^{(t)}=k_1^{(t-1)}-3$ since three 1-cells on three 2-faces of $\Delta_3$ will be not on new $B$.  We have $k_2^{(t)}=k_2^{(t-1)}-3+1$  since three 2-faces dismissed and one new. For 0-cells on $B$,
(c1) $k_0^{(t)}=k_0^{(t-1)}$ if every point on $B$ in new $D$, but (c2) $k_0^{(t)}=k_0^{(t-1)}-1$ if the forth point is in new $D-B$ (The forth point becomes an internal point. $\Delta_3$ has a 2-face on $B$, its three vertices are on $B$ at least).

(c1) $\chi(new B)=2= k_{0}^t-k_{1}^t+k_{2}^t = k_0^{t-1} -(k_1^{(t-1)}-3)+ (k_2^{(t-1)}-3+1) = k_0^{(t-1)} -(k_1^{(t-1)}-3)+ (k_2^{t-1}-2) =2+1=3$ . It is not possible since this case cannot be made for old $D$. The reason was the old $D$ have $t-1$ 3-cells is not a 3-manifold. The boundary of old $D$ is not a 2-cycle. There is a singularity point.
(c2) $\chi_(new B)=2= k_{0}^t-1-k_{1}^t+k_{2}^t = k_0^{t-1}-1 -(k_1^{t-1}-3)+ (k_2^{t-1}-3+1) = k_0^{t-1}-1 -(k_1^{t-1}-3)+ (k_2^{t-1}-2) =2=2$.  This is to make an internal point, i.e. the point removed from $B$ is in new $D-B$. There is a 2-cell on $B$ that is in a 3-cell
containing a point in $D-B$. This is not the Case 1 and Case 2. Therefore, we can still contract-removal 3-cells from $D$ to maintain $B$ is a 2-cycle.

For higher dimension, we can do the same because we only need to have two $(m-1)$-faces to have all 0-cells covered in an $m$-simplex. An $m$-simplex $\Delta_{m}$ contains $(m+1)$ $(m-1)$-faces; if there are $m$ $(m-1)$-faces will be attached to old $D$, we still have two subcases: (c1) $k_0^{t}=k_0^{t-1}$ if every point on $B$ in new $D$, but (c2) $k_0^{t}=k_0^{t-1}-1$ if the $(m+1)$-th point is in new $D-B$ (The $(m+1)$-th point becomes an internal point. $\Delta_m$ has a $(m-1)$-face on $B$, its $m$ vertices are on $B$ at least).  $\chi_{new B}=1+(-1)^{(m-1)}= \Sigma_{i=1}^{(m-1)} (-1)^{(i)} k_i^{(t)}$ . Thinking about $\Delta_m$, when it is filled in to old $D$, only its boundary affect the calculation of new $B$. The new $B$ is the connected sum of the old $B$ (the boundary of $t-1$ $m$-cells) and $\partial \Delta_m$. According to the formula of the connected sum of Euler Characteristic for compact manifolds, we will have  $\chi(new B) = \chi(old B) + \chi(\partial \Delta_m)- \chi(S^{(m-1)})$ since $B$ is $(m-1)$-dimensional {\footnote {This formula is related to connected sum. It can be induced from the
 popular theorem: If $X$ is the union of two closed sets $S1$ and $S2$, then
 $\chi(X) = \chi(S1) +\chi(S2) - \chi(S1\cap S2)$. See \cite{Hat} on page 157. We assume here the connected-sum $\#$ will only for cutting out a disk. One open disk and one unopen disk glue together will be a $S^n$ without intersection. Therefore,
 $\chi(new B) = \chi(old B) + \chi(\partial \Delta_m) - \chi(S^{(m-1)})$. See also  L.I. Nicolaescu, “The Euler Characteristic.” $https://pdfs.semanticscholar.org/42dd/bacb57f8647febb130f3adaa81b48fb1a993.pdf$ }}. This will make the $(m+1)$-th 0-cell an internal point in $(D-B)$. It is the case of (c2). The case (c1) will make old $B$ not an $(m-1)$-manifold.   So we proved the result for all dimensions.

For lower dimension, it is the same that is even easier. We will discuss it in the Appendix A.

Note that:  It is important that we only allow subdivision in this proof. We will not allow subdivision when we do the actual contraction-removal procedure or computer program. Subdivision without making an internal point, will not change the fact of having two 2-faces on $B$ for a 3-cell. We do not have to allow a subdivision. As a basic fact, any $t$ cells manifold is made by a $(t-1)$ cell manifold by connected-sum (attaching or filling) an $m$-cell. We will give examples later.

We will give a more direct proof to that this $e$ always exists adjacent to $B$ in the Appendix A for Case 1 and Case 2.

\subsection{Some Examples Related to the Proof of Theorem 4.3}

In the above proof, we did not include some examples that contain some extreme cases. Now we add them here. They might be a little tedious for some readers.

{\bf Example 4.1}\\
First, the example of subdivision of $B$ as needed. This example was shown in a separated paper for the proof of Part 1 ~\cite{Che18SummerTop}.
As an example of projecting an arc from $D'$ to $B$ when most of 0-cell are already taken by $C'$, in order to make $C'_{new}$, we may need to so subdivision of some triangles for $B$ in Euclidean space (also need to check the subdivision satisfying the condition of PL-manifolds.). See Fig. 18. \\

\begin{figure}[h]
\begin{center}
   \epsfxsize=6in
  \epsfbox{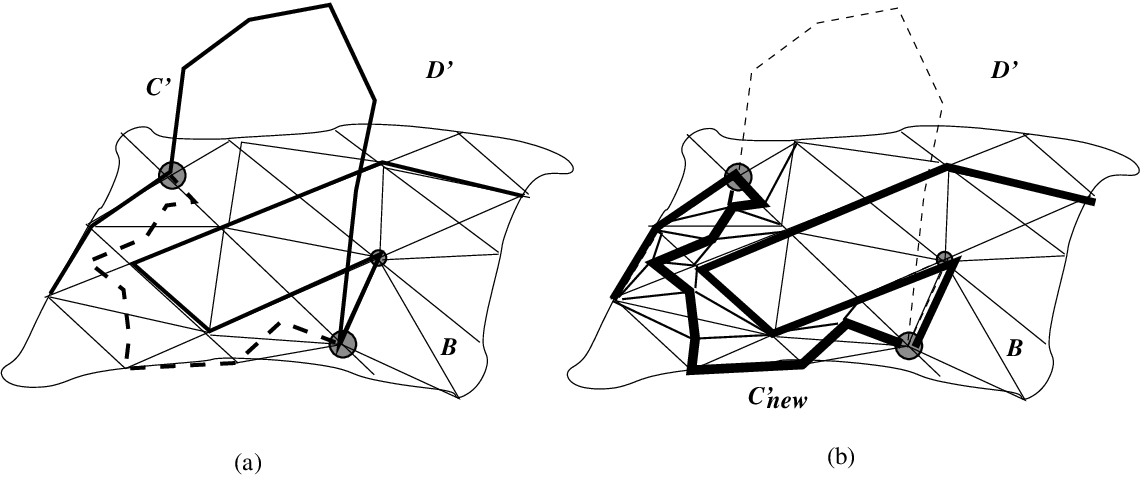}   
\caption{ Example of projecting an arc from $D'$ to $B$ if the subdivision is needed: (a)Original, and (b) the new $B$ and $C'_{new}$. }
\end{center}
\end{figure}

{\bf Example 4.2}\\
Example of moving sequence could be $C'$, $C''$, $\cdots$, $C^{(t)}$ in $B$.
The entering point and the leaving point will stay the same. each time change by one triangle in $B$. The dashed curve in the Fig. 19 shows the possible deformation.
Since any such $C^{(i-1)}$ and $C^{(i)}$  differs only the boundary of a triangle, two such triangles in $Star(b_1)$ must share a 1-cell since gradual variation or deformation is required. So the sequence may pass a "Sectors"  type. The center point could be vanished Fig. 19 (c).

The $\alpha$ (starting at $b_0$ to $D'$) in $D'-B$ will not be directly entering $b_1$, $b_2$, $b_3$, etc. in Fig. 19 (a) since $C''$ also contains $\alpha$. This is because, otherwise, it will generate a loop inside of $C''$  not only $C''$ itself.\\

\begin{figure}[h]
\begin{center}
   \epsfxsize=5in
  \epsfbox{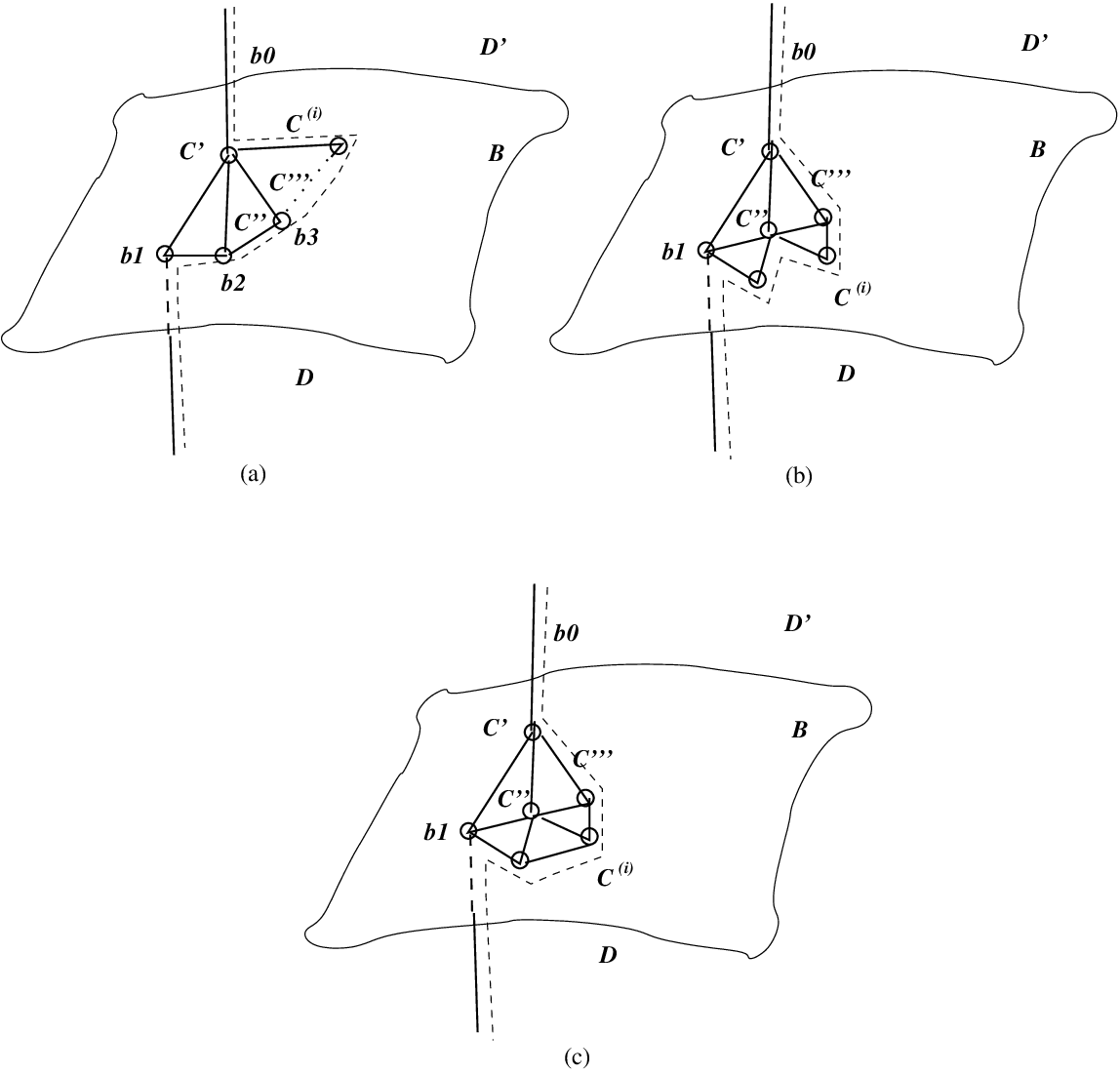}   
\caption{Check ahead for the pattern  where $Modulo2(C',C'')$ is a triangle totally in $B$. }
\end{center}
\end{figure}

{\bf Example 4.3}\\
The third example, in Case (1.2.b) and (1.2.c), when we make $C'_{new}$ for an arc in $C'$ in $D'-B$, we may use a vertex in $C''$ in $B$. We said that we need to check ahead for $C''$  in the case where such a vertex should be avoid to use in $C'_{new}$. So we can avoid to use any of the point in the triangle in $C'_{new}$. We also said to repeatedly use this strategy$\/$method finite times if it is needed. It means that there is a $C'''$ that is in the same situation in the sequence of $C'$, $C''$, $C'''$, and etc. In fact, we need to find the end of the consecutive cases, to do a unified projection.

Here is the detailed explanation for the cases in (1.2.b). We can do a similar process for the cases in (1.2.c) .

For Fig. 20,  when some points were taken by $C'_new$ in the next few moves in $\Omega$ (Fig. 20 (a)), we
can find out the last one and use the project for the last to replace the corresponding piece in $C'_new$. Remember that the difference between $C'$, $C''$, etc are
jut the triangle boundary shown in Fig. 20. When projection for $C^{t}$ in Fig. 20 (b) needs different subdivision for $B$. We will use finest one
for all concerned for these moves in $\Omega$ (meaning that each point in $C''-C'$ , $\cdots$,$C^{(t)}-C'$ were taken by $C'_new$ originally in a projection.  In order to maintain
as smaller changes as possible for $B$ in terms of subdivision, we do not have to insert too many new curves between two original $C^{(i-1)}$ to $C^{(i)}$. But we need to prove that
$C^{(-1)}_{new}$ to $C^{(t)}_{new}$ can be gradually varied. This is because that any two curves in $B$ will be gradually varied. We can always find a sub-sequence to be inserted. Using modulo2 to exchange edges
when the case of two cross-over arcs happened since $B$ is the simply connected and orientable closed 2-manifold.  In fact, another way, all such arc's projection will be changed. It is true that all cases.    ).

\begin{figure}[h]
\begin{center}
   \epsfxsize=3.5in
  \epsfbox{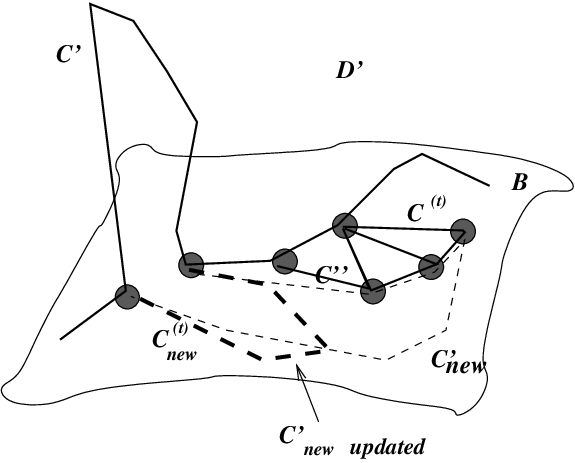}   
\caption{Points in $C''-C'$ , $\cdots$,$C^{(t)}-C'$ were taken by $C'_{new}$ originally in a projection. Find one for $C^{(t)}_{new}$ to replace it in others including $C'_{new}$. The arc was the bold-dashed in $C'_{new}$-updated.  }
\end{center}
\end{figure}

\subsection{Some Additional Discussion for the Proof of the Above Main Theorem: }

There are some other procedures and algorithms for proving the main theorem, Theorem 4.3. These algorithms may be faster in practice of contraction and deformation.   We will discuss them in other papers.

\subsection{The Main Theorem for $m$-manifolds}

Besides Theorems 4.1 and 4.2, in Theorem 4.3, we actually proved that for 3D compact smooth manifolds, the Jordan separation theorem implies the Poincare conjecture. We can do this because we have used the property of finiteness of cells in complex $K$, and
the contraction is based on a procedure described in the above proof. This type of algorithm can only be performed in finite
cases and in discrete-like manifolds. Therefore,

{\bf Corollary 4.1} For 3D closed smooth manifolds that admit the finite simplicial complex decomposition, the Jordan separation theorem (Theorem 4.2) implies the Poincare conjecture in 3D.

Using this result, we can design algorithms or procedures to show that for cell (or simplicial) complexes, a closed simply connected smooth $m$-manifold that has finite cell (or simplicial) complex decomposition (meaning that it admits a finite triangulation) is homeomorphic to $m$-sphere. We can describe this observation as follows.

The Chen-Krantz paper also presented the proof of the following theorem:

\begin{theorem}
(The Jordan Theorem of Discrete $m$-manifolds) Let $M_m$
be a simply connected discrete or PL $m$-manifold. A closed simply connected discrete $(m-1)$-manifold $B$ (with
local flatness) will separate $M_m$ into two components. Here, $M_m$ can be closed.
\end{theorem}

With the recursive assumption, $B$ can be assumed to be homeomorphic to an $(m-1)$-sphere. We can change the terminology to include
cell complexes in case of embedding to Euclidean Space.

\begin{theorem}
(The Jordan Theorem of $m$-manifolds that admit the finite triangulation) Let $M_m$
be a simply connected (finite-cell) PL $m$-manifold. A closed simply connected $(m-1)$-manifold $B$ in $M_m$ separates $M_m$ into two components. Here $M_m$ can be closed.
\end{theorem}

Followed by the proof of Theorem 4.3 in this section, we could say:

\begin{theorem}
If $M_m$ is closed in Theorem 4.5, then we can algorithmically make $M_m$ to be homeomorphic to an $m$-sphere.
\end{theorem}

Some discussions: In Step (2) of the proof of Theorem 4.3, if $B$ is not a simple cycle, then $M$ is not simply connected according to the theorem above. We will use this property to determine whether $M$ is simply connected. See~ \cite{Che17I}.

In Step (3), since $M$ contains a finite number of $m$-cells, there will always be a way to reach an end.
It is also possible to use a gradually varied "curve" or $(m-1)$-cycle on $M_m$ of $B$ to replace $B$ to reach the maximum number of removals in practice.  However, we have to keep the removals balanced to a certain point, meaning that we would always try to remove the point farthest on the boundary to the fixed point $o$, which we contract to.

The entire process in the algorithm determines the homeomorphism to the $m$-sphere. Therefore, we would like to say that
with finite cell complex triangulation, a closed-orientable compact simply connected smooth $m$-manifold is homeomorphic to an $m$-sphere.

{\bf Corollary 4.2}   With finite PL complex decomposition or triangulation, a closed-orientable compact simply connected smooth $m$-manifold is homeomorphic to an $m$-sphere.

In papers ~\cite{Che15I} and paper (II) ~\cite{Che17I}, we used filling techniques to turn a closed manifold into a deformed disk.
As an equivalent statement, we observed that: A closed-orientable $m$-simply connected manifold $M_m$ is homeomorphic to an $m$-sphere if and only if there is an $(m+1)$-deformed disk that is simply connected with a boundary of $M_m$.

\section{Existence of the Approximating Sequence of Triangulations of a Smooth Manifold}


In the above sections, we dealt with the compact smooth manifold that admits finite triangulations or PL complex decomposition
(decomposition that is diffeomorphic to the original manifold). If the manifold is not compact, then we may not have finite cell complex decomposition. We discuss this problem in the following.

According to the definition of topological manifolds, we know that they must be locally compact. This means for any
smooth manifold, there is a finite triangulation or simplicial decomposition for any given point and its neighbor.
Our strategy is to find a sequence of PL complexes that will approach the original smooth manifold in any approximating degree.

To do this, we need to use {\it the Axiom of Choice}. Since $\mathcal M$ is smooth, no point exists where the total curvature of the point
is infinite. We can select a point with the smallest curvature and then use Whitney's triangulation method to obtain
the triangulation, which must be finite. Now, the issue is that the process is not constructive since we used the axiom of choice.

For any small number $\epsilon >0$, we would like to make a triangulation based on $\epsilon >0$ where each simplex is part of
an $\epsilon-$ cube. As a result, we would have a sequence of triangulations. Based on local compactness, there must be
an $\epsilon$ to provide a finite cover for each point with the smallest curvature.

If we do not use the Axiom of Choice, we can select a sequence of $\epsilon_i >0, i=0,1,2,3,\cdots $ where $\epsilon_0$ is a small number less than 1 and  $\epsilon_i =  \epsilon_{i-1}/2$. We can use the Whitney's triangulation for each $\epsilon_i$, then we eventually can find a $\epsilon_i$ that is smaller than the smallest total curvature of points on $\mathcal M$.

It is very interesting to think about the new discussion on $Wt$-shapes in Appendix A.  When a loop of  $Wt$-shapes (Fig. \ref{fig:ManyWt1} and Fig. \ref{fig:ManyWt2}), even one $Wt$-shape  loop itself, attached (connected-sum) to a base $Wt$-shape, and each of the $Wt$-shape will attach such a loop, we do it forever, so this manifold will be not a simply connected intuitively since the collar will be getting smaller and smaller until they are closed. Is it possible that we can set up a
control parameter such that the loop attaching to the ``collar''  of $Wt$-shapes in the manifold is keeping its self a simply connected.   If we can find such a parameter, then the technique of using  $\epsilon_i$  will not
guarantee we have a simply connected smooth 3-manifold that must be homeomorphic to a 3-sphere if it is not compact.  More theoretical research could be done for finding some more interesting 3D shapes.

We will discuss the algorithm further in the next paper.

\newpage


\section{Conclusion}

The conclusion of the previous version is the following: In this paper, we wanted to make some connections between discrete mathematics and combinatorial topology, along with including a discussion on smooth manifolds.

We added further details and corrected several typos for the proof of a theorem originally published in ~\cite{Che17I}. Since our research on these topics is still ongoing, there may be some areas that warrant further study and discussion in future works. We believe that the direction of our current research is valuable for understanding some essential problems in geometric topology.

This revision added more figures to the extended version of the earlier revisions of {\it Algorithms for Deforming and Contracting Simply Connected Discrete
Closed Manifolds (III)}. In the second reversion, we added more details in the proof of the main theorem. In this revision, we also added more explanations and more figures.

Some discussions without using the general Jordan Theorem was presented in ~\cite{Che18SummerTop} where we make ``the general Jordan Theorem'' as a condition due to reason Chen-Krantz paper has not published in a journal yet.
In this new revision, some detailed analysis was provided especially the work in Appendix A (This new revision added Appendix A to add the details of two cases).
According to the Appendix A, some part of \cite{Chen-Krantz} needs to be modified. Appendix A in this paper was solely done by Li Chen.

\newpage

\section{Appendix A: The Proof of Existence of an $m$-cell of $D$ Having Two $(m-1)$-faces on Boundary $B$ in the Part 2 of the Main Theorem 4.3}

Assume that the contractional-removal have removed many 3-cells ($m$-cells) until the following two conditions hold for all 2-cells on the 2-cycle $B$ that is the boundary of $D$ in the proof of Part 2 of Theorem 4.3:

For each 3-cell $e$ having a 2-cell (2-face) $f$ on $B$, we have:

{\it Case 1}:  The 0-cell $P$ (in $e$) not in $f$ is on $B$. Or

{\it Case 2}:  The 0-cell $P$ and the 1-cell containing $P$ not in $f$ is on $B$.

We want to make sure the following statement is true: There must be a 3-cell has two 2-cells on $B$. Note that Case 2 is a special case of Case 1.

{\bf Theorem A:} There is an $m$-cell of $D$ having two $(m-1)$-faces on Boundary $B$ if for each $m$-cell $e$ having a $(m-1)$-cell (2-face) $f$ on $B$, the following conditions hold:
{\it Case 1}:  The 0-cell $P$ (in $e$) not in $f$ is on $B$. Or {\it Case 2}:  The 0-cell $P$ (in $e$) and the 1-cell containing $P$ not in $f$ is on $B$.

Before we prove Theorem A for 3D cases or higher dimensional cases, we first prove it for 2D cases. Because in 2D, Case 2 does not exist (Otherwise we admits two 1-cells on $B$ already).

{\bf Theorem B:} There is an $2$-cell of $D$ having two $1$-cell on Boundary $B$ of $D$ if for each 2-cell $e$ having a 1-cell (2-face) $f$ on $B$, the following condition holds:
: The 0-cell $P\in e$ not in $f$ is on $B$ .

{\it Proof:}

We define the distance between two cells $E$ and $E'$ is the shortest distance (the length of a shortest edge-path or 1-cell path) between two points in each $E$ and $E'$. $E$ and $E'$ can be in different dimensions.  Denote $d_{\Omega}(E,E')=\{ n | n \mbox {is the length of a shortest path from a 0-cell in $E$ to another 0-cell in $E'$\ in $\Omega$}$.

We now recall the concept of $i$-cell-distance: the length of the shortest $i$-cell path between two cells (usually 0-cells, include these two cells, these two cells are not the same) where each adjacent pair of two $i$-cells share an $(i-1)$-cell. In general we use $d_{\Omega}^{(i)}(E,E')$ for $i$-cell length that always shares an $(i-1)$-cell between $i$-cells in the path.  Use $d_{D}^{(i)}(X_0,X_1)$ to denote this distance for two points or two cells.  For instance, 1-cell distance is the graph distance. We know that each pair of 3-cells in $D$ is 2-connected.
We can identify two points on $B$ that has the largest 3-cell-distance in $D$.
We will prove that we can find $e$ nearby one of the two points that has two 2-faces on $B$.

We start with the example in two-dimensional cases. See Fig.~\ref{fig:TwoDCase}

This manifold has the following property: (a) there is no inner point, and
(b) each 1-cell on boundary has associated 2-cell that contains this 1-cell. This 2-cell has the third point on the boundary too.

\begin{figure}[h]
\begin{center}
   \epsfxsize=4.5in
  \epsfbox{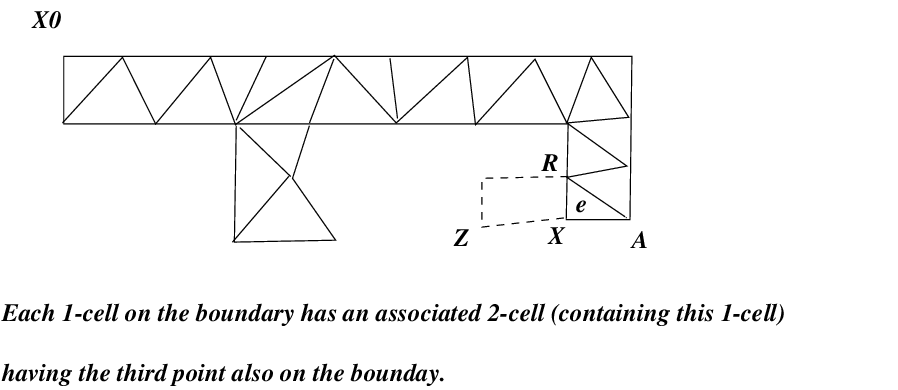}   
\caption{ A 2D example where exists a 2-cell that must have two 1-cells on the boundary $B$.} \label{fig:TwoDCase}
\end{center}
\end{figure}

{\bf We now to prove for the 2D case: the first proof}
The assumption is: For each 2-cell $e \in D$, if $e$ has a 1-cell $f$ in $B$ ($B$ is a 1-cycle) and the third point $P$ of $e$ not in $f$ is also on $B$, then there must be a 2-cell $e'$ that has  has two 1-cells on $B$.
So, we can do the contractional removal of the corresponding 2-cell $e'$ to maintain that the new boundary is still a 1-cycle after removal of $e'$ .

We will use 2-cell distance. We will see that the far-most pair of points in $D$ in 2-cell distance will help to determine such an $m$-cell, $m=2$, that has two $(m-1)$-cells on $B$. The good thing in 2D is that a 2-cell, having a 1-cell on $B$ in $D$, contains a 1-cell not on $B$ that will split $D$ into two disconnected parts. Assume that these two points are $X_0$ and $X$ (meaning that this pair has the longest 2-cell distance). Let the 2-cell containing $X$ is $e=ARX$. If $e$ does not have two 1-cells on $B$, we can assume the 1-cell $AX$ is on $B$ ($R$ is at another side of $B$ from $A$). See Fig. ~\ref{fig:TwoDCase}. Since $B$ is a 1-cycle, then there must exist another 1-cell on $B$ that contains $X$. (Any 0-cell is contained by exact two 1-cells.) This 1-cell can be named $XZ$, $Z$ is on $B$, and from $Z$ to $R$ there is a 1-cell path (since $B$ is connected). $Z$ is the point that has strictly longer 2-cell-distance comparing to $X$ from $X_0$ since any 2-cell path to $Z$ from $X_0$ must include cell $e$. So, we get the contradiction that $X$ is the far-most point from $X_0$. Therefore, there must be a 2-cell that has two 1-cells on $B$.

Thus, we proved the 2D case. The idea of proving for the 3D case can be similar as the case for 2D but having more complexity.  The above simple proof for 2D can be inserted to \cite{Chen16, Chen14} as a supplement fact for the 2D case.  

Note that in this proof, we can see that we only need to get a pair that has the local maximum distance (in $Star_D(X)$ from $X_0$) in 2-cell path of $D$. In addition, $X_0$ can be any point on $B$.

{\bf The second proof:} The first proof is simple and direct. However, the technique of the second proof will be used in the 3D and higher dimensional cases.
In fact, we have the second proof of this Theorem using just edge-distance.  Unlike it was very intuitive regrading to the first proof. The second proof needs more explanations. Let's look at the following example:

\begin{figure}[h]
\begin{center}
   \epsfxsize=4.5in
  \epsfbox{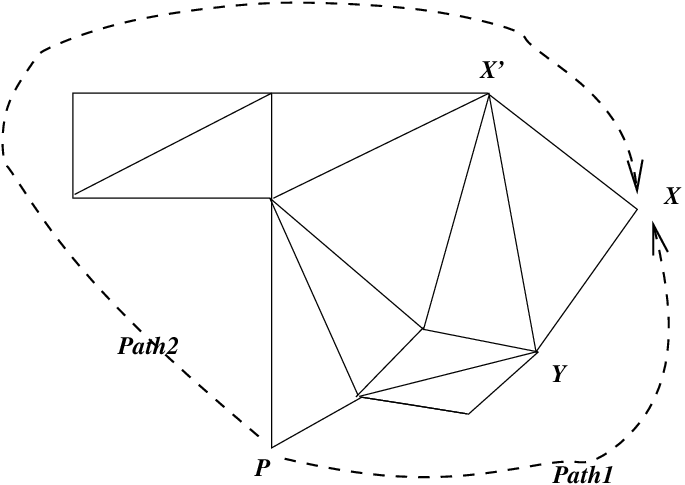}   
\caption{The second 2D example where exists a 2-cell  having two 1-cells on the boundary $B$ .} \label{fig:TwoDCase2}
\end{center}
\end{figure}

Let's check paths from $P$ to $X'$, this type of paths is the edge-path or the 1-cell path in $D$. The 1-cell distance from $P$ to $X'$ is 3 (the longest, $d_{D}(P,X')=d_{D}^{(1)}(P,X')$ reaches maximum. ), and  1-cell distance from $P$ to $X$ is also 3.  That is why the first proof is more intuitive. So we need to use another property of $B$, there is a simple 1-cell path from $X'$ to $P$ on $B$. This path could be different from the path we used to calculate the 1-cell distance in $D$. The path on $B$, $\pi(X',P)$, could be much longer. However, we can always find a pair of points $X$ and $Y$ on $\pi(X',P) =X' \cdots X Y\cdots P$ such that $d_{D}(P,X)$ is the maximum (local maximal is also fine) and  $d_{D}(P,Y) = d_{D}(P,X) -1$. We can assume that $X$, $Y$ is the last pair satisfying this condition. Since $X$ on $B$  ($B$ is an 1-cycle) has two adjacent points on $B$. We can assume they are $X''$ and $Y$.  There are only two possible cases: a)
$X''XY$ are in a 3-cell (see Fig. ~\ref{fig:TwoDCase2} ), so there are two 1-cells of a 2-cell of $D$ on $B$. We proved the statement. b) $X''X$ and $XY$ are in different 2-cells. Assume that $X''X$ is included in a 2-cell in D and the third point is denoted as $A$. $A$ is also on $B$ as assumed. Therefore, any 1-cell path in $D$ from $P$ to $X''$ must pass $X$ or $A$. It is also to say that 1-cell $XA$ split $D$ into two (or more) parts.
Since $d_{D}(P,X)$ is the maximum length. If $d_{D}(P,A)$ is also  the maximum length, then  $d_{D}(P,X'')>d_{D}(P,X)$ . This causes a contradiction.  So we must have $d_{D}(P,A)=d_{D}(P,X)-1$. (It is impossible that $d_{D}(P,A)\le d_{D}(P,X)-2$ since $d_{D}(X,A)=1$. )  Thus, if $X''A$ is on $B$ than we have the 2-cell $AX''X$ containing two 1-cells on $B$. Otherwise, $X''A$ (not on boundary) must be contained by another 2-cell $X^{(3)}X''A$. $X^{(3)}X$ must be on $B$. If not, there will be another 2-cell $e''$ that contains $X^{(3)}X$. $e''$ must contain a point that has the (absolutely) longer 1-cell distance to $P$.  So we can exam if $X^{(3)}A$ is on $B$ or not. If it is true, we
proved our theorem since there are two 1-cells on $B$ in 2-cell $X^{(3)}X''A$. Otherwise, we continue to have $X^{(4)}$, $\cdots$, $X^{(k)}$.   Since $Star_D(A)$ is finite ($D$ is finite), we can not $X^{(k)}$ forever. There must be a finite number $k$ such that 2-cell $X^{(k)}X^{(k-1)}A$ contains two 1-cells on $B$. So we completed the second proof.

(Note that it is not passible to have other paths such as $Path1$ or $Path2$ since $X'$, $X$, and $Y$ are on the boundary $B$ in Fig. \ref{fig:TwoDCase2}. )      $\hfill$ $\square$

As we discussed in the proof of Theorem A,  It can be added to \cite{Chen14_1, Chen14} as a supplement fact for the 2D case. Since in  \cite{Chen14_1, Chen14}, we discussed the general 2-cell not only triangles. In fact
we only need to prove:  There must be a 2-cell (polygon) $e_2$ where $B\cap e_2$ is homeomorphic to a 1-cell ( $B\cap e_2$ is connected on $B$ ) and  ${\partial e_2}-B\cap e_2$ is also homeomorphic to a 1-cell ( ${\partial e_2}-B\cap e_2$ is connected on $B$ )  under the condition for
each 1-cell on $B$, if it is in 2-cell $e'$ then all 0-cells in $e'$ is also on $B$.  The proof will be the same as above.     

{\bf {\it Proof of Theorem A:}}

In this proof, we mainly prove this theorem for $m=3$. The extension is very natural with very minimum changes.

In the proof of Part 2 in Section 4.2,  there are some extreme cases
that need to consider the order of contractions for 3-cells ($m$-cells) in $D$. Specifically, we could not contract a 3-cell $e$ when

{\it Case 1}: A 2-face ($f$) of $e$ is on $B$ and any of the three other 2-faces is not on $B$, but a 0-cell $P$ not in $f$ is on $B$. Or

{\it Case 2}: A 2-face ($f$) of $e$ is on $B$ and any of the three other 2-faces is not on $B$, but a 0-cell $P$ and a 1-cell (containing $P$) not in $f$ is on $B$.

\noindent It is easy to see that contracting an $e$ in Case 1 or 2 could make $B$ that is not a 2-manifold. In other words, we just could not contract this 3-cell $e$ at this ``moment.’’ We need to remove its neighboring cells before contracting this one ($e$).  This is because in such two cases, all 0-cells of $e$ are in $B$, $e$ only contain four 0-cells if $m=3$, the boundary of $e\cap B$ is not 1-cycle.

Note that these two cases usually not happen, but if it happens we will do the following special process to find $e$'s neighbor or its neighbor's neighbor, $\cdots$, to remove them beforehand.

In any of the two cases, we can use the following process: We can first search for a 3-cell that has two or more 2-faces on $B$; denote it as $e'$,  we contract $e'$ first.
We will prove that such an $e'$ always exists on $B$.

The main idea of the proof is presented as follows: On the contrary,  if there is no such a
3-cell $e \in D$ that have two (or more) 2-faces on $B$ for 3D, then $B$ only intersects $e_i$'s each of which has just one 2-face $f_i$ on $B$. Since each $f_i$ contains three 0-cells, the forth 0-cell $P_i$ of $e_i$ will be also on $B$ in both cases (Case 1 and Case 2). Again, recall the process of contraction in Section 4.2, if $P_i$ is not on $B$, meaning that $P_i$ is an inner point in $D-B$, then we can contract $e_i$ still keeping $B$ a 2-cycle. In other words, we can remove such 3-cell $e_i$ by using three other 2-faces of it to replace $f_i$.  (Note: Let $O$ is the point we want to contract $D$ to $O$. We should keep $O$ as the special point. That means if we want to keep $Star(O)$, $P_i$ cannot be $O$. But we do not have to keep $Star(O)$. It is just time saving concerns.)

For $m$-dimensional cases, the above discussion could mean for an $m$-cell $e$ with two $(m-1)$-faces on $B$ for the $m$-manifold in $n$-dimensional Euclidean space. $e$ has $(m+1)$ 0-cells. Any two $(m-1)$-faces  of $e$ will
contain all $(m+1)$ 0-cells in $e$.

How are we certain that there must be a 3-cell of $D$ having two or more 2-faces on $B$? Before the strict proof, our observation can be described here: (1) $B$ is a 2-cycle, there is a simple 2-cell path (1-connected, every pair of adjacent\/neighboring 2-cells sharing a 1-cell) on $B$ that link $f_i$ to $P_i$. (``link'' is not a mathematical term here, but just think a 2-cell containing $P_i$.)  Since each 2-cell on boundary $B$ will be included in a 3-cell in $D$, but the 0-cell of it  not in the 2-cell on $B$ will also on $B$.  $D$ is very thin in most part.    Think about that $D$ is a very thin pad  that could be bend in a flexible shape. This 2-cell path most likely to  wrap a ''corner\/edge'' 3-cell in order to reach $P_i$. Note that $f_i$ and $P_i$ are in a 3-cell in $D$.  This ''corner\/edge'' 3-cell is most likely containing two 2-faces on $B$.
(2) How do we define this ''corner\/edge'' 3-cell? We can think about it is a far-most 3-cell related to this path (having a none-empty intersection to this 2-cell path). (3) If this ''corner\/edge'' 3-cell does not contain two 2-faces on $B$, (only has one 2-face $f_j$ on $B$,) then we start with $f_j$ to find $P_j$. Repeatedly, we find a simple 2-cell path from $f_j$ to $P_j$ on $B$. This procedure cannot be forever running since $D$ only contains finite number of cells. One of the tasks is to determine the structure near the ``corner\/edge'' 3-cell.  When  this structure is clear, then we can prove this theorem.

Now we start to prove this theorem.

{\bf Case 1: }

For Case 1, we want to prove there must be a 3-cell $e'$ intersect with $B$ that has two 2-faces on $B$ under the condition: Each 2-cells on $B$ is in  Case 1 or Case 2.  For details, the Case 1 means that a 2-face $f$ of $e$ is on $B$ and a 0-cell of $e$, $P$, which is not on $f$ is also on $B$.
In other words, each 2-cell $f$ on $B$ that is contained by a 3-cell $e$ in $D$. The point in $e$ not in $f$, denoted $P$, is also on $B$.  Case 2 is a special case of Case 1.
Since $P$ is on $B$, $Star_B(P)$ contains a triangle $f_P$ is on $B$ with $P\in f_P$. Intuitively, $f_P$ is located in the other side of $f$ in $D$.  Since $B$ is a 2-cycle, from $f$ to $f_P$, there is a 2-cell path on $B$, and the path is 1-connected. We can denote the path as $\pi(f,f_P) = (f=f_0), f_1, \cdots, (f_N=f_P)$. $N$ is finite. $f$ and $f_N$ are in $Star_M(P)$.

Every $f_i$ is included in an $e_i$. There is a $P_i$ (the forth point in $e_i$) such that $P_i \in e_i$ but $P_i\notin f_i$. Intuitively, $D$ would be a very thin pad, but it could be bend.

(We can also assume that $N$ is the length of a shortest 2-cell path with 1-connectedness from $f$ (from any point of $f$) to $f_P$ on $B$. But it is not really important, we just need a simple 2-cell path.)

We consider $f$ moves along with the path when $f_i$ moves, the 3-cell $e_i$ also moves as $i=1,...,N$. Some $e_i$ could be the same. If that is the case, we have proved our statement: there is always a 3-cell adjacent\/intersect to $B$ contains two 2-faces on $B$.

In order to make the proof easier, we like to find a 0-cell $X$ on $B$ that is the far-most to $P$ in terms of the 3-cell distance. (Originally, we want that $X$ is in the path $f=f_0, f_1, \cdots, f_N=f_P$.) We can find $X$ first, then find a 2-cell $f_X$ on $B$ that contains $X$.  We can make a shortest 2-cell path, $\pi_2(f_0,f_X)$ , on $B$ from $f_0$ to $f_X$ then we can make a shortest 2-cell path $\pi_2(f_X,f_P)$ from $f_X$ to $f_P$ without using the 1-cell in the previous $\pi_2(f_0,f_X)$. It is possible this path is not a simple 2-cell path, however, we just need to make sure not have a``cross''  in $Star_B(X)$. Actually, we just need to look at the $\pi_2(f_X,f_P)$. There is no need to check the entire 2-cell path from $f_0$ to $f_P$.

We define the distance between two cells $E$ and $E'$ is the shortest distance (the length of a shortest edge-path or 1-cell path) between two points in each $E$ and $E'$. $E$ and $E'$ can be in different dimensions.  Denote $d_{\Omega}(E,E')=\{ n | n \mbox {is the length of a shortest path from a 0-cell in $E$ to another 0-cell in $E'$\ in $\Omega$}$.

We now recall the concept of $i$-cell-distance: the length of the shortest $i$-cell path between two cells (usually 0-cells, include these two cells, these two cells are not the same) where each adjacent pair of two $i$-cells share an $(i-1)$-cell. In general we use $d_{\Omega}^{(i)}(E,E')$ for $i$-cell length that always shares an $(i-1)$-cell between $i$-cells in the path.  Use $d_{D}^{(i)}(X_0,X_1)$ to denote this distance for two points or two cells.  For instance, 1-cell distance is the graph distance. We know that each pair of 3-cells in $D$ is 2-connected.
We can identify two points on $B$ that has the largest 3-cell-distance in $D$.
We will prove that we can find $e$ nearby one of the two points that has two 2-faces on $B$.

{\bf The first proof for the existence for  $Wt$-Shape using 3-cell distance:}

From the proof in 2D cases (Theorem B), in the first proof, we selected two points $X_0$ and $X$ on $B$ that have the far-most $m$-cell distance ($m=2$).
In the proof, we can see that it is not important for where $X_0$ is, we just want to determine the far-most $X$ from any $X_0$. This fact is the same as we prove for 3D.

Let $d_{D}^{(3)}(X_0,X_1)$ denote $3$-cell-distance from $X_0$ to $X_1$. Again if $X_0=X_1$, $d_{D}^{(3)}(X_0,X_1)=0$. If $X_0$ is not $X_1$ but they are in a 3-cell, $d_{D}^{(3)}(X_0,X_1)=1$. Let $a_1,\cdots, a_n$  be a 2-connected 3-cell path. Assume $X_0\in a_1$ and $X_1\in a_n$. If $n$ is the smallest, then  $d_{D}^{(3)}(X_0,X_1)=n$.  (Note that we use $n$ as the 3-cell distance, should not have a confusion to $n$-dimensional Euclidean space
two $n$'s are different.)

Once we find $X_0$ and $X$, we can have a 2-cell path $\pi$ on $B$ just denote it as  $f_i, i=0,...,N$. We usually think about $X_0$ is $P$. But for the first round of findings. We just assume that  $X_0$ and $X$ are two 0-cells. $\pi(X,X_{0})=\pi(X,P)= \{f_X=f_0, \cdots, f_N=f_P\}$. Note that $N$ is for 2-cell path while $n=d_{D}^{(3)}(X,P)$. Each $f_i$ will belong to a 3-cell $e_i$ in $D$.

If there are two or more such $X$'s (having the same 3-cell distances) we pick up the biggest $i$. (i.e. if there are multiple such point $X$, we pick the latest one in $e_0,\cdots, e_N$. ) So we make sure that $f_{k+1}$ contains a point $Y$ such that

$d_{D}^{(3)}(P,Y)\le d_{D}^{(3)}(P,X)-1\le n-1$.

(Note that, if we use  $d_{D}^{(1)}$, the edge-distance, we still can prove the same result. Here the equation will be $d_{D}^{(1)}(P,Y)=d_{D}^{(1)}(P,X)-1=n-1$. In fact,
we used $d_{D}^{(1)}$ first, but later we found $d_{D}^{(m)}$ will be more intuitive for 2D. As long as we can find such $Y$, we will be able to continue our proof. For 3D or higher dimensions, we do not need to assume we found the exact same point for $X$ in two different measures. For 2D, we can easily prove that $X$ will be the same.

We will provide the extra explanations to another case that seems more general but we will also see they are the same. )

Now we can assume that $f_0$ is $f_X$. (This is because we only care about $f_P$ to $f_X$) Let $X_0=X$. If there are multiple $X$ in $\pi_2(f_X,f_P) =\{f_0,\cdots, f_N\}$ such that $d_{D}^{(3)}(P,X)= d_{D}^{(3)}(P,X_0)=n$.  We always find the last one in the path and name it as $X$. Find $k$, $X$ is in $f_k$ and $f_{k+1}$. But $X$ is not in $f_{k+2}$.   So we make sure that $f_{k+1}$ contains a point $Y$ such that  $d_{D}^{(3)}(P,Y)\le d_{D}^{(3)}(P,X)-1 \le n-1$.  If $X$, $Y$ are already in $f_0$ and satisfying the condition. We need to extend another $f$ such that $f$ is in $Star_B(X)$ and $f\cap f_0$ is a 1-cell that contains $X$. We can call $f$ as $f_{-1}$.

Note that $f_k$ is always further than $f_{k+1}$ in 2-cell path $\pi$ on $B$. We want to find $X$ that is far most point from $P$ in 3-cell distance in $D$ for points on the 2-cell path $\pi$. And $Y$ is absolutely closer to $P$ in 3-cell distance. ($k$ could be -1. Or we can just shift one index for each $f_i$.  }

\begin{figure}[h]
\begin{center}
   \epsfxsize=5.0in
  \epsfbox{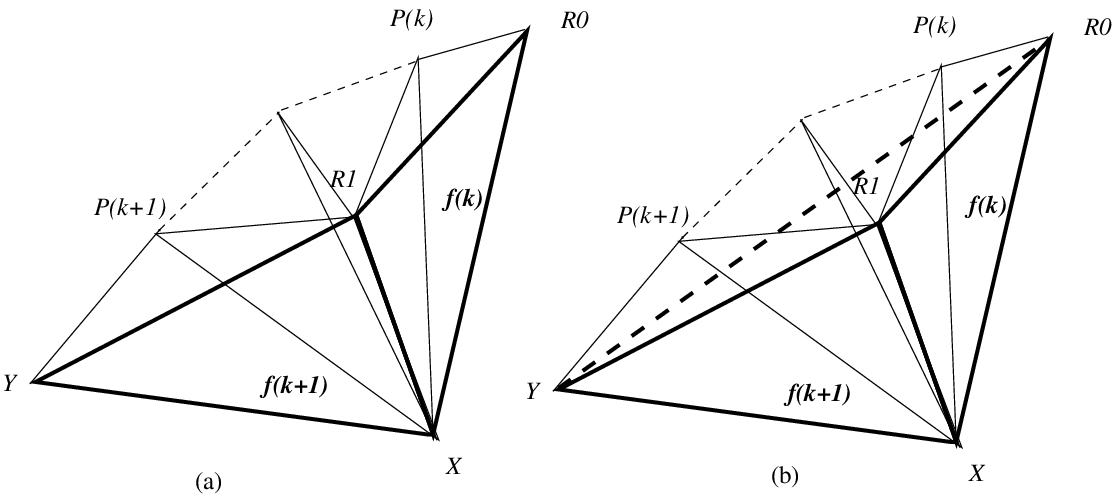}   
\caption{ A general case of a 2-cell path on $B$: basic .} \label{fig:XR0R1YTwoCells}
\end{center}
\end{figure}

In Fig. \ref{fig:XR0R1YTwoCells}, we have two consecutive 2-cell $f_k$ and $f_{k+1}$ on the 2-cell path $\pi$.

Let $f_k = R_0R_1X$ and  $f_{k+1} = R_1XY$. See Fig. \ref{fig:XR0R1YTwoCells} (a)). So $f_k$ and  $f_{k+1}$ share 1-cell $R_1X$. If $YR_0$ is a 1-cell in $D$, then we have a 3-cell $R_1R_0XY$ having two 2-cells $f_k$ and $f_{k+1}$ on $B$. See Fig. \ref{fig:XR0R1YTwoCells} (b)).  We did the proof.

Let’s assume If $YR_0$ is not a 1-cell in $D$. We will have $\{P_k,f_k\}$ and $\{P_{k+1},f_{k+1}\}$ are two 3-cells. It is possible that  $P_k= P_{k+1}$, see Fig. \ref{fig:XR0R1YTwoCells} (a).

There could be a sequence of $P^{(0)}_k=P_k,\cdots, P^{(t)}_{k}=P_{k+1}$ .   Each of these points is a vertex of a 3-cell between $f_k$ and $f_{k+1}$. Each pair $P^{(i)}_k$ and $P^{(i+1)}_k$ is 1-cell. And $P^{(i)}_k P^{(i+1)}_kXR_1$ is a 3-cell. $P^{(0)}_k=P_k$ and $P^{(t)}_{k}=P_{k+1}$ are on $B$.

The idea of the on-going proof is the following:

Since $X$ is on $B$, $Star_B(X)$ is a disk and $Link_B(X)$ is 1-cycle: $Y \to R_1\to R_0\to \cdots \to Y$.  If 2-cell $P_{k+1}YX$ is not on $B$ , we done since we have two 2-cells on $B$ in 3-cell $YXP_{k+1}R_1$.
If $P_{k+1}YX$ is not on $B$, then there is $Z$ (in $Link_B(X)$) such that $XYZ$ is a 2-cell on $B$. So $Link_B(X)= Y \to R_1\to R_0\to \cdots \to Z\to Y $ . Since 2-cell $XYZ$ is on $B$, there will be a $P_Z$ on $B$ such that
$YXZP_Z$ is a 3-cell in $D$. Note, $P_Z$ could be $P_{k+1}$.  See Fig. \ref{fig:YXZP_Z} (a).

\begin{figure}[h]
\begin{center}
   \epsfxsize=5.0in
  \epsfbox{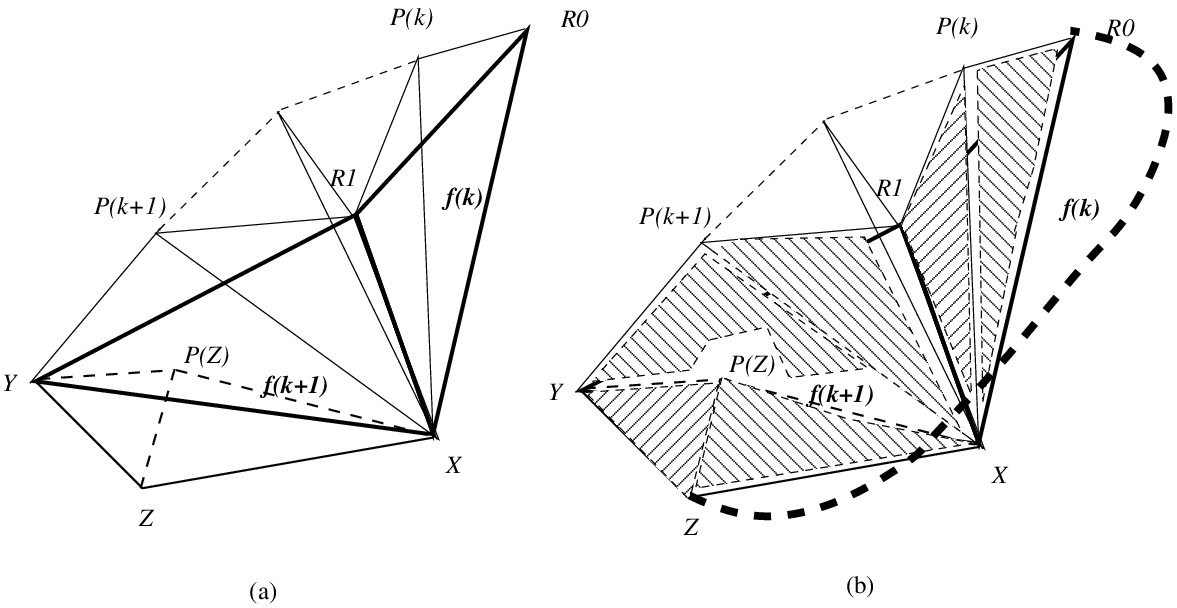}   
\caption{ A general case of a 2-cell path on $B$: more details.} \label{fig:YXZP_Z}
\end{center}
\end{figure}

Could $P_Z$ be $R_0$?  No this is because $YXZP_Z$ is a 3-cell. If it is true, $YP_Z=YR_0$ is a 1-cell in $D$. We already assumed that $YR_0$ is not a 1-cell in $D$. We also can see that $Z$ is not $R_0$, otherwise, $YR_0$ is a 1-cell.

Therefore, we have one more point $Z$ on link $Link_B(X)$. We can extend it to $Z_1$ toward to $R_0$. See Fig. \ref{fig:Z1P_Z1} (a).  So on so forth,  we can have $Z_2 , \cdots, Z_t$, and $Z_t=R_0$. If $Z_1$ is $R_0$, it is possible. See Fig. \ref{fig:Z1P_Z1} (b). So is $t$ any finite number.

\begin{figure}[h]
\begin{center}
   \epsfxsize=5.0in
  \epsfbox{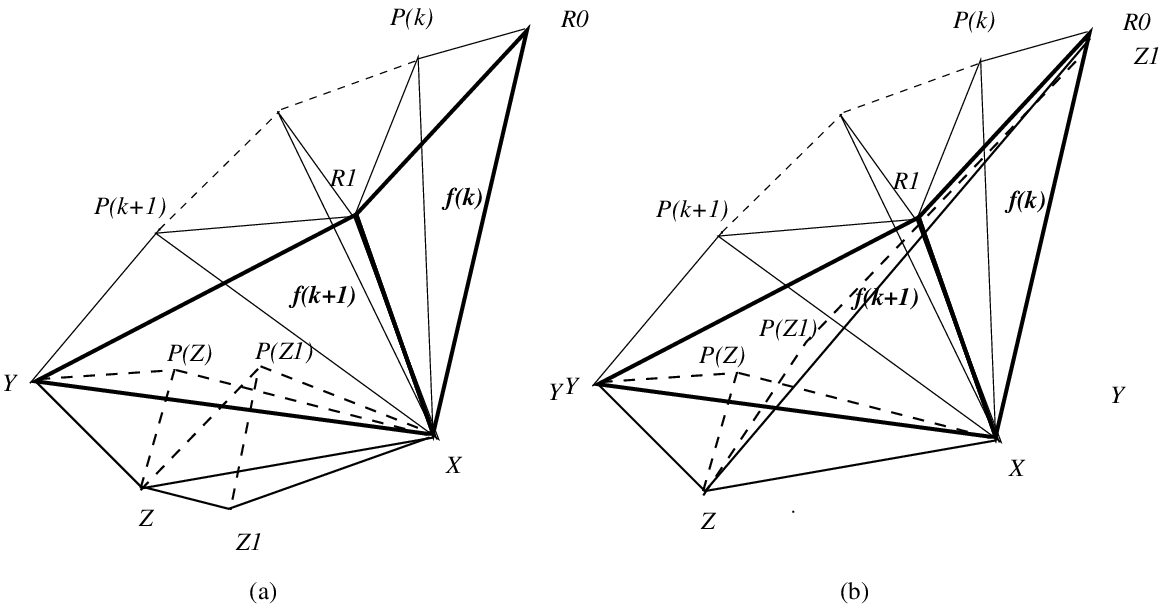}   
\caption{ A general case of a 2-cell path on $B$: could surround $X$ to reach $R_0$ without associating a 3-cell that has two 2-cells on $B$.} \label{fig:YXZP_Z}
\end{center}
\end{figure}

\begin{figure}[h]
\begin{center}
   \epsfxsize=3.0in
  \epsfbox{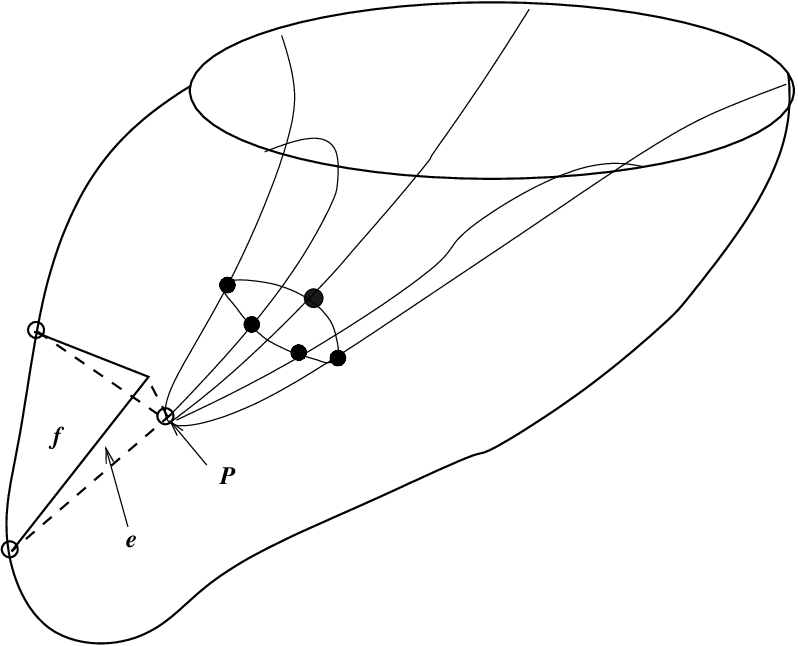}   
\caption{The Wowotou structure: A structure centered at $X$. There could be several points like $P$.  } \label{fig:Wt}
\end{center}
\end{figure}

This structure centered at $X$ is a Wowotou structure if there is no 3-cell containing 2-faces on $B$. All point $Z$, $Z_1,\cdots,Z_t$ and $P_k$, $P_{k+1}$ are also on $B$. See Fig. \ref{fig:Wt}.
(A Wowotou bun is the steam bread that has the same thickness of a corn shape. )

\begin{figure}[h]
\begin{center}
   \epsfxsize=3.0in
  \epsfbox{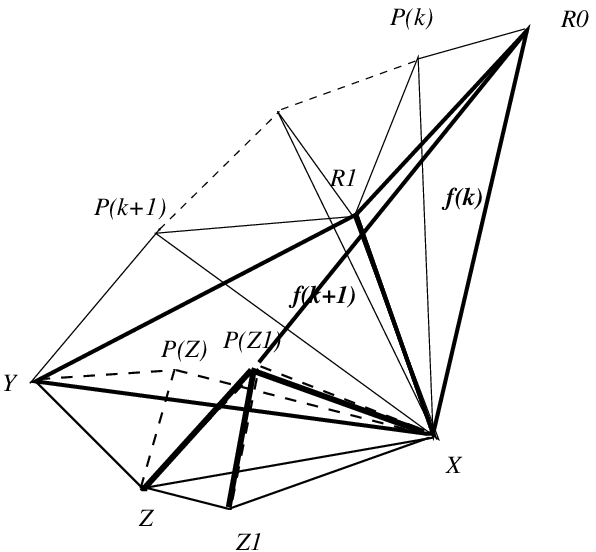}   
\caption{ 1-cell $P(Z_t)R0$ (or $P(Z)R0$) is on $B$ (in $Link_B(X)$). Then there is a 3-cell having two 2-faces on $B$.} \label{fig:6}
\end{center}
\end{figure}

It is also possible that $P(Z)$ or any $P(Z_t)$ directly links to $R_0$. and 1-cell $P(Z)R0$ or 1-cell $P(Z_t)R0$ is on $B$ (in $Link_B(X)$). This case will lead a 3-cell having two 2-faces on $B$. We did our proof.
See Fig. \ref{fig:6}.

We can see that there are only two possible choices: $Z$ or a $Z_t$ is linked to $R_0$ with a 1-cell; or, $P(Z)$ or a $P(Z_t)$ is linked to $R_0$ with a 1-cell.  As long as we have a simple 2-cell path $\pi(f,P_f)$ in $B$ from 2-cell $f$ to 0-cell $P_f$   where $f$ and $P_f$ form a 3-cell in $D$, we can find: (1) a 3-cell of $D$ that has two 2-cells on $B$, or (2) a Wowotou structure.

Thus, our conclusion is that the structure surrounding $X$ is either a ``Wowotou'' or a 3-cell containing two 2-cells on $B$. We will later prove: If there is no such a 3-cell that contains two 2-faces on $B$, $D$ must
contain infinitive number of ``Wowotous.''

We have used the 3-cell distance to find a ``Wowotou'' or ``Wotou'' structure. We call it a $Wt$-shape in the following text.  In fact, we can use 1-cell distance (typical  graph distance or edge-distance) to find this structure too. Just like
we did for the proof of Theorem B. Originally, we did try this method then we changed to the 3-cell distance that provide us the clear way to find $Wt$-Shape. Now we come back to use 1-cell distance(edge-distance, $d_D$ or $d_D^{(1)}$) to determine a $Wt$-Shape.

{\bf The second proof for the existence for  $Wt$-Shape using 1-cell distance:}

Still find two 0-cells $X_0$ and $X$ on $B$ that have the longest 1-cell distance in $D$. We can find a 2-cell path $\pi$ on $B$ connecting both $X$ and $X_0$. Assume $\pi$ on $B$ just denote it as  $f_i, i=0,...,N$. $\pi$ is 1-connected. $\pi(X,X_{0})=\pi(X,P)= \{f_X=f_0, \cdots, f_N=f_P\}$ where $f_X$ contains $X$ and $f_N=f_P$ contains $P=X_0$.  Note that $N$ is for 2-cell path while $n=d_{D}^{(1)}(X,P)$. Each $f_i$ will belong to a 3-cell $e_i$ in $D$. $n$ is the longest distance from $P$ to every point on $B$. For every point $Q$ in path $\pi$, we have $d_{D}(Q,P)=d_{D}^{(1)}(Q,P)\le n$. It is obvious that for $Q\in f_X$, $d_{D}(Q,P)=1$ if $P\ne Q$. $f_k$ is moving closer to $P$ in both $\pi$ and 1-cell distance in $D$ generally or eventually. So there must be a
$X'\in f_k$ in $\pi$ such that $d_{D}(X',P)=n$ and $f_{k+1}$ contains a point $Y$ such that $d_{D}(Y,P)=n-1$. Since $f_k$ and $f_{k+1}$ share  a 1-cell $(U,V)$ on $B$. If this 1-cell contain $X'$ then we see that $X'$ and $Y$ are in $f_{k+1}$. If  $X'$ is not on this 1-cell, we can see that $d_D(U,P)=n$ or $d_D(U,P)=n-1$ since $d_D(U,X')=1$. So we can let $U=X'$ if $d_D(V,P)=n$, or $Y=U$ if $d_D(U,P)=n-1$. So we can rename $X'$ to be $X$. We got the exact the same two 2-cells on $B$ as
we use  $d_{D}^{(3)}$: $f_k$ and $f_{k+1}$, See Fig.  25.  Therefore, we use exact strategy as we use above to extend the possible 2-cells from $Y$ to $R_0$ based on $Link_B(X)= \{R_1,Y, \cdots, R_0, R_1\}$. We can prove the same thing: If no 3-cells containing two 2-cells on $B$ are found, there will be a $Wt$-shape about $X$. First, we check if $R_0Y$ is a 1-cell. If it is true, we proved the statement. If not, there will be several 3-cells in between $f_k$ and $f_{k+1}$, and each of them contains 1-cell $XR_1$.  We have a sequence of  points  $P^{(0)}_k=P_k,\cdots, P^{(t)}_{k}=P_{k+1}$ as shown in Fig. \ref{fig:XR0R1YTwoCells}.  In the same way, we assume that $Link_B(X)= \{R_1,Y, Z=Z_0, \cdots,Z_T=R_0, R_1\}$. Then, start at $Z=Z_0$. $YZ_0X$ is a 2-cell on $B$, it was contained in a 3-cell and that has a forth point on $B$. See Fig. 26 and Fig. 27.

So we will have the same result as we  did in the first proof for the existence for  $Wt$-Shape using 3-cell distance.  (Not that we only need the local maximal for identify a $Wt$-Shape.)

We repeat this:
In such a case, $P_{R_0ZX}$ is wrapped by 3-cells bounded by $Star_B(X)$.  Such types of $P_{R_0ZX}$ could be several points. So the structure is like a WoWoTou.
A far-most point $X$ to $P$ that can be set on a 2-cell path $\pi(f,f_P)$ will have such a WoWoTou ($Wt$-) structure or it indicates a 3-cell containing two 2-faces on $B$.

In other words, Case 1 could contain a $Wt$-shape (WoWoTou), a circled “v” shape (the cone shape) that contains far-most
$X$ (3-cell distance) with all neighbors that could be not containing two 2-faces on $B$ in a
3-cell containing $X$. However, if it dose not contain such a $Wt$-shape, then it must indicate a 3-cells that contains two 2-faces on $B$.

\begin{figure}[h]
\begin{center}
   \epsfxsize=3 in
  \epsfbox{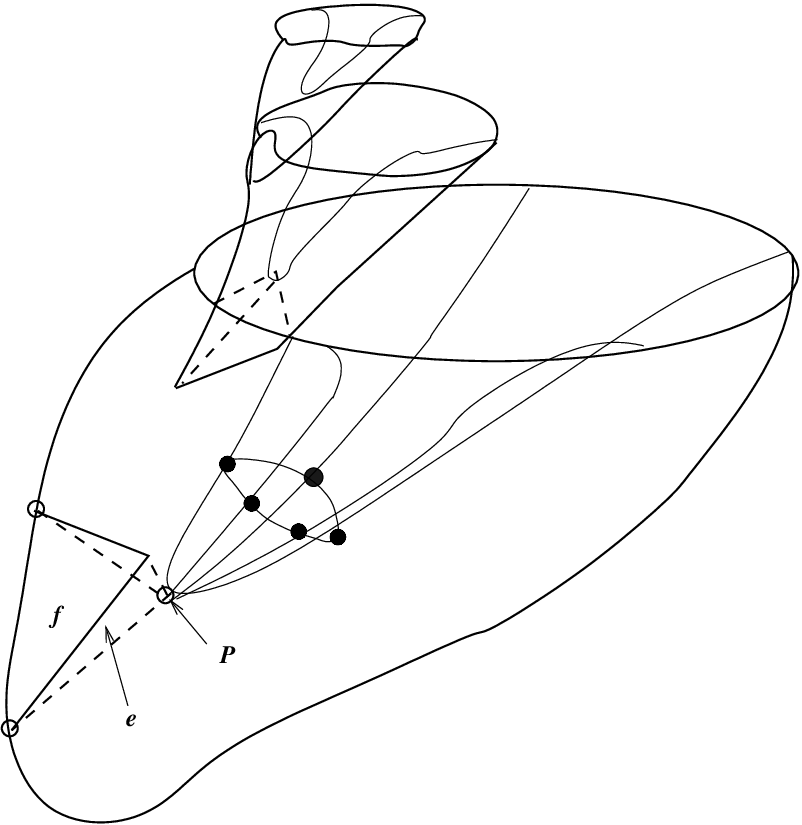}   
\caption{Many $W_t$-shapes (Wowotou structures) on a base $Wt$-shape. } \label{fig:ManyWt}
\end{center}
\end{figure}

From $f$ and $P_f$, we can determine a 2-cell path $\pi$ that contains $X$. If we cannot determine a 3-cell that contains two 2-cells on $B$. $X$ must be in a 2-cell $f_X$ on $B$, the 3-cell $e_X$ containing $f_X$ must have a forth point $P_{f_X}\notin f_X$.  Assume $e_X$ does not two 2-cells on $B$. We will have a 2-cell path $\pi_X$ from $f_X$ to $P_{f_X}$. $\pi_X$ contains a far-most point $X_1$ in 3-cell distance (or 1-cell distance) to  $P_{f_X}$, and $X_1$ is contained in a $Wt$-shape (a Wowotou structure).  If we go on, we will have a sequence $X=X_0, X_1, \cdots, X_n, \cdots$.   See Fig. \ref{fig:ManyWt}. It will be gone forever.

It is impossible since we only have finite number of cells in $D$ unless $X_i = X_j$ (or $f_{X_i} = f_{X_j}$) for some $i$ and $j$. See Fig. \ref{fig:ManyWt1}. If we follow the ``Arrow A,'' it must contain a handle. So it is not simply-connected for $D$. So the path of $Wt$-shapes only can go ``Arrow B'' in Fig. \ref{fig:ManyWt1}.  It could be still forming a handle unless one side of the path of $Wt$-shapes attaching the bottom   $Wt$-shape. See Fig. \ref{fig:ManyWt2}.

\begin{figure}[h]
\begin{center}
   \epsfxsize=3 in
  \epsfbox{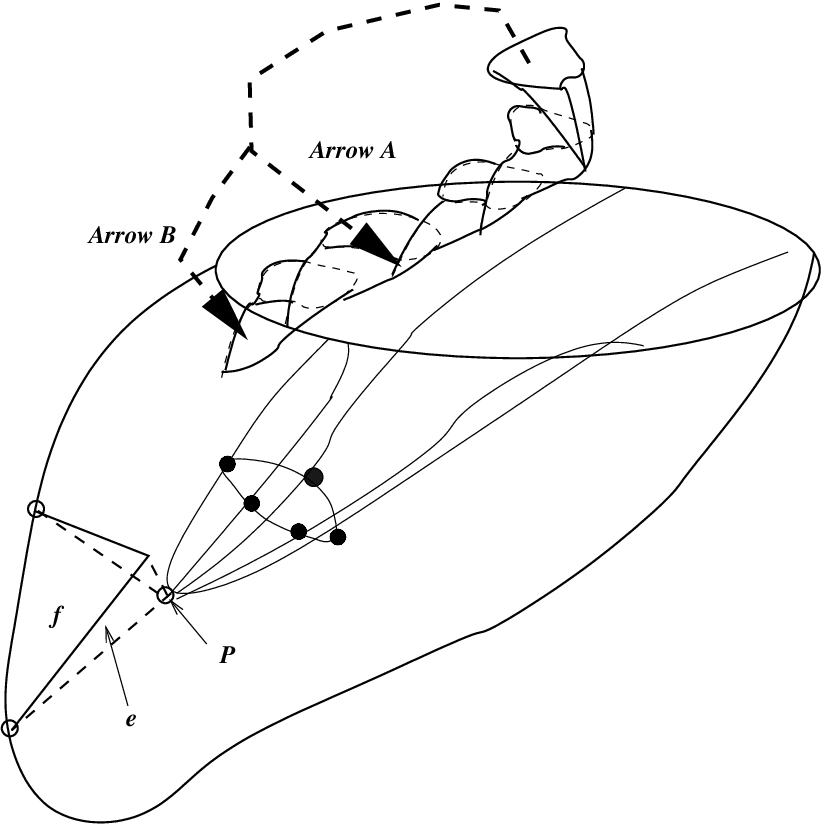}   
\caption{Looped $W_t$-shapes (Wowotou structures). } \label{fig:ManyWt1}
\end{center}
\end{figure}

\noindent (It is a little hard to imagine without these figures.)

\begin{figure}[h]
\begin{center}
   \epsfxsize=3 in
  \epsfbox{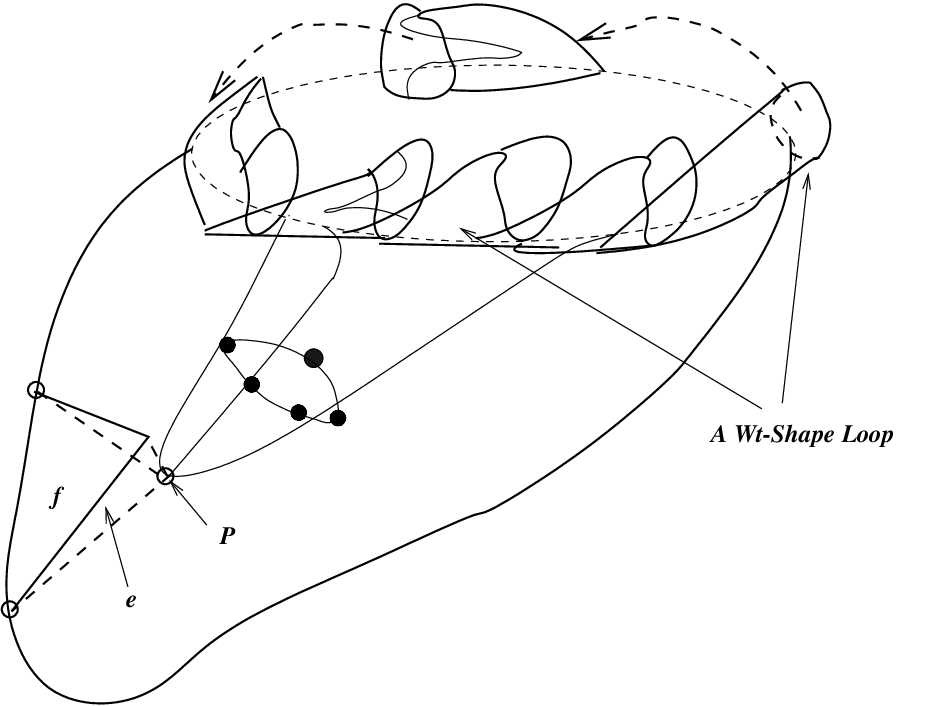}   
\caption{Looped $Wt$-shapes (Wowotou structures). The $Wt$-shapes can be embedded to a repeated sequence.} \label{fig:ManyWt2}
\end{center}
\end{figure}

So, the sequence of the $Wt$-shapes that could be embedded to a repeated sequence. Fig. \ref{fig:ManyWt2} showed it could be true.
(Again: If such a case is not possible, then we will have an infinitive number of $Wt$-shapes when
  there is no 3-cell in $D$ that has two 2-faces on $B$.  Since $D$ is finite, so it is a contradiction.)

Before we continue the proof, we summarize the result we obtained so far:
Recall from the beginning of Case 1 of the proof, we first assumed a 2-cell $f$, a 0-cell $P$ in a 3-cell $e$ in $D$ where $f$, $P$ are also on $B$. If $e$ contains two 2-cells(2-faces) on $B$, $f$ will be able to connect $P$ directly through another 2-cell containing $P$ in $e$. If $e$ does not contains two 2-cell (2-faces) on $B$, there will be a sequence (path) of 2-cells on $B$ connect $f$ and $f(P)$ where $f(P)$ containing $P$ is a 2-cell on $B$. This path can be 1-connected (each adjacent pair of two 2-cells shares 1-cell in the path).
If this path does not associating\/intersecting a 3-cell $e'$ that contains two 2-faces on $B$, (meaning that there is no 3-cell $e'$ intersecting this path or its neighbors), then we will find a $X$ on $B$ that is far-most to $P$ in terms of 3-cell distance (or 1-cell distance). And the local structure of $D$ centered at $X$ is a $Wt$-shape. In the proof above, the local far-most would be fine too.  Note that, not every path from $f$ and $f(P)$ will pass a local far-most point on $B$, but we can search all points on $B$ to find a point $X$ that has the far-most 3-cell distance (or 1-cell distance) to $P$ on $B$ and make a 2-cell path to $P$. (``local'' means a deformed disk or a neighborhood of a point on $B$, especially $Star_B(X)$ and its surrounding points on $B$. Could be for entire $B$.)
In a simple 2-cell path link $X$ to $P$ on $B$, we can determine a $Y$ on $B$ adjacent to $X$ (on $B$) that has a smaller 3-cell distance (or 1-cell distance) to $P$. Otherwise we use $Y$ to replace $X$ in the same 2-cell path. (This is because $Y$ has the same distance to $P$ as $X$ to $P$ since $X$ is already locally maximum). So 1-cell $XY$ included in at least one 2-cell $f_X$ in the path on $B$. The third 0-cell of $f_X$ is denoted as $R_1$. 1-cell $XR_1$ will be included in another 2-cell on $B$ since $B$ is a 2-cycle ($B$ is closed without a boundary curve) except 2-cell $XYR_1$. This new 2-cell will include a new 0-cell we denoted it as $R_0$.  Both $R_0$ and $R_1$ to $P$ will have smaller or equal 3-cell distance (or 1-cell distance) comparing to the 3-cell distance ((or 1-cell distance)) from $X$ to $P$.   We got all conditions matched again to the beginning of the proof of Case 1. We can determine a $Wt$-shape around $X$, or we can find a 3-cell containing two 2-cells on $B$.
We can easily get a shortest simple 2-cell path from $XYR_0$ to $f$. So we can form a 2-cell path on $B$ $\pi = \{f,\cdots,f_X,\cdots,f_P\}$ where $\{f,\cdots,f_X,\}$ and $\{f_X,\cdots,f_P\}$ are simple paths. We can also assume that $\{f_X,\cdots,f_P\}$ is a shortest 2-cell path on $B$. But we have not used these properties.

$f(X)=f_X$ on $B$ will be included in a 3-cell $e_X$ in $D$. $e_X$ has a $P_{f_X}$ on $B$ not in $f(X)$.  As assumed, $e_X$ does not have two 3-cells on $B$, then we will repeatedly find $X_1$ that is (locally far-most to $P_{f_X}$ on a simple 2-cell path from  $f(X)$ to
$P_{f_X}$. Or we first find the locally far-most point $X_1$ to $P_{f_X}$ on $B$ in terms of 3-cell (or 1-cell) distance just like we did above. We can find another $Wt$-shape centered at $X_1$ or find a 3-cell containing two 2-cells on $B$. If we continue this process, we will find $X_1, X_2, \cdots, X_n, \cdots,...$. However, $n$ cannot be infinitive since $D$ has finite number of cells.

In order to simplify the discussion, we can assume that we will get a set of 2-cell paths (sequences) on $B$:  $\pi_0 (f,P) = \{f,\cdots,f_1(X),\cdots,f(P)\}$ ,   $\pi_1 (f_1,P_1) = \{f_1,\cdots,f_2(X_1),\cdots,f(P_1)\}$, $\cdots\cdots$,
 $\pi_n (f_n,P_n) = \{f_n,\cdots,f_{n+1}(X_{n}),\cdots,f(P_n)\}$, $\cdots\cdots$ . $e_i=f_i\cup P_i$ forms a 3-cell in $D$. Note that $\pi_i (f_i,P_i)$ might not be a simple path but $\{f_i,\cdots,f_{i+1}(X_i)\}$ and $\{f_{i+1}(X_i),\cdots,f(P_i)\}$
 are simple paths. $f_{i+1}(X_i)=f_{i+1}$ indicates that  $f_{i+1}$ contains a 0-cell $(X_i)$ on $B$ that 0-cell $X_i$ is the far-most to $P_i$ in terms of 3-cell distance (or 1-cell distance) in $D$ where $e_{i+1}$ contains both $X_i\in f_{i+1}=f_{i+1}(X_i)$ and $P_{i+1}=P_{f_{i+1}}$.

As we showed in Fig. \ref{fig:ManyWt1} and Fig. \ref{fig:ManyWt2}, it is possible that $f_i = f_j, i\le j$.
Even thought this figure could be true, we can still prove our statement: There must exist a 3-cell in $D$ containing two 2-faces on $B$ as follows:


Assume such a case is possible in  Fig. \ref{fig:ManyWt2}, $f_i=f_j, i \le j $. (There could be a possibility that one repeated itself.)  Since $D$ is simply connected,
the entire finite sequence will be forming a collar-like structure  based on the $Wt$-shape containing $X$ as drawn in Fig. \ref{fig:ManyWt2}. For any $<f_k, P_k>$ pair in the sequence of $Wt$-shapes, where $\pi_k$ denotes an 2-cell path from $f_k$ to $P_k$ (meaning to a 2-cell containing $P_k$). $\pi_k$ contains $f_{k+1}=f_{k+1}(X_k)$.

Let us look at the details of a $Wt$-shape for any $k, i\le k\le j$.

We wished we could obtain a statement: For any 2-cell path $\pi = \{f,\cdots,f_X,\cdots,f_P\}$ on $B$, if the far-most point $X$ in on the path is not a local far-most point locally on $B$, $\pi$ must intersect with  a 3-cell that contains two 2-cells on $B$. This is not true since Fig. \ref{fig:ManyWt2}. But intuitively, every ``collar'' of a $Wt$-shape will be wrapped with $Wt$-shapes. This will generate an infinitive growth to preserve simple connectedness.   Therefore, we want to find a local far-most point to another 2-cell in a $Wt$-shape except from the original $X_k$.

\begin{figure}[h]
\begin{center}
   \epsfxsize=3 in
   \epsfbox{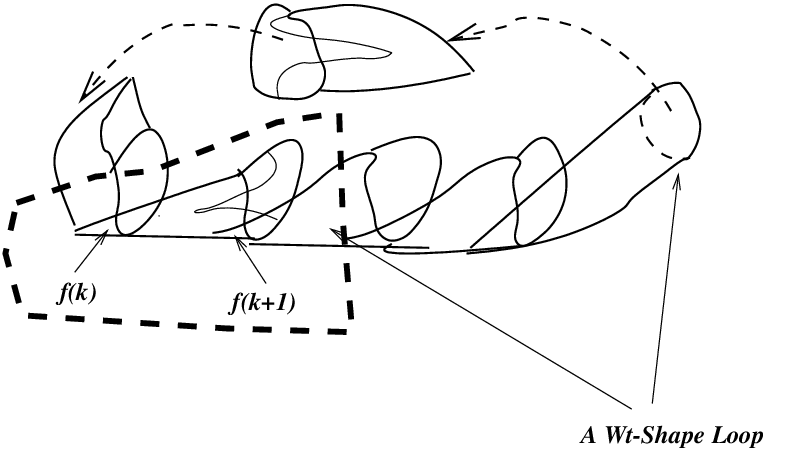}   
\caption{ $Wt$-shape for index $k$.} \label{fig:Wt8_1}
\end{center}
\end{figure}

\begin{figure}[h]
\begin{center}
   \epsfxsize=2.5 in
   \epsfbox{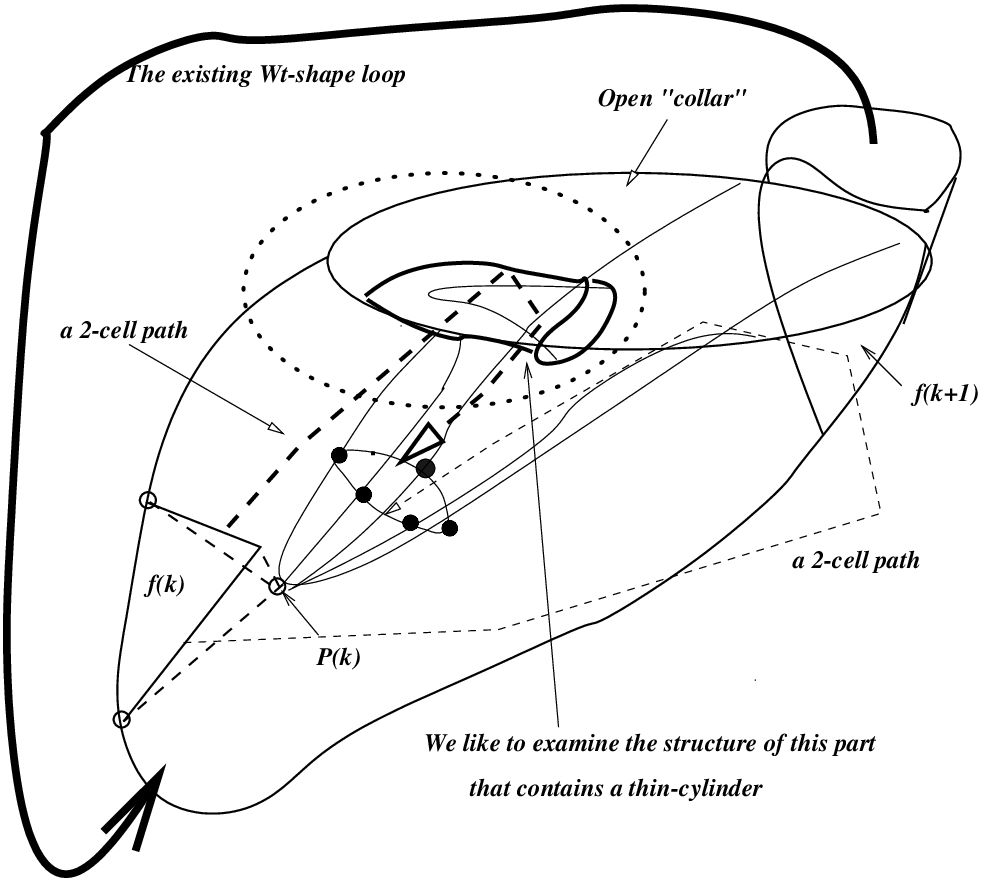}   
\caption{A part of a new $Wt$-shape on the base of $Wt_k$.} \label{fig:NewOnWt_k_0}
\end{center}
\end{figure}

Now we can have an interesting and intuitive observation first:
For $f(X_{k+1})=f_{k+1}$,  we want to (we can) examine a 2-cell path $\pi'$ that is not containing $f_{k+1}$ from $f_k$ to $P_k$.  This is because we can always make a path in $Star_B{f_{k+1}}$ not passing $f_{k+1}$.
Even though this path is a bit long in possibilities, we want find another $Wt$-shape, this new shape will grow also if the property we are looking for is not true (i.e. there is a 3-cell containing two 2-cells on $B$).
See Fig.  \ref{fig:Wt8_1}.
If growing, the sequence must be the collar-like
structure in the $Wt$-shape based on $f_k$ (not the old $f$), such a structure will be recursively run forever. It is contradiction to that $D$ only contains finite number cells.  So we proved our statement.
The problem is that not every $\pi$ path from $f_k$ to $P_k$ will determine a $Wt$-shape due to a local far-most point $X$ (or $X_i$)is needed in the path. However, similar to the proof of finding the end point $X$ or $X_i$ of a $Wt$-shape above, we can see that it will determine a part of a new $Wt$-shape  (even can not determine the place of an ending point $X$) that wraps---with empty space inside of the wrapping, like a thin-cylinder. See Fig. \ref{fig:NewOnWt_k_0} if there is no 3-cell containing two 2-cells on $B$.

The  ``collar'' based on the $Wt$-shape for
$f_k$ (we can call the $Wt_k$ shape) is open for reaching $P_k$($P(k)$). And it is  still wrapped by a ``thin cylinder'' with thickness $=1$.  Since  the ``thin cylinder''  must have exactly an end point or the bottom part that is closed, otherwise $D$ is not simply connected (no closed ends or two closed ends in $D$ will make D is not simply connected). So it is a new $Wt$-shape. So we proved our case intuitively: We have a new $Wt$-shape that is not anyone in the $Wt$-shape loop we found previously.


In fact, the property of simple connectedness will keep the repeated pattern on the collar. This completes the proof for Case 1.

Due to the fact that Case 2 could contain a case that is Case 1 (Derive a case that is case 1); otherwise, we should prove for Case 2 first. We will see it in the following.

Since Case 2 is a special case of Case 1, we actually do not need the specifically prove for Case 2. However, Case 2 is usually easier than Case 1 by adding a constraint. Therefore, it is interesting to see how Case 2 is locally configured.

{\bf Case 2:}  Case 2 is usually simpler than Case 1 since it has an extra constraint. However, it contains a subcase that could be Case 1, we need to face Case 1 first.  Due to the fact of Case 2 special case of Case 1, it might be not necessary. Anyway, we include the complete proof here.

\begin{figure}[h]
\begin{center}
   \epsfxsize=5.5in
  \epsfbox{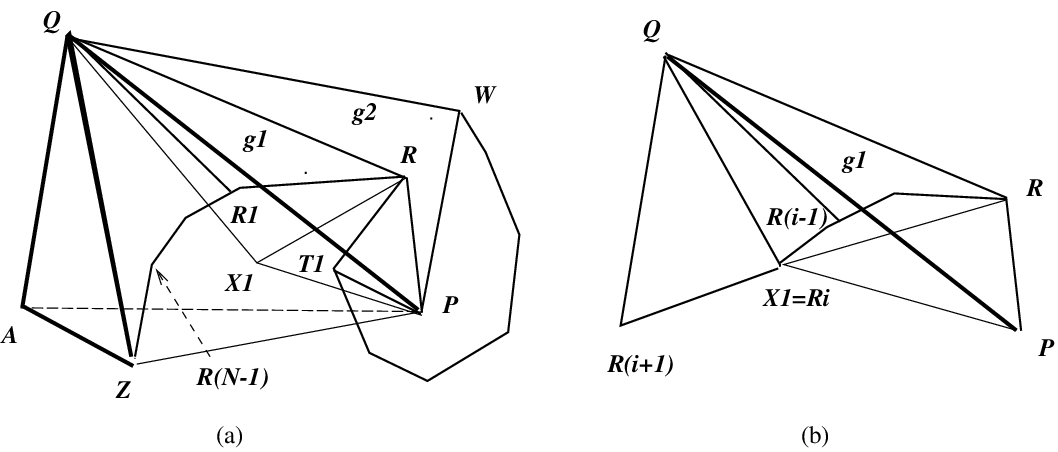}   
\caption{Case 2: (a) General,       (b) Detailed.  } \label{fig:old21}
\end{center}
\end{figure}

For Case 2, a $f$ of $e$ is on $B$ and a 1-cell not on $f$ is also on $B$.  Assume this 1-cell is $[Q,P]$ where $Q\in f$ and $P\notin f$. Since  $[Q,P]$  is on $B$, there are two 2-cells on $B$ (not $f$) sharing $[Q,P]$.
(each 1-cell will be included in exact two 2-cells $g_1$ and $g_2$ if $B$ is a 2-cycle.)  $Star_{M}(P)$ is homeomorphic to a 3-disk (a closed 3-ball);  $Star_{M}(P)$ contains 3-cell $e$, $g_1$, and $g_2$.
$Q$ is the point (0-cell); $Star_B(Q)$ is a 2-disk containing $f$, $P$, $g_1$, and $2_2$.  From $g_1$ to $f$,  there must be a 2-cell path where each adjacent pair of 2-cells shares
a 1-cell (called 1-connected~\cite{Chen15}) with the end-point $Q$. See Fig. \ref {fig:old21}
From $g_2$ to $f$,  there is another 2-cell path where each adjacent pair shares a 1-cell with the end-point $Q$ too.  Assume $P$, $Q$, $A$, $Z$ are four vertices in $e$, $\Delta(P,Z,Q)$ is 2-cell of $e$ in $D$. If  $\Delta(P,Z,Q)$ is $g_1$ (or $g_2$) , then we have done our proof (for Case 2) since $g_1$ is on $B$ and $g_1$ is in $e$.

If  $\Delta(P,Z,Q)$ is not
$g_1$ (and not $g_2$), assume $R$ is the third 0-cell of  $g_1$ , there must be a 2-cell path in $Star_B(Q)$ shares a 1-cell $[QR]$ with $g_1$, and  shares a 1-cell $[QB]$ with $f$.  This path is 1-connected.
The path is a set of triangles with common point $Q$. The other points can be named as $P, R=R_0, R_1, \cdots, R_N=B$ on $Link_B(Q)$.
$N$ is finite since $M$ contains only a finite number of cells.
$g_1=\Delta(P,Q,R)$ on $B$ must be in a 3-cell ($\Delta_3$) in $D$ (each 2-cell must be contained by two 3-cells in $M$; one is in $D$). If $R_1$ is in $\Delta_3$, we have two 2-faces $\Delta(QRR_1)$ and $g_1$ of a 3-cell ($\Delta_3$) that are on $B$. We did our proof.

If $R_1$ is not in $\Delta_3$, there is a point $X_1$ in a 3-cell $\Delta_3=\Delta_3(P,Q,R,X_1)$. $X_1$ is not $R_1$ as we talked above. There are four possibilities (subcases):

(i)Since $X_1$ in 3-cell that has one 2-face $QRP$ on $B$, $X_1$ must be on $B$ as assumed. (Otherwise, we just can remove this 3-cell.) $X_1$ is on $Link_B(Q)$ (Most like case).

(ii)$X_1$ is on $Link_B(P)$.  Let us assume the $[RP]$ is contained by 2-cell $g_1$ and another 2-cell $PRT_1$. So we want to examine the case $X_1$ is in between $T_1$ to $R_1$ on $Link_B(R)$.

(iii)$X_1$ is on $Link_B(R)$. The only possibility is that $T_1=X_1$ since the 2-cells $X_1PR$ and $T_1PR$ could not partially overlapping except just sharing 1-cell $PR$.  In addition, the 2-cell $T_1PR$ could not be included in 3-cell $R_1X_1QPR$. When $T_1=X_1$, we have two 2-faces of the 3-cell $R_1X_1QPR$ are on $B$. So we did it.

(iv)$X_1$ has its own $Star_B(X_1)$ do not contain $Q$ or $P$ or $R$ (very unusual).

If $X_1$ is not one of four subcases, then $X_1$ is not on $B$. Thus, $X_1\in  D-B$,  $g_1$=$QRP$ on $B$ can be replaced  with the inner point $X_1$ in the same 3-cell.

In fact, Sub-Case (i) and Sub-Case (ii) are similar. We now prove for Case (i).

For (i), if $X_1$ is $R_2$ on $B$, since $QRX_1$ is a triangle in $D$, $X_1=R_2$ makes $RX_1R_1$ to be a triangle. So  $QRR_1R_2$ is a 3-cell  So it contains two 2-faces on $B$.
We done.

So assume that $X_1$ is not $R_2$.  Could we have $X_1=R_3$? or $X_1=R_i$ where  $2<i<N$ ?

(i.a) If $X_1=R_N=Z$ , then $QZPR$ is a 3-cell triangle. $g_1$ and $QZ$ will satisfy Case 2.
 The arc (path) of $Z R_{N-1}\cdots R_1RP$ will be shortened to $Z R_{N-1}\cdots R_1R$. We can find $Y_1$ such that $QBR_{N-1}Y_1$ is a 3-cell. Repeat this process: we have (a) $Y_1$ is in $D-B$, or (b)
 we either find two 2-faces of an $e_i$ on the $B$ or the arc in $Link_B(Q)$ is shortened. Until we get to only  2 or 1 point in the arc between new ''Z'' and ''P'', we proved the statement.

(i.b) If $X_1=R_i$
where  $2<i<N$, So $\sigma_1$ ($QRP$) and $QR_{i}$ are in a 3-cell $QRPR_i$ that is Case 2.
 The arc length is much shortened (at least by 1 strictly). Since $QR_{i-1}R_i$ is a triangle on $B$ as assumed, We find $Y_1$ such that $QR_{i-1}R_iY_1$ is a 3-cell including $QR_{i-1}R_i$.
 As we repeat this process, we will get a 3-cell having two 2-faces on $B$ (since $N$ is finite). If $Y_1$ is not on $B$, we done.
 If $Y_1$ is on $B$, it must be one of above four cases. We repeat the above process (from other direction) if it is the case (i) (or (ii)). Until this arc length is reduced to 1 or 2. We done our proof for this case. If it is one of other three cases. We will present the solutions for each later.

\begin{figure}[h]
\begin{center}
   \epsfxsize=5.5in
  \epsfbox{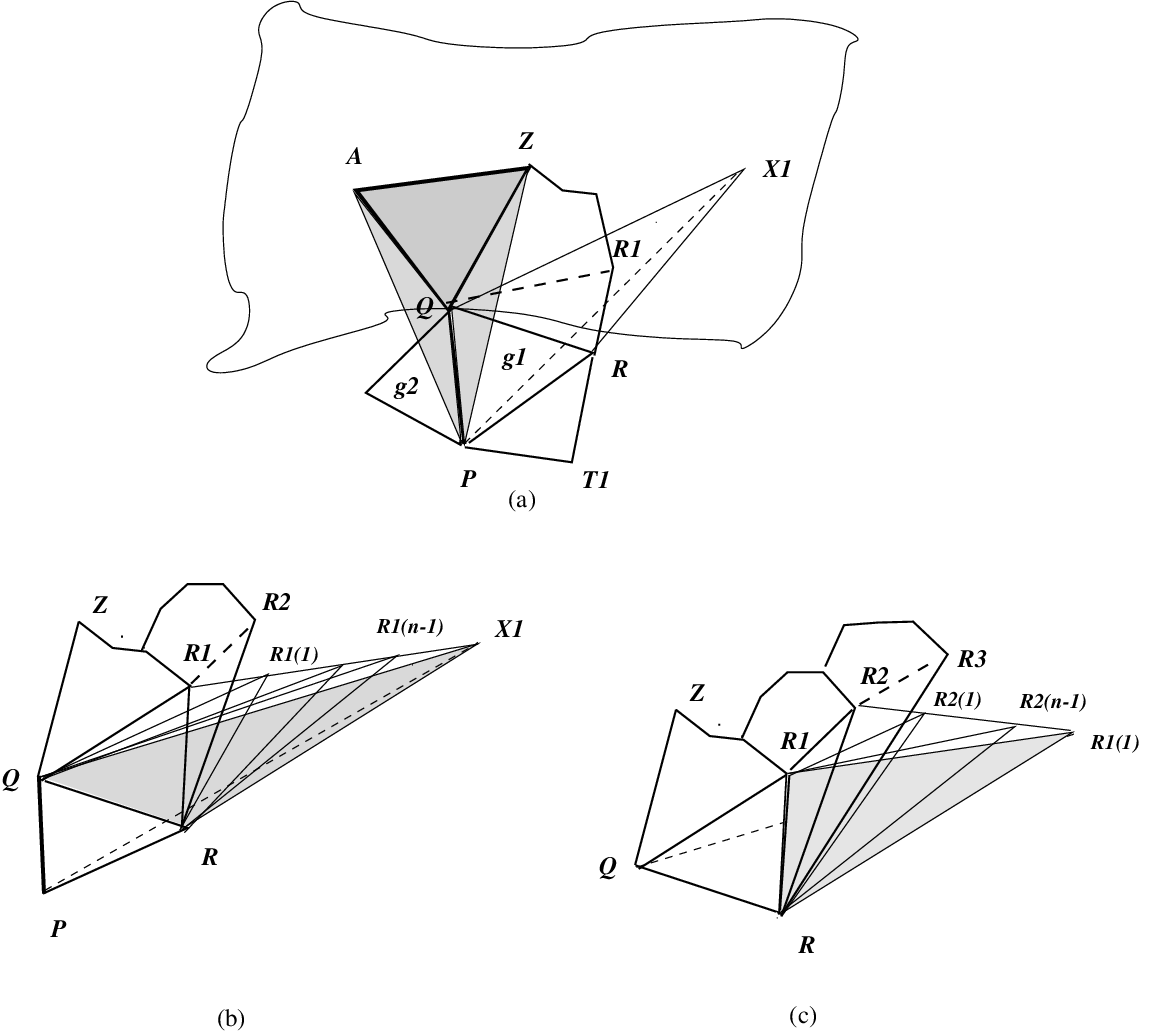}   
\caption{ Special Subcases  }\label{fig:old22}
\end{center}
\end{figure}

For (ii), (ii)$X_1$ is on $Link_B(P)$. 1-cell $RX_1$ is in $Star_B(R )$. Still a case 2. The arc from $P$ to $R_1$ in $Link_B(P)$ cycle contains $X_1$. Assume $T_1$ is the third 0-cell of the 2-cell in $Start_B(P)$ containing 1-cell $PR$. Then if $X_1=T_1$, the 3-cell $QRPT_1$ contains two 2-faces on $B$. Done. Now, if the $X_1=T_i$, $1<i<n$, is in the middle of an arc $\pi(T_1,R_1)$ that is assumed as $T_1,T_2,\cdots,T_n=R_1$. So we got exact the case like Sub-Case (i). $RX_1$ is on $B$ instead of $QX_1$ in the previous case.  It was proved above.  See Fig. \ref {fig:old22} .

Since Sub-Case (iii) was already discussed, we now only need to explore the case (iv).

For (iv), $X_1$ is not on  $Link_B(Q)$,  $Link_B(R)$, and  $Link_B(P)$. This is a special case of Case 1. If we can prove Case 1 without using Case 2, we get a proof. We already did it above.

So, we complete the proof for Case 2.

In other words, there is 1-cell path on $B$ that is shortest from $R_1$ to $X_1$ must passes $R_1^{(1)}$. However, there will be another $R_1^{{(1)}^{(1)}}$, etc. unless it is on $B$. This means if we
have Case 2, there will be a 3-cell near by that contains two 2-faces On $B$.

We have complete the proof for Case (iv).

To summarize, we will see we got a shortened list (arc) of $R_i$'s, using mathematical induction, $X_1$ or $Y_1$ cannot be any $R_i$ otherwise, there will be an $e_i$ including two 2-faces on $B$.
So we proved our statement. If $X_1$ or $Y_1$ in $D-B$, then it makes a contradiction (of the assumption) for the general condition of there is no inner point of $D$ in a 3-cell $e_i$ that contains a 2-face on $B$.

Now, we have completed the proof for this theorem especially for the two special cases in The main theorem in Section 4.  When we have an case $m>3$, we just need to find two special $(m-1)$-faces $f_k$ and $f_{k+1}$
to find $Wt$-shape in higher dimension or to determine an $m$-cell in $D$ that has two $(m-1)$-faces directly on $B$. The purpose of finding $Wt$-shape is also for the purpose to mathematically prove that such an $m$-cell exists.     $\diamond $

For 2D case, it more easy to see. We can also use 4 vertices theorem to find the biggest curvature point to do the contraction. Always contract the cell that has the biggest curvature point. In 3-manifold, the principle is
the same. The further most point is the one that has the biggest local curvature. We can always contract a cell that contains the point. If there is no such a point, there will be a handle structure, it will make $B$ is not
simply connected.

\newpage

\section{Appendix B: The Proving Path of the Paper and Some Remarks for Readers}

In this appendix, we provide some explanations of our attempt and ideas. In order to keep it simple and intuitive, the materials here are less formal.

In this paper, we make the connection between two famous problems: The general Jordan separation problem and the Poincare conjecture for 3D manifolds.
Even thought most mathematicians believed that the Poincare conjecture was solved, but there are still some top geometers who are still questioning this. \footnote{See Youtube videos for more information.}
It is fair to say that at least for some steps in the current proof using geometric analysis methods  still contain some part that are not totally constructive. Plus, such a proof was still very long and
not elementary. Therefore, a topological and constructive proof of the Poincare conjecture is remaining very important and expected.

This paper algorithmically proved that if the general Jordan separation property holds for a 3D compact manifold, then the Poincare conjecture is true.
In fact, Chen and Krantz proved the general Jordan separation property for all dimensions for discrete manifolds~\cite{Chen-Krantz}. The discrete manifold is a generalization of
the piecewise-linear manifold.

In order to make easier to understand this paper, we first list several related concepts and results:

{\bf The Jordan Curve Theorem:}  A simply closed curve on a simply connected plane will separate the plane into two components. Each component is (path-) connected. The plane can be replaced a closed 2-manifold, the
theorem is still valid.

{\bf The Jordan-Schoenflies Theorem:} Each component of above the above theorem is homeomorphic to an open 2-ball (or 2-disk) for a closed 2-manifold. The simple closed curve is the common boundary of the two components.

{\bf The General Jordan-Schoenflies Theorem:} Embedding an $(n-1)$ sphere $S^{(n-1)}$ locally flatly in an $n$ sphere $S^{n}$, then it decomposes $S^{n}$ into two components. In addition, the embedded $S^{(n-1)}$ is the common boundary of the two components and each component is homeomorphic to the $n$-ball.

{\bf The Poincare Conjecture:} Every simply connected compact smooth closed 3D manifold $(M_3)$ is homeomorphic to $S^3$.

{\bf The General Jordan Separation Property:} Every $(n-1)$-cycle on a simply connected compact smooth closed $n$-manifold will separate this $n$-manifold into two components.

So we can see that if the Poincare Conjecture holds with additional to the General Jordan-Schoenflies Theorem, then we can first make the homeomorphic mapping from $M_3$ to $S^3$, then use the General Jordan-Schoenflies Theorem,
we can obtain The General Jordan Separation Property for 3D manifolds. Of course, we might need to add the condition of local flatness.

We can see that the above five theorems, conjecture, an properties are very closed related.

It is very natural to ask whether or not we can directly prove the General Jordan Separation Property, and use it to prove The Poincare Conjecture.

The problem is that the proofs of The General Jordan-Schoenflies Theorem contained some gaps. They need the Poincare Conjecture's assertion to validate the proof.  Therefore, the Poincare Conjecture
can be used to prove The General Jordan Separation Property.

This paper attempts to use The General Jordan Separation Property to prove the Poincare Conjecture. The proving path is the following:

(1) We need to first prove the each component is simply connected. This is important to to the following proof since without simple connectedness, we could not guarantee the contraction is valid.
(2) We need to prove each component is homeomorphic to a 3-cell. In order to prove this, we designed a special way to contract a special 3D manifold to be a 3-cell. After we did this, we realized that this idea was
used by Whitehead in his false proof of the Poincare Conjecture. It was also possible that we some how heard this idea before but trigged later.  Each time of contraction, we need to keep the boundary of the manifold
to be 2-cycle.

One much see that we must only have a finite number of cells to be contracted.  $M_3$ must only contain finite cells. Fortunately, it is well known that for a compact smooth manifold there exist finite triangulations.

\end{document}